\documentclass[11pt]{amsart}

\usepackage{amsmath}
\usepackage{amsthm}
\usepackage{amsfonts}
\usepackage{amssymb}
\usepackage[dvipsnames]{xcolor}
\definecolor{darkgreen}{rgb}{0,0.45,0}
\definecolor{darkred}{rgb}{0.75,0,0}
\definecolor{darkblue}{rgb}{0,0,0.6}
\usepackage[colorlinks,citecolor=darkgreen,linkcolor=darkblue,urlcolor=darkred]{hyperref}
\usepackage{tikz}
\usepackage{tikz-cd}
\usepackage{stmaryrd}
\usepackage{float}
\usepackage{bbm}

\usepackage{setspace}
\usepackage{ arydshln}

\newcommand{\iso}{\cong}
\renewcommand{\to}{\rightarrow}

\newcommand{\Q}{\mathbb{Q}}
\newcommand{\C}{\mathbb{C}}
\newcommand{\R}{\mathbb{R}}

\newcommand{\Z}{\mathbb{Z}}

\newcommand{\A}{\mathbb{A}}

\newcommand{\p}{\mathfrak{p}}
\newcommand{\g}{\mathfrak{g}}

\newcommand{\M}{\mathcal{M}}

\renewcommand{\H}{\mathcal{H}}
\renewcommand{\O}{\mathcal{O}}
\newcommand{\E}{\mathcal{E}}

\newcommand{\rank}{\operatorname{rank}}
\newcommand{\End}{\operatorname{End}}
\newcommand{\Res}{\operatorname{Res}}
\newcommand{\Hom}{\operatorname{Hom}}
\newcommand{\Gal}{\operatorname{Gal}}
\newcommand{\im}{\operatorname{im}}

\newcommand{\an}{\mathrm{an}}

\newcommand{\tors}{\mathrm{tors}}

\newcommand{\SL}{\operatorname{SL}}
\newcommand{\GL}{\operatorname{GL}}
\newcommand{\Frob}{\operatorname{Frob}}
\newcommand{\Tr}{\operatorname{Tr}}

\newcommand{\sgn}{\operatorname{sgn}}

\newcommand{\Ad}{\operatorname{Ad}}
\newcommand{\Ind}{\operatorname{Ind}}

\theoremstyle{plain}
\newtheorem{theorem}{Theorem}[section]

\newtheorem{proposition}[theorem]{Proposition}
\newtheorem{lemma}[theorem]{Lemma}
\newtheorem{corollary}[theorem]{Corollary}

\newtheorem{conjecture}[theorem]{Conjecture}

\theoremstyle{definition}
\newtheorem{definition}[theorem]{Definition}
\newtheorem*{definition*}{Definition}
\newtheorem{example}[theorem]{Example}

\newtheorem{remark}[theorem]{Remark}

\newtheorem{algorithm}[theorem]{Algorithm}

\numberwithin{equation}{section}
\numberwithin{table}{section}

\begin{document}

\title{Motivic action on coherent cohomology of Hilbert modular varieties}
\author{Aleksander Horawa}
\date{\today}
\maketitle

\begin{abstract}
	We propose an action of a certain motivic cohomology group on the coherent cohomology of Hilbert modular varieties, extending conjectures of Venkatesh, Prasanna, and Harris. The action is described in two ways: on cohomology modulo $p$ and over $\C$, and we conjecture that they both lift to an action on cohomology with integral coefficients. The conjecture is supported by theoretical evidence based on Stark's conjecture on special values of Artin $L$-functions, and by numerical evidence in base change cases.
\end{abstract}

\section{Introduction}

A surprising property of the cohomology of locally symmetric spaces is that Hecke operators can act on multiple cohomological degrees with the same eigenvalues. One can observe this by a standard dimension count, but this does little to explain the phenomenon. In a series of papers, Venkatesh and his collaborators propose an arithmetic reason for this: a hidden degree-shifting action of a certain motivic cohomology group.

Initially, Prasanna--Venkatesh \cite{Prasanna_Venkatesh} and Venkatesh \cite{Venkatesh:Derived_Hecke} developed these conjectures for singular cohomology of locally symmetric spaces. Later, Harris--Venkatesh~\cite{Harris_Venkatesh} observed similar behavior for coherent cohomology of the Hodge bundle on the modular curve. See also \cite{Marcil, Darmon_Harris_Rotger_Venkatesh} for more evidence for their conjecture. Connections to derived Galois deformation theory and modularity lifting were also explored by Galatius--Venkatesh~\cite{Galatius_Venkatesh}. For a general introduction to this subject, see~\cite{Venkatesh:Takagi, Venkatesh:ICM}.

In this paper, we propose analogous conjectures for coherent cohomology of the Hodge bundle on Hilbert modular varieties. To give a more precise statement, we first set up some notation.

Let $F$ be a totally real extension of $\Q$ of degree $d$ and let $f$ be a parallel weight one, cuspidal, normalized Hilbert modular eigenform for $F$, with Fourier coefficients in the ring of integers $\O_{E}$ of a number field~${E}$. One can identify $f$ with a section of the Hodge bundle $\omega$ on a Hilbert modular variety~$X$:
$$f \in H^0(X, \omega) \otimes \O_{E}.$$
More specifically, we consider an integral model $X$ of the toroidal compactification of the open Hilbert modular variety with good reduction away from primes dividing the discriminant of~$F$ and the conductor of $f$. While this choice is not canonical, the resulting cohomology groups are independent of the choice of~$X$.

The action of the Hecke algebra extends to higher cohomology groups $H^i(X, \omega) \otimes \O_{E}$ and we may consider the subspace on which the Hecke algebra acts with the same eigenvalues as on $f$, which we denote by $H^i(X, \omega)_f$. It follows from~\cite{Su} that
\begin{equation}\label{eqn:dim_coh}
	\rank H^i(X, \omega)_f = \binom{d}{i}
\end{equation}
(c.f.\ Corollary~\ref{cor:dim_of_coh_group}). There is a motivic cohomology group $U_f$ associated to $f$, which is an $\O_E$-module of rank $d = [F:\Q]$ (Corollary~\ref{cor:Stark_basis_HMF}); explicitly, it is the {\em Stark unit group}~\cite{StarkII} for the trace zero adjoint representation of $f$. We conjecture that there is a degree-shifting action of its dual $U_f^\vee$ on the cohomology space $H^\ast(X, \omega)_f$ which makes $H^\ast(X, \omega)_f$ a~module of rank one over the exterior algebra $\bigwedge^\ast U_f^\vee$, generated by $f \in H^0(X, \omega)_f$.

We can describe this action in two ways: modulo $p$ and over $\C$. Let $\mathfrak p$ be a prime of $\O_{E}$, $n \geq 1$ be an integer, and $\iota \colon E \hookrightarrow \C$ be an embedding. We show that there is:
\begin{enumerate}
	\item a map 
	$$\bigoplus\limits_{j=1}^d  U_{f, j}^{\mathfrak p^n} \to \displaystyle U_f^\vee \otimes \O_{E}/\p^n $$ 
	for some free $\O_{E}/\p^n$-modules $U_{f, j}^{\p^n}$ of rank one (Proposition~\ref{prop:dual_units_mod_p}), and define an action of $U_{f, j}^{\mathfrak p^n}$ on $H^\ast(X, \omega)_f \otimes \O_{E}/\p^n$ by derived Hecke operators (Definition \ref{def:action_mod_p}),\label{intro:actionp}
	\vspace{0.3cm}
	\item an isomorphism 
	$$\displaystyle \bigoplus\limits_{j=1}^d U_{f, j}^{\C}   \overset \iso\to U_f^\vee \otimes \C$$ 
	for some one-dimensional $\C$-vector spaces $U_{f, j}^{\C}$ (Proposition~\ref{prop:dual_units_over_C}), and define an action of $U_{f, j}^\C$ on $H^\ast(X, \omega)_f \otimes \C$  by partial complex conjugation $z_j \mapsto \overline{z_j}$ (Definition~\ref{def:action}). \label{intro:actionC}
\end{enumerate}

The following conjecture predicts that these actions come from a single ``motivic'' action that is defined rationally or even integrally.

\begin{conjecture}[Conjectures~\ref{conj:modulo}, \ref{conj:motivic_action}]\label{conj:A}
	There is a graded action $\star$ of the exterior algebra $\bigwedge^\ast U_f^\vee$ on $H^\ast(X, \omega)_f$ such that:
	\begin{enumerate}
		\item the action of $\bigwedge^\ast U_f^\vee \otimes \O_{E}/\p^n$ is the same as that in (\ref{intro:actionp}) above, up to $\GL_n(\O_{E})$ ambiguity,
		\item the action of $\bigwedge^\ast U_f^\vee \otimes \C$ is the same as that in (\ref{intro:actionC}) above, up to $\GL_n(E)$ ambiguity.
	\end{enumerate}
	Moreover, $H^\ast(X, \omega)_f$ is generated by $f \in H^0(X, \omega)_f$ over $ \bigwedge^\ast U_f^\vee$.
\end{conjecture}

The conjectures will be stated precisely in the main body of the paper.

Part (1) is a generalization of the main conjecture of Harris and Venkatesh~\cite[Conjecture~3.1]{Harris_Venkatesh}. It should be seen as a first step towards establishing a $p$-adic conjecture, similar to Venkatesh's conjecture~\cite{Venkatesh:Derived_Hecke}. In fact, our original motivation to study the Stark unit group $U_f$ for Hilbert modular forms was to generalize the conjecture of Darmon--Lauder--Rotger~\cite{Darmon_Lauder_Rotger} to elliptic curves over totally real fields. A $p$-adic version of Conjecture~\ref{conj:A} may explain the appearance of $p$-adic logarithms of Stark units therein.

Part (2) is similar to the main conjecture of Prasanna and Venkatesh~\cite[Conjecture~1.2.1]{Prasanna_Venkatesh} but in the coherent (as opposed to singular) cohomology setting. We discuss the precise relationship in Appendix~\ref{app:Prasanna-Venkatesh}. As far as we know, it is new even when $F = \Q$. In the Hilbert case, it is also closely related to the study of period invariants attached to Hilbert modular forms at the infinite places. Such period invariants had previously been defined by Shimura \cite{Shimura_HMF, Shimura:on_crit_val}, Harris \cite{Harris:JAMS, Harris_periods_I, Harris_periods_II}, and Ichino--Prasanna \cite{Ichino_Prasanna} in cases where the weight of~$f$ is at least two at some of the infinite places. The parallel weight one case is different because the form does not transfer to a quaternion algebra ramified at any infinite place, so the periods at infinite places do not admit a simple interpretation as periods of a holomorphic differential form on a Shimura curve, or even as ratios of periods of holomorphic forms on quaternionic Shimura varieties. Instead, we give specific linear combinations of the higher coherent cohomology classes which we expect to be rational in coherent cohomology. The expressions involve logarithms of units which is natural because the adjoint $L$-value is non-critical at $s=1$ in this case, so one should expect the periods to be of ``Beilinson-type".

Part (2) of the conjecture admits a natural generalization to partial weight one Hilbert modular forms, which we discuss in Appendix~\ref{app:Prasanna-Venkatesh}. In that case, however, the motivic cohomology group in question does not admit an interpretation as a unit group. 

These conjectures lead to many interesting questions about potential generalizations to other reductive groups which we are currently pursing elsewhere. We were also recently made aware of the forthcoming work of Gyujin Oh and Stanislav Ivanov Atanasov on this topic.

Next, we give a more explicit versions of Conjecture~\ref{conj:A} in the cases $[F:\Q] = 1$ and $[F:\Q] = 2$ and summarize our evidence for them. For simplicity, we assume that the automorphic representation associated to $f$ is not supercuspidal at $p = 2$ (this assumption avoids a potential factor of $\sqrt{2}$ and we expect it to be unnecessary; see Remark~\ref{rmk:adjoint_conductor}).

{\bf The case $[F:\Q] = 1$: modular curves.}  When $[F:\Q] = 1$,  $X$ is a modular curve and $f$ is a classical modular form of weight one. This is the situation considered by Harris--Venkatesh~\cite{Harris_Venkatesh} and Conjecture~\ref{conj:A}~(1) specializes to their conjecture. Conjecture~\ref{conj:A}~(2) is its archimedean version and follows from Stark's conjecture on special values of Artin $L$-functions.

\begin{theorem}[Corollary~\ref{cor:conj_for_F=Q}]\label{thmB}
	Let $f$ be a modular form of weight one. If $f$ does not have CM or the Fourier coefficients of $f$ are not rational, assume Stark's conjecture~\ref{conj:Stark}. Then Conjecture~\ref{conj:A}~(2) is true and has the following explicit form: there is an action $\star$ of $\bigwedge^\ast U_f^\vee \otimes E$ on $H^\ast(X, \omega)_f$ such that given $u_f^\vee \in U_f^\vee$, the action:
	$$H^0(X, \omega)_f \overset {u_f^\vee \star} \longrightarrow H^1(X, \omega)_f$$
	is given by
	$$f \mapsto \frac{\omega_f^{\infty}}{\log |u_{f}|},$$
	where
	$$\omega_f^\infty = f(- \overline z) y \frac{dz \wedge d\overline z}{y^2} \in H^1(X_\C, \omega)_f$$
	and $u_f \in U_L$ is a unit in the splitting field $L$ of the adjoint Artin representation of $f$ associated to $u_f^\vee$.
\end{theorem}

In fact, the rationality of $\frac{\omega_f^\infty}{\log|u_f|}$ is equivalent to Stark's conjecture for the trace 0 adjoint representation of~$f$.

{\bf The case $[F:\Q] = 2$: Hilbert modular surfaces.} When $[F:\Q] = 2$, $X$ is a Hilbert modular surface and $f$ is a Hilbert modular form in two variables $z_1, z_2$. We give an explication of Conjecture~\ref{conj:A}~(2) in this case and summarize our evidence for it.

Corollary~\ref{cor:dim_of_coh_group} gives an explicit basis for $H^\ast(X, \omega)_f \otimes \C$:
\begin{align*}
	f & \in H^0(X, \omega)_f \\
	\omega_f^{\sigma_1}, \omega_f^{\sigma_2} & \in H^1(X, \omega)_f \otimes \C \\
	\omega_f^{\sigma_1, \sigma_2} & \in H^2(X, \omega)_f \otimes \C,
\end{align*}
where we choose a fundamental unit $\epsilon$ such that $\epsilon_1 < 0$, $\epsilon_2 > 0$ and let:
\begin{align}
	\omega_f^{\sigma_1} & = f(\epsilon_1 \overline{z_1}, \epsilon_2 z_2) y_1 \, \frac{dz_1 \wedge d\overline {z_1}}{y_1^2}, \label{eqn:omega_f^1} \\
	\omega_f^{\sigma_2} & = f(\epsilon_2 z_1, \epsilon_1 \overline{z_2}) y_2 \, \frac{dz_2 \wedge d\overline {z_2}}{y_2^2}, \label{eqn:omega_f^2} \\
	\omega_f^{\sigma_1, \sigma_2} & = f(- \overline{z_1}, - \overline{z_2}) y_1 y_2\, \frac{dz_1 \wedge d\overline {z_1}}{y_1^2} \,  \frac{dz_2 \wedge d\overline {z_2}}{y_2^2}.
\end{align}

Conjecture~\ref{conj:A}~(2) gives explicit linear combinations of these cohomology classes which should be $E$-rational in cohomology. Specifically, there are four units $u_{11}, u_{12}, u_{21}, u_{22} \in U_L$ in the splitting field $L$ of the adjoint Artin representation of $f$, and we can form the Stark regulator matrix:
\begin{equation}
	R_f = \begin{pmatrix}
		\log|\tau (u_{11}) | & \log|\tau(u_{12})| \\
		\log|\tau(u_{21})| & \log|\tau(u_{22}) |
	\end{pmatrix},
\end{equation}
where $\tau \colon L \hookrightarrow \C$ is a complex embedding of $L$. We show that there is an explicit basis $u_1^\vee, u_2^\vee$ of $U_f^\vee \otimes E$ such that the action of $u_1^\vee$ and $u_2^\vee$ is given by:
\begin{align}
	u_1^\vee \star f & = \frac{ \log|\tau(u_{22}) | \cdot \omega_f^{\sigma_1} - \log|\tau (u_{21}) | \cdot \omega_f^{\sigma_2}}{\det R_f} \in H^1(X, \omega)_f \otimes \C, \\
	u_2^\vee \star f & =\frac{- \log|\tau (u_{12}) | \cdot \omega_f^{\sigma_1} +\log|\tau (u_{11}) |\cdot \omega_f^{\sigma_2}}{\det R_f}  \in H^1(X, \omega)_f \otimes \C
\end{align}
and the action of $u_1^\vee \wedge u_2^\vee$ is given by:
\begin{align}
	(u_1^\vee \wedge u_2^\vee) \star f & = \frac{\omega_f^{\sigma_1, \sigma_2}}{\det R_f}  \in H^2(X, \omega)_f \otimes \C.
\end{align}

We then have the following explicit version of Conjecture~\ref{conj:A}~(2) for $[F:\Q] = 2$.

\begin{conjecture}[Conjecture~\ref{conj:Harris_periods}]\label{conj:B}
	\leavevmode
	\begin{enumerate}
		\item[(a)] A basis of $H^1(X, \omega)_f$ is given by:
		\begin{align*}
			\frac{ \log|\tau(u_{22}) | \cdot \omega_f^{\sigma_1} - \log|\tau (u_{21}) | \cdot \omega_f^{\sigma_2}}{\det R_f}, \\
			\frac{- \log|\tau (u_{12}) | \cdot \omega_f^{\sigma_1} +\log|\tau (u_{11}) |\cdot \omega_f^{\sigma_2}}{\det R_f}.
		\end{align*}
		\item[(b)] A basis of $H^2(X, \omega)_f$ is given by:
		\begin{align*}
			\frac{\omega_f^{\sigma_1, \sigma_2}}{\det R_f}.
		\end{align*}
	\end{enumerate}
\end{conjecture}

A previous version of the manuscript incorrectly assumed that the matrix of the isomorphism $U_f^\vee \otimes \C \iso \bigoplus\limits_{j=1}^d U_{f,j}^\C$ is diagonal in certain natural bases. This led to a different rationality statement, namely that some multiples of $\omega_f^{\sigma_1}$ and $\omega_f^{\sigma_2}$ are rational. We would like to thank the anonymous referee for the previous version and Gyujin Oh for pointing out that this claim may be false in general.

We next summarize our evidence for this conjecture. The theoretical evidence in the case $[F:\Q] = 2$ is summarized in the following theorem.

\begin{theorem}[Corollary~\ref{cor:d=2det}, Corollary~\ref{cor:Stark_implies_conj}] \label{thmD}
	If the Fourier coefficients of $f$ are not rational, assume Stark's conjecture~\ref{conj:Stark}.
	\begin{enumerate}
		\item[(a)] The determinant of the basis $u_1^\vee \star f$, $u_2^\vee \star f$ is $E$-rational, i.e.\
		$$(u_1^\vee \star f) \wedge (u_2^\vee \star f) \in \wedge^2 H^1(X, \omega)_f \subseteq  \wedge^2 H^1(X, \omega)_f \otimes \C.$$
		\item[(b)] The cohomology class $(u_1^\vee \wedge u_2^\vee) \star f$ is $E$-rational, i.e.\ belongs to $H^2(X, \omega)_f$.
	\end{enumerate}
\end{theorem}

In fact, the rationality of $(u_1^\vee \wedge u_2^\vee) \star f$ is equivalent to Stark's conjecture for the trace 0 adjoint representation of $f$. Therefore, we may think of Conjecture~\ref{conj:B} as a refinement of Stark's conjecture for this representation. We thank Samit Dasgupta for suggesting this phrasing.

See Section~\ref{sec:Stark_conj} for generalizations of these results and further evidence in the case $[F:\Q] > 2$.

{\bf Numerical evidence.} The next goal of the paper is to verify the rationality of the classes
\begin{equation}
	u_1^\vee \star f, u_2^\vee \star f \in H^1(X, \omega)_f \otimes \C
\end{equation}
numerically. These cohomology classes are a linear combination of $\omega_f^{\sigma_1}$, $\omega_f^{\sigma_2}$, which are defined in equations~\eqref{eqn:omega_f^1}~\eqref{eqn:omega_f^2} as Dolbeault classes. We identify them with sheaf cohomology classes via the Dolbeault and the GAGA theorems. To check that they are ${E}$-rational is to show that the resulting sheaf cohomology classes come from base change of cohomology classes in~$H^1(X, \omega)_f$. The translation between Dolbeault and sheaf cohomology is not explicit enough to yield a satisfactory criterion for rationality. Worse yet, there seems to be no natural automorphic criterion to verify rationality. Indeed, the integral representations of Rankin--Selberg or triple product $L$-functions for Hilbert modular forms only involve cohomology classes $\omega_f^{J}$ where $J$ is the set of places where $f$ is dominant (see~\cite{Harris_periods_I} for details). Since parallel weight one forms are never dominant at any place, the cohomology classes we are interested in do not feature in these integral representations.

Instead, we consider an embedded modular curve $\iota \colon C \hookrightarrow X$ and check computationally in some cases that the restriction of $u_i^\vee \star f$ for $i=1,2$ to $C$ is rational, i.e.\
\begin{equation}
	\iota^\ast \left( u_i^\vee \star f  \right) \in H^1(C, \iota^\ast \omega) \otimes {E}.
\end{equation}
The drawback of this approach is that this restriction is non-zero only if the Hilbert modular form~$f$ is the base change of a modular form over $\Q$ (see, for example, Proposition~\ref{prop:Oda}). Let us hence assume that $f$ is the base change of a weight one modular form $f_0$. Then Conjecture~\ref{conj:B}~(a) can be restated in the simpler form (Conjecture~\ref{conj:real_quad}): the classes
\begin{align}
	\frac{\omega_f^{\sigma_1} + \omega_f^{\sigma_2}}{\log|u_{f_0}|} & \in H^1(X, \omega)_f \otimes \C,\\
	\frac{\omega_f^{\sigma_1} - \omega_f^{\sigma_2}}{\log|u_{f_0}^F|} & \in H^1(X, \omega)_f \otimes \C
\end{align}
belong to the rational structure $H^1(X, \omega)_f$, where $u_{f_0}$ is the unit associated to the adjoint representation of $f_0$ and $u_{f_0}^F$ is a unit associated to a twist of the adjoint representation of $f_0$. Finally, we check that this conjecture is equivalent to the single rationality statement:
\begin{equation}\label{eqn:restriction_intro}
	\frac{\iota^\ast(\omega_f^{\sigma_1})}{\log|u_{f_0}|} \in H^1(C, \iota^\ast \omega) \otimes E \subseteq H^1(C, \iota^\ast \omega) \otimes \C
\end{equation}
as long as $\iota^\ast(\omega_f^{\sigma_1}) \neq 0$ (c.f.\ Conjecture~\ref{conj:cons_of_conj_quadratic}, Proposition~\ref{prop:d=2_equivalent}).

We develop an algorithm to compute the trace of this cohomology class, i.e.\ an integral on the modular curve $C(\C)$ (see Conjecture~\ref{conj:cons_of_conj_quadratic}). We use results of Nelson~\cite{Nelson} to derive an expression for this integral (Theorem~\ref{thm:explicit_formula_for_integral}) which may be of independent interest. To use it, we give explicit formulas for the $q$-expansion of $f$ at other cusps when the level of $f$ is square-free (Theorem~\ref{thm:explicit_qexp_other_cusps}), generalizing results of Asai~\cite{Asai:Fourier}. Finally, we compute the integral numerically up to at least 15 digits of accuracy to give evidence for equation~\eqref{eqn:restriction_intro} in several cases (Tables~\ref{table}, \ref{table2}). 

The paper is organized as follows:
\begin{itemize}\setlength{\itemindent}{-.5cm}
	\item Section~\ref{sec:Stark} briefly discusses Stark's conjecture, introduces the unit group $U_f$, computes its rank, and gives a relation to a motivic cohomology group.
	\item Section~\ref{sec:der_Hecke} introduces the {\em derived Hecke action} and the generalization of the conjecture of Harris and Venkatesh~\cite{Harris_Venkatesh} to the Hilbert modular case (Conjecture~\ref{conj:A}~(1)).
	\item Section~\ref{sec:arch_real} introduces partial complex conjugation operators on cohomology and the archimedean conjectures (Conjectures~\ref{conj:A}(2) and~\ref{conj:B}).
	\item Section~\ref{sec:Stark_conj} discusses how the results of Stark and Tate give evidence for the archimedean conjecture, proving Theorems~\ref{thmB} and \ref{thmD}.
	\item Section~\ref{sec:BC} discusses base change cases, proves Theorems~\ref{thm:explicit_formula_for_integral} and~\ref{thm:explicit_qexp_other_cusps}, and provides numerical evidence for the archimedean conjecture.
	\item Appendix~\ref{app:Prasanna-Venkatesh} explains how Conjecture~\ref{conj:A}~(2) fits in the framework of Prasanna--Venkatesh~\cite{Prasanna_Venkatesh} and gives a version of this conjecture for partial weight one Hilbert modular forms.
\end{itemize}

Sections~\ref{sec:der_Hecke} and~\ref{sec:arch_real} are independent of one another and hence may be read in any order. The reader who wants to understand the full statements of the two conjectures as fast as possible may just skim Section~\ref{subsec:Stark_for_HMF} and proceed directly to these two sections.

\addtocontents{toc}{\protect\setcounter{tocdepth}{1}}

\subsection*{Acknowledgements}

First and foremost, I am grateful to my advisor, Kartik Prasanna, for countless helpful discussions, encouragement to pursue numerical evidence for the conjectures, and comments on various drafts of the paper. I would also like to thank Henri Darmon, Michael Harris, and Akshay Venkatesh for discussions related to a previous version of the paper. Finally, I would like to thank the referee for carefully reading the paper and providing helpful comments.

This work was supported by the National Science Foundation grant DMS-2001293, an Allen Shields Fellowship at the Department of Mathematics of the University of Michigan, and the Rackham Predoctoral Fellowship at the University of Michigan.

\addtocontents{toc}{\protect\setcounter{tocdepth}{2}}

\newpage
\begin{spacing}{.98}
\tableofcontents
\end{spacing}

\newpage

\section{Stark units and Stark's conjecture}\label{sec:Stark}

The goal of this section is to introduce the unit group $U_f$ mentioned in the introduction, compute its rank, and discuss its relation to motivic cohomology. We start by briefly recalling the definition of Stark units and Stark's conjecture. We then compute the unit group explicitly in the case of Hilbert modular forms.

\subsection{Stark units}

We follow~\cite{StarkII} to introduce the group of {\em Stark units} associated to an Artin representation. We caution the reader that the representations in loc.\ cit.\ are right representations, whereas we consider left representations, which leads to some discrepancies in notation. See also Dasgupta's excellent survey \cite{Dasgupta}.

Consider any Artin representation, i.e.\ a representation of the absolute Galois group $G_{\Q}$ which factors through a finite Galois extension $L$ of $\Q$:
\begin{center}
	\begin{tikzcd}
		G_{\Q} \ar[rr, "\varrho"] \ar[rd, swap, "\Res_L"] &&  \GL(M) \\
		& G_{L/\Q} \ar[ru, swap, "\varrho"]
	\end{tikzcd}
\end{center}
where $M$ is a free $\O_{E}$-module of rank $n$ and $E$ is a finite extension of $\Q$. We often write $G$ for the Galois group  $G_{L/\Q}$ and $U_L$ for the group of units of $\O_L$. 

\begin{definition}\label{def:Stark_units}
	The group of {\em Stark units} associated to $\varrho \colon G_{L/\Q} \to \GL(M)$ is:
	$$U_L[\varrho] = \Hom_{\O_E[G]}(M, U_L \otimes_\Z \O_E).$$
\end{definition}

We will soon check that $U_L[\varrho]$ depends only on $\varrho$ and not on the choice of $L$. To describe the group $U_L[\varrho]$ in more detail, we first need to understand the structure of $U_L$ as a $G_{L/\Q}$-module.

Fix an embedding $\tau \colon L \hookrightarrow \C$ which induces a complex conjugation $c_0$ of $L$. Note that $\rank U_L + 1 = \# (G/\langle c_0 \rangle)$ by Dirichlet's units theorem.

\begin{lemma}[{Minkowski's unit theorem, \cite[Lemma 2]{StarkII}}]\label{lemma:Minkowski}
	There is a unit $\epsilon$ of $L$, fixed by~$c_0$, such that there is only one relation among the $\rank U_L + 1$ units $\epsilon^{\sigma^{-1}}$ for $\sigma \in G/\langle c_0 \rangle$, and this relation is
	$$\prod_{\sigma \in G/ \langle c_0\rangle} \epsilon^{\sigma^{-1}} = \pm 1.$$
\end{lemma}

\begin{definition}\label{def:Minkowski}
	A unit whose existence is guaranteed by Lemma~\ref{lemma:Minkowski} is called a {\em Minkowski unit} of~$L$ with respect to $\tau \colon L \hookrightarrow \C$.
\end{definition}

\begin{corollary}\label{cor:U_L_via_Mink_units}
	The $\log$ map induces a $G$-equivariant isomorphism:
	\begin{align*}
		U_L/U_L^\tors & \overset \iso \longrightarrow \frac{\Z[\log(|\tau(\epsilon^{\sigma^{-1}})|) \ | \ \sigma \in G/\langle c_0 \rangle]}{\left\langle \sum\limits_{\sigma \in G/\langle c_0 \rangle} \log(|\tau(\epsilon^{\sigma^{-1}})|) \right\rangle},
	\end{align*}
	(the numerator on the right hand side is the free abelian group in those variables) and there is also a $G$-equivariant isomorphism:
	\begin{align*}
		\Ind_{\langle c_0 \rangle}^G \Z & \overset\iso\to \Z[\log(|\tau(\epsilon^{\sigma^{-1}})|) \ | \ \sigma \in G/\langle c_0 \rangle], \\
		(f \colon G/\langle c_0 \rangle \to \Z) & \mapsto \sum_{\sigma \in G/\langle c_0 \rangle} f(\sigma \langle c_0 \rangle) [\log(|\tau(\epsilon^{\sigma^{-1}})|)].
	\end{align*}
	In particular,
	$$U_L/U_L^{\mathrm{tors}} \iso \Ind_{\langle c_0 \rangle}^{G} \Z - \Z \text{ as a representation of }G = G_{L/\Q}.$$
\end{corollary}

We now compute the rank of $U_L[\varrho]$ and find a natural basis for $U_L[\varrho] \otimes_{\O_E} E$, given a basis of $M_E = M \otimes_{\O_E} E$. Let
\begin{equation}
	a = \dim_E M_E^{\langle c_0 \rangle}.
\end{equation}
Note that $a = (\mathrm{Tr}\varrho(1) + \mathrm{Tr} \varrho(c_0))/2$, so since any two complex conjugations of $L$ are conjugate, this number is independent of the choice of $c_0$. We write $b = n - a$ where $n = \dim_E M_E$.

\begin{proposition}\label{prop:Stark_basis}
	Suppose $\varrho$ does not contain a copy of the trivial representation. Then 
	$$U_L[\varrho] \otimes E \iso (M_E^{\langle c_0 \rangle})^\vee$$ 
	and hence $\rank U_L[\varrho] = a$. 
	
	Moreover, if $m_1, \ldots, m_a$ is a basis of $M_E^{\langle c_0 \rangle}$ and we complete it to a basis $m_1, \ldots, m_n$ of $M_E$ such that $\varrho(c_0)=\begin{pmatrix}
		I_a & 0 \\
		0 & -I_b
	\end{pmatrix}$ in this basis, then the corresponding basis of $U_L[\varrho] \otimes_{\O_E} E$ consists of the homomorphisms $\varphi_1, \ldots, \varphi_a$ defined by:
	\begin{equation}
		\varphi_i(m_j) = \prod_{\sigma \in G} (\epsilon^{\sigma^{-1}})^{a_{ij}(\sigma)} \in U_L \otimes E,
	\end{equation}
	where 
	$$\varrho(\sigma) = (a_{ij}(\sigma))_{i,j} \text{ in the basis }m_1, \ldots, m_n.$$
\end{proposition}
\begin{proof}
	We have that
	\begin{align*}
		U_L[\varrho] \otimes_{\O_E} E & = \Hom_{E[G]}(M_E, U_L \otimes_\Z E) \\
		& = \Hom_{E[G]} \left( M_E, \Ind_{\langle c_0 \rangle}^G E - E \right) & \text{Corollary~\ref{cor:U_L_via_Mink_units}} \\
		& = \Hom_{E[G]} \left(M_E, \Ind_{\langle c_0 \rangle}^G E\right) & \text{$\varrho$ does not contain the trivial rep.} \\
		& = \Hom_{E[\langle c_0 \rangle]}(M_E, E) & \text{Frobenius reciprocity} \\
		& = (M_E^{\langle c_0 \rangle})^\vee.
	\end{align*}
	Now, pick a basis $m_1, \ldots, m_n$ of $M$ such that  $\varrho(c_0) = \begin{pmatrix}
		I_a & 0 \\
		0 & -I_b
	\end{pmatrix}$ in it. By definition of the matrix ~$(a_{ij}(\sigma))_{i,j}$, 
	$$\varrho(\sigma) m_j = \sum_{k=1}^n a_{kj}(\sigma) m_k.$$
	Hence a map $\varphi \in \Hom_{\O_E}(M, U_L \otimes_\Z \O_E)$ is $G$-equivariant if and only if:
	\begin{equation}\label{eqn:G-equiv}
		(\varphi(m_j))^\tau =  \varphi(\varrho(\tau) m_j) = \varphi\left( \sum\limits_{k=1}^n a_{kj}(\tau) m_k  \right) =  \prod_{k=1}^n \varphi(m_k)^{a_{kj}(\tau)}
	\end{equation}
	(where the group of units is written multiplicatively).
	
	We check that each $\varphi_i$ defined above satisfies this equation. Let 
	$$u_{ij} =  \prod_{\sigma \in G} (\epsilon^{\sigma^{-1}})^{a_{ij}(\sigma)} \in U_L \otimes \O_E.$$ Then:
	\begin{align*}
		u_{ij}^\tau & = \left(\prod_{\sigma \in G} (\epsilon^{\sigma^{-1}})^{a_{ij}(\sigma)} \right)^\tau \\
		& = \prod_{\sigma \in G} (\epsilon^{\tau \sigma^{-1} })^{a_{ij}(\sigma)} \\
		& = \prod_{\sigma' \in G} (\epsilon^{(\sigma')^{-1}})^{a_{ij}(\sigma' \tau)}  & \text{for }\sigma' = \sigma \tau^{-1} \\
		& = \prod_{\sigma' \in G} (\epsilon^{(\sigma')^{-1}})^{\sum\limits_{k=1}^n a_{ik}(\sigma') a_{kj}(\tau)} \\
		& = \prod_{k=1}^n {\underbrace{\left(\prod_{\sigma \in G} (\epsilon^{\sigma^{-1}})^{a_{ik}(\sigma)}\right)}_{u_{ik}}}^{a_{kj}(\tau)} & \text{for }\sigma' = \sigma \\
		& = \prod_{k=1}^n u_{ik}^{a_{kj}(\tau)}.
	\end{align*}
	This shows that the functions $\varphi_i$ given by $\varphi_i(m_j) = u_{ij}$ are $G$-equivariant~(\ref{eqn:G-equiv}). Indeed:
	$$\varphi_i(m_j)^\tau = u_{ij}^\tau = \prod_{k=1}^n u_{ik}^{a_{kj}(\tau)} = \prod_{k=1}^n \varphi_i(m_k) ^{a_{kj}(\tau)}.$$
	Hence $\varphi_1, \ldots, \varphi_a \in U_L[\varrho]$.
	
	Tracing through the isomorphism
	$$U_L[\varrho] \otimes_{\O_E} E \iso (M_E^{\langle c_0 \rangle})^\vee$$
	established above, we see that
	$$\varphi_i \mapsto m_i^\vee \quad \text{for }i = 1, \ldots, a,$$
	where $m_i^\vee$ is a basis of $M_E^\vee$ dual to the basis $m_i$ of $M_E$. Since this is an isomorphism and $m_1, \ldots, m_a$ is a basis of $M_E^{\langle c_0 \rangle}$, $\varphi_1, \ldots, \varphi_a$ is a basis of $U_L[\varrho] \otimes E$.
\end{proof}

\begin{corollary}
	Suppose $\varrho \colon G_\Q \to \GL(M)$ is an Artin representation. Then $U_L[\varrho] \otimes E$ is independent of the choice of splitting field $L/\Q$.
\end{corollary}
\begin{proof}
	For an extension $L'/L$, the natural inclusion $U_L \hookrightarrow U_{L'}$ induces an inclusion $U_L[\varrho] \to U_{L'}[\varrho']$. By~Proposition~\ref{prop:Stark_basis}, $\dim U_L[\varrho] \otimes E = \dim U_{L'}[\varrho'] \otimes E$, which completes the proof.
\end{proof}

We will later be interested in the reduction of $U_L[\varrho]$ modulo $\p^n$ for a prime $\p$ of $E$. For now, we just remark that the following follows from Proposition~\ref{prop:Stark_basis}.

\begin{corollary}
	Let $t = \# U_L^{\rm tors}$ and $p$ be a prime not dividing $t$. Then $U_L[\varrho] \otimes_\Z \Z_p[\frac{1}{t}]$ is a free $\O_E \otimes_\Z \Z_p[\frac{1}{t}]$-module of rank $d$. Hence, for a prime $\mathfrak p$ of $E$ above $p$, $U_L[\varrho] \otimes \O_E/\p^n$ is a free $(\O_E/\p^n)$-module of rank $d$.
\end{corollary}
\begin{proof}
	This follows immediately from Proposition~\ref{prop:Stark_basis} and the structure theorem for modules over PIDs.
\end{proof}

\subsection{Stark's conjecture~\cite{StarkII, Tate:Stark}}\label{subsec:Stark_conj}

We give a brief summary of the results and conjectures on special values of Artin $L$-functions.

For any Artin representation $\varrho \colon G_{L/\Q} \to \GL(M)$
where $M$ is an $n$-dimensional $E$-vector space and an embedding $E \hookrightarrow \C$, we consider the $L$-function $L(s, \varrho)$ of $\varrho$. If we need to make the embedding $\iota \colon E \hookrightarrow \C$ explicit, we write $L(s, \varrho, \iota)$ for $L(s, \varrho)$.

The completed $L$-function is then:
\begin{equation}
	\Lambda(s, \varrho) = \left(\frac{f_\varrho}{\pi^n} \right)^{s/2} \Gamma(s/2)^a \Gamma((s+1)/2)^b L(s, \varrho)
\end{equation}
where:
\begin{align}
	f_\varrho & = \text{Artin conductor of }\varrho, \\
	a & = \dim_E M_E^{\langle c_0 \rangle},  & \text{(as above)}\\
	b & = n-a.
\end{align}

It satisfies a functional equation of the form:
\begin{align*}
	\Lambda(1- s, \overline{\varrho}) = W(\varrho) \Lambda(s, \varrho)
\end{align*}
where $|W(\varrho)| = 1$.

Stark gives a formula for the special value of $L$ at $s = 1$ (or, equivalently, the residue of the pole at $s = 0$). Associated to the units $u_{ij}$ in Proposition~\ref{prop:Stark_basis} is a regulator defined in terms of their logarithms.

Fix an embedding $\tau \colon L \hookrightarrow \C$ and let $c_0$ be the complex conjugation associated to $\tau$. Define
\begin{align*}
	\log \colon \C \iso L \otimes_{\tau} \C & \to \R \\
	z & \mapsto \log|z|
\end{align*}
and extend it linearly to
\begin{align*}
	\log \colon (L \otimes_{\tau} \C) \otimes ({E} \otimes_{\iota} \C)  & \to \C \\
	z \otimes \lambda & \mapsto \lambda \log|z|.
\end{align*}
Thus for $x \otimes y \in L \otimes {E}$, we write
\begin{equation}\label{eqn:log_def}
	\log|\tau \otimes \iota(x \otimes y)| = \iota(y) \cdot \log|\tau(x)| \in \C.
\end{equation}
We often make the choice of embeddings $\iota$ and/or $\tau$ implicit in the notation and write simply $\log|\tau(-)|$ or $\log|(-)|$ for $\log|\tau \otimes \iota (-)|$.

\begin{definition}\label{def:Stark_regulator}
	The {\em Stark regulator matrix} associated to $\varrho$ (and the embeddings $\tau \colon L \hookrightarrow \C$ and $\iota \colon E \hookrightarrow \C$) is
	$$R(\varrho) = (|\log(\tau \otimes \iota (u_{ij}))| )_{1 \leq i,j \leq a}.$$
	
	Abstractly, there is a perfect pairing
	\begin{align*}
		U_L[\varrho] \times M^{c_0} & \to \C \\
		(\varphi, m) & \mapsto \log(|\tau \otimes \iota(\varphi(m))|)
	\end{align*}
	via Proposition~\ref{prop:Stark_basis} and $R(\varrho)$ is the matrix of this pairing.
\end{definition}

\begin{conjecture}[Stark, {\cite{StarkII, Tate:Stark}}]\label{conj:Stark}
	If $\varrho$ does not contain the trivial representation, then
	$$L(1, \varrho) = \frac{W(\overline \varrho) 2^a \pi^b}{f_\varrho^{1/2}} \cdot \theta(\overline\varrho) \cdot \det R(\overline\varrho),$$
	for some $\theta(\overline \varrho) \in \Q(\Tr\varrho)^\times$, where $\Q(\Tr\varrho)$ is the field generated by the values of the character of $\varrho$.
\end{conjecture}

\begin{remark}
	The assumption that $\varrho$ does not contain the trivial representation is completely innocuous. Indeed, $L(s,\chi_{1,L}) = \zeta_L(s)$, so the value at $s = 1$ is given by the class number formula for $L$. Moreover, $L(s, \varrho_1 \oplus \varrho_2) = L(s, \varrho_1) \cdot L(s, \varrho_2)$.
\end{remark}

Stark's conjecture is known for representations with rational characters.

\begin{theorem}[Stark~{\cite[Theorem 1]{StarkII}}, Tate~{\cite[Corollary II.7.4]{Tate:Stark}}]\label{thm:Stark}
	Conjecture~\ref{conj:Stark} is true for representations $\varrho$ whose characters take rational values.
\end{theorem}

\subsection{Stark units for Hilbert modular forms}\label{subsec:Stark_for_HMF}

We now discuss Stark units for Artin representations associated to weight one Hilbert modular forms. Let $F$ be a totally real field. By \cite{Rogawski_Tunnell}, normalized weight one Hilbert modular eigenforms $f$ with Fourier coefficients in~$\O_{E_f}$ correspond to 2-dimensional odd irreducible Artin representations
\begin{center}
	\begin{tikzcd}
		G_{F} \ar[rr, "\varrho_f"] \ar[rd, swap, "\Res_L"] &&  \GL(M) \\
		& G_{L/F} \ar[ru, swap, "\varrho_f"]
	\end{tikzcd}
\end{center}
where $M$ is a $\O_{E}$-module of rank 2 and $E$ is a finite extension of $E_f$. By enlarging $L$ if necessary, we may assume that $L$ is Galois over $\Q$. We write $G = G_{L/\Q}$ and $G' = G_{L/F}$ for simplicity. 

As in the previous section, fix an embedding $\tau \colon L \hookrightarrow \C$ which induces a complex conjugation~$c_0$ of~$L$. Note that $c_0$ necessarily lies in $G'$ because $F$ is totally real. Since $\varrho_f$ is an odd representation,
$$\varrho_f(c_0) \text{ is conjugate to } \begin{pmatrix}
	1 & 0 \\
	0 & -1
\end{pmatrix}.$$

Consider the {\em adjoint representation of $\varrho$}, i.e.\
\begin{align*}
	\Ad \varrho_f \colon G_{L/F} & \to \GL(\End(M)) \\
	\sigma & \mapsto (T \mapsto \varrho(\sigma) T \varrho(\sigma)^{-1}).
\end{align*}
We note that if $T$ has trace 0, then so does $\varrho(\sigma) T \varrho(\sigma)^{-1}$. The representation is hence reducible, and we define the {\em trace zero adjoint representation} as
$$\Ad^0 \varrho_f  \colon G_{L/F} \to \GL(\End^0(M)),$$
where $\End^0(M) = \{T \colon M \to M \ | \ \Tr T = 0 \}$. This is a 3-dimensional representation.

Choosing a basis of $M_E$ such that $\varrho(c_0) = \begin{pmatrix}
	1 & 0 \\
	0 & -1
\end{pmatrix}$, we see that
$$(\Ad \varrho)(c_0) \begin{pmatrix}
	a & b \\
	c & -a
\end{pmatrix} =  \begin{pmatrix}
	1 & 0 \\
	0 & -1
\end{pmatrix}  \begin{pmatrix}
	a & b \\
	c & -a
\end{pmatrix} \begin{pmatrix}
	1 & 0 \\
	0 & -1
\end{pmatrix} = \begin{pmatrix}
	a & -b \\
	-c & -a
\end{pmatrix}.$$
Hence $\rank \left( (\Ad^0 \varrho_f)^{\langle c_0 \rangle}\right) = 1$.

\begin{definition}\label{def:Stark_units_ad}
	Let $U_L$ be the units of $L$ and $\O = \O_{E}$ be the ring over which $\varrho_f$ is defined. The group of {\em Stark units} associated to $f$ is
	$$U_f = \Hom_{\O[G_{L/F}]}(\Ad^0 \varrho_f, U_L \otimes_\Z \O).$$
	We sometimes write $\Ad^\ast \varrho = \Hom_{\O[G_{L/F}]}(\Ad^0 \varrho_f, \O)$, so that $U_f = \Ad^\ast \varrho \otimes_{\Z[G_{L/F}]} U_L$.
\end{definition}

Write $\sigma_1, \ldots, \sigma_d \in G$ for representatives of $G/G'$. Having fixed an embedding $\tau \colon L \hookrightarrow \C$, we have embeddings $\tau_j = \tau^{\sigma_j} \colon L \hookrightarrow \C$. We sometimes identify $\sigma_j$ with the embedding $\tau_j|_F \colon F \hookrightarrow \R$. We write $c_j = \sigma_j c_0 \sigma_j^{-1}$ for the complex conjugation associated to $\tau_j$.

\begin{corollary}\label{cor:Stark_basis_HMF}
	Suppose that $\varrho_f$ is irreducible. Then: $$U_f = U_L[\Ind_{G'}^G  \Ad^0 \varrho_f],$$
	is the group of Stark units associated to the $3d$-dimensional Artin representation $\Ind_{G'}^G  \Ad^0 \varrho_f$. Therefore,
	$$U_f \iso ((\Ind_{G'}^G  \Ad^0 \varrho_f)^{\langle c_0 \rangle})^{\vee} \iso \bigoplus_{j=1}^d ((\Ad^0 \varrho_f)^{\langle c_j \rangle})^\vee$$
	and hence
	$$\rank U_f = d.$$
	
	Moreover, for each $j=0, \ldots, d$, fix a basis $m_{1, j}, m_{2, j}, m_{3,j}$ of $\Ad^0 \varrho_f$ such that $\varrho(c_j) = \begin{pmatrix}
		I_1 & 0 \\
		0 & -I_2
	\end{pmatrix}$ in this basis, and consider the basis 
	$$\{\sigma_j m_{i, j} \ | \ j = 1,\ldots, d, i = 1,2,3\} \text{ of }\Ind_{G'}^G \Ad^0 \varrho.$$ 
	Let $a^0(\sigma)$ be the matrix of $\Ad^0 \varrho_f(\sigma)$ in the basis $\{m_{0, i}\}$ and write $P_{j}$ for the change of basis matrix from $\{m_{i,0}\}$ to $\{m_{i,j}\}$. Then there is a basis $\varphi_1, \ldots, \varphi_d$ of $U_f$ defined by Proposition~\ref{prop:Stark_basis} such that
	$$\varphi_j(\sigma_k m_{1, k}) = \prod_{\sigma' \in G'} \left(\epsilon^{( \sigma_k \sigma' \sigma_j^{-1} )^{-1}}\right)^{  (P_k a^0(\sigma') P_j^{-1})_{11}  }$$
\end{corollary}
\begin{proof}
	We have that:
	\begin{align*}
		U_f	& = \Hom_{\O[G']}(\Ad^0 \varrho_f, U_L \otimes_\Z \O) \\
		& = \Hom_{\O[G]}(\Ind_{G'}^G \Ad^0 \varrho_f, U_L \otimes_\Z \O) & \text{Frobenius reciprocity} \\
		& = U_L[\Ind_{G'}^G  \Ad^0 \varrho_f].
	\end{align*}
	
	Since $\varrho_f$ is irreducible, $\Ad^0 \varrho_f$ does not contain a copy of the trivial representation. We may hence apply Proposition~\ref{prop:Stark_basis} to the Artin representation $\Ind_{G'}^G \Ad^0 \varrho_f$ to get the result. Finally:
	\begin{align*}
		(\Ind_{G'}^G \Ad^0 \varrho_f)^{c_0}  & = \bigoplus_{j=1}^d \left(\sigma_j \Ad^0 \varrho_f  \right)^{c_0} \\
		& = \bigoplus_{j=1}^d \left(\Ad^0 \varrho_f   \right)^{c_j},
	\end{align*}
	completing the proof of the first part.
	
	It remains to prove the final assertion. To compute the action of an element $\sigma \in G$ on $\sigma_j \Ad^0 \varrho$, we find $\sigma_{k}$ and $\sigma' \in G'$ such that $\sigma \sigma_j = \sigma_{k} \sigma'$ and send
	$$\sigma_j m \mapsto \sigma_{k} (\sigma' m) \in \sigma_{j'}  \Ad^0 \varrho_f.$$
	
	By Proposition~\ref{prop:Stark_basis}, for $1 \leq j, k \leq d$:
	$$\varphi_j(\sigma_{k} m_{1,k}) = u_{jk} = \prod_{\sigma \in G} (\epsilon^{\sigma^{-1}})^{a_{jk}(\sigma)},$$
	where $a_{jk}(\sigma)$ is the matrix of $\Ind_{G'}^G \Ad^0 \varrho(\sigma)$ in the chosen basis. Then for $1 \leq j, k \leq d$:
	$$a_{jk}(\sigma) = \begin{cases}
		(P_k a^0(\sigma') P_j^{-1})_{11} & \text{if } \sigma_{k}^{-1} \sigma \sigma_{j} = \sigma'\text{ for some }\sigma' \in G', \\
		0 & \text{otherwise}.
	\end{cases}$$
	
	Therefore:
	$$u_{jk} = \prod_{\sigma' \in G'} \left(\epsilon^{ (\sigma_k \sigma' \sigma_j^{-1})^{-1} }  \right)^{a_{jk}(\sigma_k \sigma' \sigma_j^{-1})} =  \prod_{\sigma' \in G'} \left(\epsilon^{( \sigma_k \sigma' \sigma_j^{-1} )^{-1}}\right)^{  (P_k a^0( \sigma' ) P_j^{-1})_{11}  },$$
	as claimed.
\end{proof}

\begin{remark}
	The decomposition in Corollary~\ref{cor:Stark_basis_HMF} generalizes to any {\em plectic} Artin representation~\cite{Nekovar_Scholl_plectic}, i.e.\ an Artin representation of $G_F$ for a totally real field $F$. We have not used anything specific to Hilbert modular forms.
\end{remark}

\begin{remark}
	There is also a description of $U_f$ similar to \cite{Darmon_Lauder_Rotger}. For a chosen prime~$\mathfrak p$ of $F$, for each $\varphi_\sigma$, we may consider the component of $\varphi_\sigma(\Ad^0\varrho_f) \subseteq U_L$ on which a chosen Frobenius $\mathrm{Frob}_{\mathfrak p} \in G_{L/F}$ acts by $\alpha/\beta$ where $\alpha$ and $\beta$ are the ordered eigenvalues $\varrho_f(\mathrm{Frob}_{\mathfrak p})$. As in loc.\ cit.\, this space should be one-dimensional under extra assumptions; for example, that $\alpha \neq -\beta$. This description may be useful when considering a $p$-adic analogue of the conjecture, but we omit this here entirely.
\end{remark}

\subsection{Stark's conjecture for Hilbert modular forms}

We now state Stark's conjecture for the trace zero adjoint representation associated to a Hilbert modular form of parallel weight one. 

\begin{definition}\label{def:Stark_reg_HMF}
	The {\em Stark regulator matrix} associated to (the trace zero adjoint representation of) $f$ is
	$$R_f = ( \log( |  u_{jk}  |  )  )_{1 \leq j, k \leq d},$$
	with
	$$u_{jk} =  \prod_{\sigma' \in G'} \left(\epsilon^{( \sigma_k \sigma' \sigma_j^{-1} )^{-1}}\right)^{  (P_k a^0(\sigma') P_j^{-1})_{11}  }$$
	(notation as in Corollary~\ref{cor:Stark_basis_HMF}). If we need to specify $f$, we write $u_{jk}^f$ for $u_{jk}$.
\end{definition}

\begin{proposition}\label{prop:Stark_for_HMF}
	Stark's conjecture~\ref{conj:Stark} for $\Ad^0\varrho_f$ is equivalent to the statement:
	$$L(1, \overline{\Ad^0 \varrho_f}) \sim_{E^\times} \frac{\pi^{2d}}{f_\varrho^{1/2}} \cdot  \det R_{f},$$
	where $f_\varrho$ is the conductor of $\varrho =\Ind_{G'}^G \Ad^0 \varrho_f$.
\end{proposition}

\begin{remark}
	In Section~\ref{sec:Stark_conj}, we will relate the adjoint $L$-function to the Petersson inner product of $f$. This will give evidence for our archimedean conjecture (Conjecture~\ref{conj:motivic_action}).
\end{remark}

\subsection{Examples}

The Stark unit group can be determined explicitly in many cases. We provide a few illustrative examples.

\begin{example}[Heegner units]
	The first example of Stark units comes from the theory of elliptic units. 
	
	Let $F = \Q$ and $K/\Q$ be an imaginary quadratic extension. For any Dirichlet character $\chi \colon  G_{H/K} \to \C^\times$ of~$K$, where $H/K$ is an abelian extension, there is an associated weight one form $f = \theta_\chi$, the {\em theta function} of $\chi$, such that
	$$L(s, \chi) = L(s, f).$$
	The Artin representation $\varrho$ associated to $f$ is the 2-dimensional representation:
	$$\varrho_f = \Ind_{G_{H/K}}^{G_{H/\Q}} \chi = \{\phi \colon G_{H/\Q} \to \C \ | \ \phi(\sigma \tau) = \chi(\sigma) \phi(\tau) \quad \text{ for $\sigma \in G_{H/K}$} \}.$$
	For the non-trivial element $c \in G_{K/\Q}$, we can define a character $\chi^c(\sigma) = \chi(c \sigma c)$.  Writing $1$ for the trivial representation and $\chi_0$ for $\chi \cdot (\chi^c)^{-1}$, we see that
	$$\Ad^0 \varrho_f \iso 1 \oplus \Ind_{G_{H/K}}^{G_{H/\Q}}  \chi_0.$$
	Since the unit group does not contain a copy of the trivial representation, this shows that
	$$U_f \iso U_H[\chi_0],$$
	the $\chi_0$-isotypic component of the units of $H$. For a Minkowski unit $\epsilon \in \O_H^\times$, the unit associated to~$f$ is:
	$$u_f = u_{\chi_0} =  \prod_{\sigma \in G_{H/K}} (\epsilon^{\sigma^{-1}})^{\chi_0(\sigma)}.$$
	In literature, this unit is often written additively as $u_{\chi_0} = \sum\limits_{\sigma \in G_{H/K}} \chi_0(\sigma)^{-1} u^\sigma \in U_H[\chi_0]$. {\em Elliptic units}, constructed using singular values of modular functions, provide an explicit construction of Minkowski units $u \in \O_H^\times$, and hence of Stark units $u_f$.
	
	The {\em logarithms} of these units appear as special values of the $L$-function of $\chi_0$, via Kronecker's second limit formula. This also has a $p$-adic analogue: the {\em $p$-adic logarithm} of $u_{\chi_0}$ accounts for the special value of the Katz $p$-adic $L$-function evaluated at the finite order character $\chi_0^{-1}$, which is outside of the range of interpolation~\cite[10.4, 10.5]{Katz:p-adic_int}. More generally, Darmon, Lauder, and Rotger conjecture \cite[Conjecture ES]{Darmon_Lauder_Rotger} that $p$-adic logarithms of other Stark units associated to weight one modular forms appear in a formula for values of triple product $p$-adic $L$-functions outside the range of interpolation.
\end{example}

The following example is suitable for computations in the case $F = \Q$. In fact, it is the example where Harris--Venkatesh~\cite{Harris_Venkatesh} perform their computations. It is also a simple example where our archimedean conjecture (Conjecture~\ref{conj:motivic_action}) can be proved (Corollary~\ref{cor:conj_for_F=Q}).

\begin{example}[Units in cubic fields, $F = \Q$]\label{ex:cubic_fields}
	This example is discussed in \cite[Sec.~5.6]{Harris_Venkatesh}, but we recall it here in detail to provide context for the generalizations to $[F:\Q] = 2$ we make below.
	
	Let $K$ be a cubic field of signature $[1,1]$ and write $L$ for the Galois closure of $K$. Then $G_{L/\Q} \iso S_3$ and we may assume that $K$ is the fixed field $L^{(12)}$ of the action of the cycle $(12) \in S_3$ on $L$. 
	
	To give a 2-dimensional representation of $G_{L/\Q}$, we need to give a 2-dimensional representation of~$S_3$. There is a unique irreducible 2-dimensional representation: the {\em regular representation} $\varrho \colon G_{\Q} \twoheadrightarrow S_3 \to \GL(M) \iso \GL_2(\Z)$, obtained by considering the action of~$S_3$ on
	$$M = \left\{(x_1, x_2, x_3) \in \Z^3\ \middle| \ \sum x_i = 0 \right\}$$
	by permuting the coordinates.
	
	In the basis $e_1 = (1, 0, -1)$, $e_2 = (0,1,-1)$ of $M$, we have that:
	\begin{align*}
		\sigma = (12) & \mapsto S = \begin{pmatrix}
			0 & 1 \\
			1 & 0
		\end{pmatrix}, \\
		\tau = (123) & \mapsto T = \begin{pmatrix}
			-1 & -1 \\
			1 & 0
		\end{pmatrix}.
	\end{align*}
	Note that $\varrho$ is an {\em odd} Galois representation since $\det S = -1$. Therefore, there is a weight one modular form $f$ corresponding to $\varrho$.

	Recall that
	$$U_f = \Hom_{G_{L/\Q}}(\Ad^0 \varrho, U_L).$$
	Lemma 5.7 in \cite{Harris_Venkatesh} shows that
	\begin{align}
		U_f \otimes \Z\left[\frac{1}{6} \right] & \overset \iso \to U_K^{(1)} \otimes \Z\left[\frac{1}{6} \right] \\
		(\varphi \colon \Ad^0 \varrho \to U_L) & \mapsto \varphi(S), \nonumber
	\end{align}
	where $U_K^{(1)}$ are the norm 1 units of $K$.
	
	We recall the proof here. By definition
	$$\Ad^0 \varrho \iso \End^0(M),$$
	with the action of $S_3$ on the right hand side given by conjugation. Note that each element of $S_3$ gives an element of $\End(M)$ and we may use the $S_3$-invariant projection
	\begin{align*}
		\End(M) & \to \End^0(M) \\
		A & \mapsto A - (1/2) \mathrm{Tr}(A)  
	\end{align*}
	to get a spanning set for $\Hom^0(M,M)$ this way. Since the lengths of cycles are conjugation-invariant, we see that
	$$\Hom^0(M,M) \iso \mathrm{span}(\text{images of }(123), (132)) \oplus \mathrm{span}(\text{images of }(12), (13), (23)).$$
	One checks that $\mathrm{span}(\text{images of }(123), (132)) = \Z[e]$, where $e \in \Hom^0(M,M)$ is the common image of $(123)$ and $(132)$. We write $W = \mathrm{span}(\text{images of }(12), (13), (23))$. Hence
	$$\Ad^0 \varrho \iso \Z[e] \oplus W.$$
	Now, for any $S_3$-representation $V$:
	\begin{itemize}
		\item $\Hom_{S_3}(\Z[e], V) = V^{\sgn}$, the $\sgn$-isotypic part of $V$,
		\item $\Hom_{S_3}(W, V) \iso \{v \in V^{(12)} \ | \ v + (123)v + (132)v = 0\}$ via $\varphi \mapsto \varphi(S)$.
	\end{itemize}
	This shows that:
	$$U_f \iso U_L^{\sgn} \oplus U_K^{(1)}$$
	Since $\Q(\sqrt{\mathrm{disc}(L)}) = L^{\langle (123) \rangle}$, $U_L^{\sgn} = U_{\Q\left(\sqrt{\mathrm{disc}(L)}\right)}$ is a finite group of order at most 6. Hence
	$$U_f \otimes \Z \left[\frac{1}{6}\right] \iso U_K^{(1)} \otimes  \Z \left[\frac{1}{6}\right].$$
\end{example}

The following is the simplest example of explicit Stark units over real quadratic fields. It is the base change of Example~\ref{ex:cubic_fields} to a real quadratic field and one of the examples in which we will do the numerical computations later on.

\begin{example}[Units in cubic extensions of $F$ for ${[F:\Q] = 2}$]\label{ex:base_change_units}
	Consider $K$ as in Example~\ref{ex:cubic_fields} and consider a quadratic extension $F$ of $\Q$. Then $KF$ is a cubic extension of $F$ of signature $[2,2]$:
	
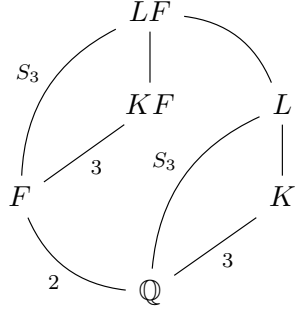
\begin{center}
		\begin{tikzcd}
			& LF \\
			& KF \ar[u, no head] & L \ar[lu, no head, bend right] \\
			F \ar[ruu, no head, bend left, "S_3"] \ar[ru, no head, swap, "3"] & & K \ar[u, no head] \\
			& \Q \ar[lu, no head, bend left, "2"] \ar[no head, swap, ru, "3"] \ar[ruu, no head, bend left, "S_3"]
		\end{tikzcd}
	\end{center}
	As above, we consider the Galois representation $\varrho \colon G_{LF/F} \iso S_3 \to \GL(M)$. If $f$ is the weight~one Hilbert modular form associated to $\varrho$, then one can check that
	\begin{align*}
		U_{f} \otimes \Z[1/6] & \iso \Hom_{ G_{LF/F}}(\Ad^0 \varrho, U_{LF}) \otimes \Z[1/6] \\
		& \iso U_{LF}^{\mathrm{sgn}} \otimes \Z[1/6] \oplus \{u \in U_{KF} \ | \ N_{F}^{KF} u = 1 \} \otimes \Z[1/6] \\
		& \iso (\{u \in U_{K} \ | \ N_{\Q}^{K} u = 1 \} \oplus \{u \in U_{KF} \ | \ u^\sigma = u^{-1}, N^{KF}_{F} u = 1\}) \otimes \Z[1/6],
	\end{align*}
	where we write $\Gal(F/\Q) = \langle \sigma \rangle$.
	The Hilbert modular form $f$ is the base change of the modular form $f_0$ associated to $K$ in the previous example. We will later prove a more general result of this form in Corollary~\ref{cor:units_for_real_quadratic_BC}.
\end{example}

Finally, we present the ``simplest'' non base change example where explicit Stark units are available over real quadratic fields. It is a direct analogue of Example~\ref{ex:cubic_fields}, but the Galois theory is more complicated.

\begin{example}[Units in cubic extensions of $F$ for ${[F:\Q]} = 2$, non base change]\label{ex:cubic_ext_units}
	We generalize Example~\ref{ex:cubic_fields} to the case $[F:\Q] = 2$ and a cubic extension $K$ of $F$ of signature $[2,2]$:
	\begin{center}
		\begin{tikzcd}
			& & L \\
			K \ar[rru, no head, bend left, "12"] \\
			F \ar[u, no head, "3"] \ar[rruu, swap, no head, "S_3 \times S_3"] \\
			& \Q \ar[lu, no head, bend left, "2"] \ar[uuur, no head, bend right] \ar[ruuu, no head, bend right, swap, "S_3 \wr C_2"]
		\end{tikzcd}
	\end{center}
	We may assume that $K = L^{S_3 \times \langle (12) \rangle}$. Consider the representation
	\begin{align*}
		\varrho = \sgn \boxtimes \mathrm{reg} \colon S_3^2 & \to \GL_2(\Z), \\
		(\sigma, (12)) & \mapsto \sgn(\sigma) \cdot \begin{pmatrix}
			0 & 1 \\
			1 & 0
		\end{pmatrix}, \\
		(\sigma, (123)) & \mapsto \sgn(\sigma) \cdot \begin{pmatrix}
			-1 & 1 \\
			-1 & 0
		\end{pmatrix}.
	\end{align*}
	Then $\varrho$ corresponds to a Hilbert modular form $f$ of parallel weight~one.
	
	As before,
	$$\Ad^0 \varrho \iso \Z[e] \oplus W$$
	and for any $S_3^2$-representation $V$:
	\begin{itemize}
		\item $\Hom_{S_3^2}(\Z[e], V) \iso V^{\sgn \boxtimes \sgn}$,
		\item $\Hom_{S_3^2}(W, V) \iso \{v \in V^{S_3 \times (12)} \ | \ v + (1,(123))v + (1, (132)) v = 0 \}$ with the isomorphism given by sending $\varphi \colon W \to V$ to $\varphi(S)$.
	\end{itemize}
	
	Therefore,
	\begin{align*}
		U_f & = \Hom_{S_3^2}(\Ad^0 \varrho, U_L) \\
		& = \Hom(\Z[e], U_L) \oplus \Hom(W, U_L) \\
		& = (U_L^{\sgn \boxtimes \sgn}) \oplus \{u \in U_K \ | \ N_F^K u = 1\}.
	\end{align*}
	
	We claim that the group $U_L^{\sgn \boxtimes \sgn}$ is torsion. If $u \in U_L^{\sgn \boxtimes \sgn}$, then $u$ is fixed by a subgroup $H \subseteq S_3^2$ of order 18 of elements $(\sigma, \sigma')$ such that $\sgn(\sigma) = \sgn(\sigma')$. One can check that $L^H = F(\sqrt{{\rm disc}(L/F)})$, which is a CM extension of $F$. Therefore, if $\Gal(F(\sqrt{\mathrm{disc(L/F)}})/F) \iso  
	\langle\tau \rangle$,
	$$U_L^{\sgn \boxtimes \sgn} \iso (U_{F(\sqrt{\mathrm{disc(L/F)}})})^{\tau = -1}.$$
	Since $F(\sqrt{\mathrm{disc(L/F)}})/F$ is CM, the ranks of the two unit groups are equal. On the other hand, if $u \in U_{F(\sqrt{\mathrm{disc(L/F)}})}^{\tau = -1}$ was a non-torsion element, then $u$ would generate an infinite subgroup of $U_{F(\sqrt{\mathrm{disc(L/F)}})}^{\tau = -1}$ which does not belong to $U_F$. This is a contradiction.
	
	Finally, let $N$ be the order of the torsion group $U_L^{\sgn \boxtimes \sgn}$. Then:
	\begin{equation}
		U_f \otimes \Z[1/N] \iso \{u \in U_K \ | \ N_F^K u = 1\} \otimes \Z[1/N].
	\end{equation}
	As expected by Corollary~\ref{cor:Stark_basis_HMF} this is a group of rank 2. In terms of the notation of Definition~\ref{def:Stark_reg_HMF}, the units $u_{11}$, $u_{12}$ give a basis of the last space. Identifying the units $u_{21}$, $u_{22}$ seems more difficult.
\end{example}

\subsection{Comparison with motivic cohomology}

This section is not used in the remainder of this paper. The general conjectures of Venkatesh~\cite{Venkatesh:Derived_Hecke} predict the action of the dual of a motivic cohomology group associated to the coadjoint motive of~$f$. We identify this motivic cohomology group with the group of Stark units $U_f$, analogously to~\cite[Sec.~2.8]{Harris_Venkatesh}. Some of this section is based on standard conjectures. 

\subsubsection{Motivic cohomology}

Let $k$ be any number field and $\O_k$ be its ring of integers. (In general, $\O_k$ could be any Dedekind domain and $k$ its field of fractions). Let $E$ be a field of characteristic 0.

For any Chow motive $M$ defined over $k$ with coefficients in $E$, we may define motivic cohomology groups (cf.\ \cite{Bloch} or \cite[Definition 3.4]{Mazza_Veovodsky_Weibel})
$$H^r_{\mathcal M_k}(M, E(n))$$
which are equipped with specialization maps to various cohomology theories, including \'etale cohomology:
$$H^r_{\mathcal M_k}(M, E(n)) \otimes E_\p \to H^r_{\textnormal{\'et}}(M, E_\p(n)).$$

Scholl~\cite[Theorem 1.1.6]{Scholl} proved that these have a subspace of {\em integral classes}
$$H^r_{\mathcal M_{\O_k}}(M, E(n)) \subseteq H^r_{\mathcal M_k}(M, E(n)).$$

We will be concerned with the case $r = 1$, $n = 1$. For the trivial motive $M = k$, conjecturally:
\begin{equation}
	H^1_{\M_{\O_k}}(k, E(1)) \iso U_k \otimes E.
\end{equation}
This statement is certainly true in all realizations; see, for example, \cite[4.3]{Nekovar_Beilinson_survey} or \cite[Corollary 4.2]{Mazza_Veovodsky_Weibel}.

\subsubsection{Motivic cohomology of the coadjoint motive}

Conjecturally, there is a 3-dimensional Chow motive $M_{\mathrm{coad}}$ with coefficients in~$E$, the {\em coadjoint motive of~$f$}, associated to the dual of the trace zero adjoint representation, $\Ad^\ast \varrho_f$. By definition, for any prime $\p$ of $E$, its $\p$-adic \'etale realization is isomorphic to:
$$H^\bullet_\lambda(M_{\mathrm{coad}} \times_\Q \overline \Q, E_\p) \iso \Ad^\ast \varrho_f \otimes_E E_\p$$
(concentrated in cohomological degree 0). Without loss of generality, we assume that $M_{\mathrm{coad}}$ is defined over $F$ (and not just $\overline{F_\lambda}$). 

\begin{remark}
	Motives associated to Hilbert modular forms were constructed in \cite{Blasius_Rogawski} in some cases where the weights are cohomological. Since weight one Hilbert modular forms are not cohomological, there is no known construction of the motive, but we assume that (at least) the coadjoint motive exists.
\end{remark}

According to~\cite{Venkatesh:Derived_Hecke, Prasanna_Venkatesh}, we should consider the motivic cohomology group
$$H^1_{\mathcal M_{\O_F}}(M_{\mathrm{coad}}, E(1)).$$
There is a natural map
$$H^1_{\mathcal M_{\O_F}}(M_{\mathrm{coad}}, E(1)) \to H^1_{\mathcal M_{\O_L}}(M_{\mathrm{coad}}, E(1))^{G_{L/F}}$$
and we will work with the codomain instead. According to \cite[(2.8)]{Harris_Venkatesh}, this map should be an isomorphism. In the proof of Proposition~\ref{prop:motivic_coh} below, we check this in the \'etale realization (the induced map is denoted by $i$).

For a prime $\p$ of $E$, the $\p$-adic \'etale realization map:
$$H^1_{\M_{\O_F}}(M_{\mathrm{coad}}, E(1)) \otimes \O_\p \to H^1_f(F, (\Ad^\ast \varrho_f \otimes \O_\p)(1))$$
is conjecturally an isomorphism~\cite[5.3(ii)]{Bloch_Kato}. Here, $H^1_f$ denotes the Bloch--Kato Selmer group~\cite{Bloch_Kato}. (We apologize for the clash of notation with the Hilbert modular form $f$ and hope that this does not cause confusion.) We compute the last group.

\begin{proposition}\label{prop:motivic_coh}
	We have that
	$$H^1_f(F, (\Ad^\ast \varrho_f \otimes \O_\p)(1)) \iso U_f \otimes \Q \otimes \O_\p$$
	for all $\p$ such that $N\p$ is coprime to $[L:F]$.
\end{proposition}
\begin{proof} 
	This argument is adapted from \cite[Lemma 4.5]{Harris_Venkatesh}. We claim that
	$$H^1_f(G_F, \Ad^\ast \varrho_f \otimes \O_\p) \iso (U_L \otimes \Q \otimes \Ad^\ast \varrho_f \otimes \O_\p)^{G_{L/F}}.$$
	Recall that $(U_L \otimes \Ad^\ast \varrho_f)^{G_{L/F}} = U_f$ by definition, so this will prove the proposition.
	
	We write $\Ad^\ast \varrho_\p$ for $\Ad^\ast \varrho_f \otimes \O_\p$ for simplicity. The (global) Bloch--Kato Selmer group $H^1_f$ is defined by the short exact sequence:
	\begin{center}
		\begin{tikzcd}
			0 \ar[r] & H^1_f(F, \Ad^\ast \varrho_\p(1)) \ar[r] & H^1(F, \Ad^\ast \varrho_\p(1)) \ar[r] & \bigoplus\limits_{v} \frac{H^1(F_v, \Ad^\ast \varrho_\p(1))}{H^1_f(F_v, \Ad^\ast \varrho_\p(1))}.
		\end{tikzcd}
	\end{center}
	where $H^1_f(F_v, \Ad^\ast \varrho_\p(1))$ are the local Bloch--Kato Selmer groups. The restriction maps to the subgroup $G_{\overline L/ L} \subseteq G_{\overline{F} /F}$ give a commutative diagram
	\begin{center}
		\begin{tikzcd}
			0 \ar[r] & H^1_f(F, \Ad^\ast \varrho_\p(1)) \ar[r] \ar[d, "i"] & H^1(F, \Ad^\ast \varrho_\p(1)) \ar[r] \ar[d, "j"] & \bigoplus\limits_{v} \frac{H^1(F_v, \Ad^\ast \varrho_\p(1))}{H^1_f(F_v, \Ad^\ast \varrho_\p(1))} \ar[d, "k"] \\
			0 \ar[r] & \left( H^1_f(L, \Ad^\ast \varrho_\p(1)) \right)^{G_{L/F}} \ar[r] & \left(H^1(L, \Ad^\ast \varrho_\p(1)) \right)^{G_{L/F}} \ar[r] & \left( \bigoplus\limits_{w} \frac{H^1(L_w, \Ad^\ast \varrho_\p(1))}{H^1_f(L_w, \Ad^\ast \varrho_\p(1))} \right)^{G_{L/F}}
		\end{tikzcd}
	\end{center}
	with exact rows. Since $\Ad^\ast \varrho_\p(1)$ is trivial as a $G_{\overline{L}/L}$-representation, we have that:
	\begin{align*}
		\left( H^1_f(L, \Ad^\ast \varrho_\p(1)) \right)^{G_{L/F}} & \iso  \left(\Ad^\ast \varrho_\p \otimes_{\O_\p} H^1(L, \O_\p(1)) \right)^{G_{L/F}} \\
		& \iso \left(\Ad^\ast \varrho_\p \otimes_{\O_\p} U_L \otimes \O_\p \otimes \Q \right)^{G_{L/F}},
	\end{align*}
	so we just need to show that the map $i$ is an isomorphism.
	
	Since $N\p$ is coprime to $[L:F]$, the restriction map $j$ is an isomorphism by a general group cohomology result~\cite[{I.2.4}]{Serre:Galois_coh}. By the Snake Lemma, this shows that $i$ is also injective.
	
	To show that it is surjective, we must show that $k$ is injective. In fact, for a place $w$ of $L$ above a place $v$ of $F$, the restriction map
	$$\frac{H^1(F_v, \Ad^\ast \varrho_\p(1))}{H^1_f(F_v, \Ad^\ast \varrho_\p(1))} \to \frac{H^1(L_w, \Ad^\ast \varrho_\p(1))}{H^1_f(L_w, \Ad^\ast \varrho_\p(1))}$$
	is split by the corestriction map divided by $[L_w : F_v]$ (since $[L_w:F_v]$ is invertible in $\O_\p$).
\end{proof}

\section{Derived Hecke operators on the special fiber}\label{sec:der_Hecke}

Let:
\begin{itemize}
	\item $f$ be a normalized Hilbert modular eigenform of parallel weight one, new of level $\mathfrak N$, with coefficients in the ring $\O_{E_f}$;
	\item  $\varrho = \varrho_f$ be the associated Artin representation, defined over $\O = \O_E$ where $E$ is a finite extension of $E_f$;
	\item $U_f$ be the group of {\em Stark units}, which has rank $d = [F:\Q]$ over $\O$;
	\item $\p$ be a prime of $\O_E$ such that $(p) = \p \cap {\Q}$ has good reduction in $F$ and $p$ is coprime to $\mathfrak N$, and let $k = \O_{E}/\p^n$.
\end{itemize}

We consider a smooth, compact, integral model $X = X_{1}(\mathfrak N)$ for the Hilbert modular variety associated to $F$ and the level~$\Gamma_1(\mathfrak N)$ (the level of $f$). Such integral models for the toroidal compactifications with the level structures considered here were developed in \cite{Dimitrov_Gamma_1_n}, following the standard methods of Rapoport~\cite{Rapoport}. They are defined over $\Z[1/N_{F/\Q} \mathfrak N]$, where $N_{F/\Q}$ denotes the norm from $F$ to $\Q$. See also \cite{Chai}, \cite{Dimitrov_arithmetic_aspects_of_HMF}, or \cite{Goren} for surveys on Hilbert modular varieties and Hilbert modular forms.

Let $\omega$ be the Hodge bundle on the integral Hilbert modular surface $X_{\Z[1/N_{F/\Q} \mathfrak N]}$, so that
$$f \in H^0(X_{\Z[1/ N_{F/\Q}\mathfrak N]}, \omega) \otimes_\Z \O_{E_f}.$$
In this section we construct an action of $U_f^\vee \otimes_{\O_{E}} k$ on the cohomology space
$$(H^\ast(X_{\Z[1/ N_{F/\Q}\mathfrak N]}, \omega) \otimes_\Z \O_{E})_{f} \otimes_{\O_{E}} k \iso H^\ast(X_k, \omega)_{f} $$
via derived Hecke operators on the special fiber and conjecture that it lifts to $\O_{E}$. This is an analogue of the Harris--Venkatesh conjecture~\cite{Harris_Venkatesh} for the coherent cohomology of the Hodge bundle on Hilbert modular varieties.

Recall (c.f.\ Section~\ref{subsec:Stark_for_HMF}) that the Artin representation associated to $f$ factors through a finite Galois extension $L/F$ and has coefficients in the integers $\O_E$ of a number field $E$, i.e.\ $\varrho_f \colon \Gal(L/F) \to \GL_2(\O_E)$. Let $q > 5$ be a prime and $\mathfrak q$ be a prime of $F$ above it such that $N\mathfrak q \equiv 1\ (p^n)$. We fix a choice of a prime ideal $\mathfrak Q$ of $L$ above $\mathfrak q$. We write $G' = \Gal(L/F)$ and $G = \Gal(L/\Q)$.

This configuration is summarized by the following diagram:
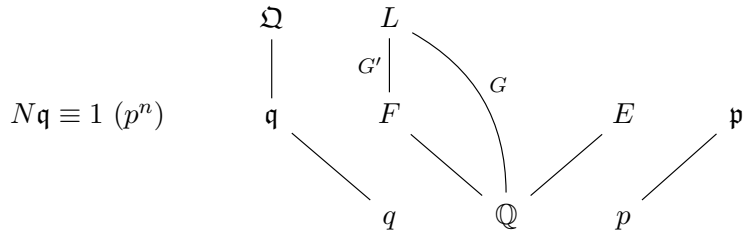
\begin{center}
	\begin{tikzcd}
		& \mathfrak Q & L \\
		N\mathfrak q \equiv 1\ (p^n)& \mathfrak q  \ar[u, no head] & F \ar[u, no head, "G'"] & & {E} & \mathfrak p \\
		& & q \ar[lu, no head]&  \Q \ar[lu, no head] \ar[ru, no head]  \ar[luu, no head, bend right, swap, "G"]& p \ar[ru, no head]
	\end{tikzcd}
\end{center}

We will describe:
\begin{itemize}
	\item a {\em map} 
	$$\theta^\vee_{\mathfrak q} \colon \bigoplus_{\sigma \in G/G'} U_{f, \sigma}^\vee \to U_f^\vee \otimes k$$
	in Section~\ref{subsec:red_Stark} (Proposition~\ref{prop:dual_units_mod_p});
	
	\item an action of the domain via {\em derived Hecke operators}: 
	$$T_{\sigma \mathfrak q, z} \colon H^q(X_k, \omega)_f \to H^{q+1}(X_k, \omega)_f$$ 
	associated to $z \in U_{f, \sigma}^\vee$ in Sections~\ref{subsec:Shimura_class}, ~\ref{subsec:derived_Hecke} (Definition~\ref{def:action_mod_p});
\end{itemize}
and conjecture that the resulting action of $U_f^\vee \otimes k$ lifts to characteristic 0 in Section~\ref{subsec:conj_mod_p} (Conjecture~\ref{conj:modulo}).

\subsection{Dual Stark units mod $\p^n$}\label{subsec:red_Stark}

We start by describing the group $U_f^\vee \otimes_{\O_{E}} k$. The description will depend on a choice of a Taylor--Wiles prime $\mathfrak q$ of $F$.

\subsubsection{Taylor--Wiles primes}

Suppose $\p$ is a prime of ${E}$ above $p$ and for any $n$, consider
$$k = \O_{E}/\p^n.$$

\begin{definition}
	A {\em Taylor--Wiles} prime for $f$ of level $n \geq 1$ consists of the following data:
	\begin{enumerate}
		\item a prime $\mathfrak q$ of of $F$, relatively prime to the level of $f$, such that $N \mathfrak q \equiv 1\ (p^n)$,
		\item a choice $(\alpha, \beta) \in \mathbb F_\p^2$ with $\alpha \neq \beta$ such that $$\overline{\varrho}(\mathrm{Frob}_{\mathfrak q}) = \begin{pmatrix}
			\alpha & 0 \\
			0 & \beta
		\end{pmatrix},$$
		where $\overline \varrho$ is the reduction of $\varrho$ modulo $\mathfrak p$.
	\end{enumerate}
\end{definition}

If $\mathfrak q$ is a Taylor--Wiles prime, $(\O_F/\mathfrak q)^\times$ contains a subgroup $\Delta \iso \Z/p^n \Z$ of size $p^n$. We frequently denote it by $(\O_F/\mathfrak q)^\times_{p^n}$.

We also write
\begin{equation}
	k\langle 1 \rangle_{\mathfrak q} = k \otimes (\O_F/\mathfrak q)^\times_{p^n}, \quad k\langle -1 \rangle_{\mathfrak q} = \Hom((\O_F/\mathfrak q)^\times_{p^n}, k),
\end{equation}
both non-canonically isomorphic to $k$. When the underlying prime $\mathfrak q$ is clear, we drop it from the notation.

Finally, for any $\Z$-module $M$, we write
\begin{equation}
	M \langle m \rangle = M \otimes_\Z k\langle m \rangle \quad \text{for } m = \pm 1.
\end{equation}
For example, $\mathbb F_p \langle 1 \rangle$ is canonically identified with a quotient of $(\O_F/\mathfrak q)^\times$ of size $p$.

\subsubsection{Reduction of dual Stark units at a Taylor--Wiles prime}

Let $\mathfrak Q$ be a prime of $L$ above a Taylor--Wiles prime $\mathfrak q$ of $F$. Let
$$\Frob_{\mathfrak Q} = \Frob_{\mathfrak Q/ \mathfrak q} \in G_{L/F} \subseteq G_{L/\Q}$$
be the Frobenius automorphism associated to the prime $\mathfrak Q$ above $\mathfrak q$.

\begin{lemma}\label{lemma:modp_pairing}
	For any Artin representation $\varrho_0 \colon G_{L/\Q} \to \GL(M_0)$ where $M_0$ is an $\O_{E}$-module, there is a natural pairing
	\begin{align*}
		(U_L[\varrho_0] \otimes k) \times (M_0^{\Frob_{\mathfrak Q}} \otimes k) & \to k\langle 1 \rangle \\
		(\varphi, m) & \mapsto \textnormal{reduction of }\varphi(m).
	\end{align*}
\end{lemma}
\begin{proof}
	For $\varphi \in U_L[\varrho_0]$ and $m \in M_0^{\Frob_{\mathfrak Q}}$, we have
	$$\varphi(m) \in (U_L \otimes k)^{\Frob_{\mathfrak Q}}.$$
	The composition
	$$U_L \hookrightarrow U_{L_{\mathfrak Q}} \twoheadrightarrow U_{L_{\mathfrak Q}}/ (1 + \mathfrak Q) \iso \mathbb F_{\mathfrak Q}^\times$$
	induces a reduction map
	$$(U_L \otimes k)^{\Frob_{\mathfrak Q}} \to (\mathbb F^\times_{\mathfrak Q} \otimes k)^{\Frob_{\mathfrak Q}} \iso k\langle 1 \rangle,$$
	where we recall that $ k\langle 1 \rangle = k \otimes (\O_F/\mathfrak q)_{p^n}^\times$.
\end{proof}

\begin{remark}
	We think of the {\em reduction map} as a {\em discrete logarithm}. Then this lemma is the discrete analogue of Lemma~\ref{lemma:pairing}, where the actual logarithm will be used. To generalize this result $p$-adically, one would use a $p$-adic logarithm.
\end{remark}

\begin{proposition}\label{prop:dual_units_mod_p}
	Let $\varrho \colon G' = G_{L/F} \to \GL(M)$ be the Artin representation associated to~$f$. Recall the notation $G = G_{L/\Q}$. Then there is a natural map:
	$$\theta^\vee_{\mathfrak q} \colon \bigoplus_{\sigma \in G/G'} (\Ad^0 M \otimes k)^{\Frob_{\sigma \mathfrak Q/\sigma \mathfrak q}} \otimes k \langle -1 \rangle \to  U_f^\vee \otimes k$$
	where the domain is a direct sum of free $k$-modules of rank 1.
\end{proposition}

We will later use the shorthand $U_{f,\sigma}^\vee = (\Ad^0 M \otimes k)^{\Frob_{\sigma \mathfrak Q/\sigma \mathfrak q}} \otimes k \langle -1 \rangle$. In the notation of the introduction, $U_{f, \sigma_i}^\vee = U_{f, i}^{\mathfrak p^n}$ if we label the representatives of $G/G'$ by $\sigma_1, \ldots, \sigma_d$.

\begin{proof}
	Applying Lemma~\ref{lemma:modp_pairing} to $\varrho_0 = \Ind_{G'}^G \Ad^0 \varrho$, we see that there is a pairing:
	$$(U_f \otimes k) \times (M_0^{\Frob_{\mathfrak Q}} \otimes k \langle - 1 \rangle) \to k,$$
	which induces a map
	$$(M_0^{\Frob_{\mathfrak Q}} \otimes k \langle - 1 \rangle) \to (U_f^\vee \otimes k).$$
	Then
	\begin{align*}
		M_0^{\Frob_{\mathfrak Q}} & = (\Ind_{G'}^{G} \Ad^0 M)^{\Frob_{\mathfrak Q}} \\
		& = \left(\bigoplus_{\sigma \in G/G'} \sigma \Ad^0 M\right)^{\Frob_{\mathfrak Q}} \\
		& = \bigoplus_{\sigma \in G/G'} (\Ad^0 M)^{\Frob_{\sigma \mathfrak Q/\sigma \mathfrak q}}
	\end{align*}
	because $\sigma \Frob_{\mathfrak Q/\mathfrak q} \sigma^{-1} = \Frob_{\sigma \mathfrak Q/\sigma \mathfrak q} \in G'$.
	
	Finally, using the basis such that $\varrho(\Frob_{\mathfrak Q}) = \begin{pmatrix}
		\alpha & \\
		& \beta
	\end{pmatrix}$ for $\alpha \neq \beta$, we have that
	$$\Ad^0 \varrho( \Frob_{\mathfrak Q}) = \begin{pmatrix}
		\frac{\alpha}{\beta} \\
		& \frac{\beta}{\alpha} \\
		& & 1
	\end{pmatrix}.$$
	Since $\alpha \neq \beta$, this shows that $(\Ad^0 M)^{\Frob_{\sigma \mathfrak Q/\sigma q}}$ has rank 1.
\end{proof}

We finally recast this in the language of~\cite[Section 2.9]{Harris_Venkatesh}. For any $\mathfrak Q$, we may consider the element
$$e_{\mathfrak Q} = \varrho(\mathrm{Frob}_{\mathfrak Q}) - (1/2) \mathrm{Tr}  \varrho(\mathrm{Frob}_{\mathfrak Q}) \in \Ad^0 \varrho.$$
Note that for all $g \in G_{L/F}$,
$$e_{g \mathfrak Q} = \Ad(\varrho(g)) e_{\mathfrak Q}.$$
Therefore:
$$\Ad^0(\Frob_{\mathfrak Q}) e_{\mathfrak Q} = e_{\Frob_{\mathfrak Q} \mathfrak Q} = e_{\mathfrak Q},$$
showing that
$$e_{\mathfrak Q} \in (\Ad^0 \varrho)^{\Frob_{\mathfrak Q}}.$$

By Proposition~\ref{prop:dual_units_mod_p}, this choice defines a map
\begin{equation}\label{eqn:dual_units_mod_p}
	\theta^\vee_{\mathfrak q} \colon \bigoplus_{\sigma \in G/G'} k\langle -1 \rangle \to U_f^\vee \otimes k.
\end{equation}
When $F = \Q$, this recovers the map $\theta_{q}^\vee$ from \cite[Section 2.9]{Harris_Venkatesh}.

\subsection{The Shimura class}\label{subsec:Shimura_class}

We consider two level structures: for an ideal $\mathfrak N \subseteq \O_F$,
\begin{align*}	\Gamma_0(\frak N) & = \left\{\begin{pmatrix}
		a & b \\
		c & d
	\end{pmatrix} \in \GL(\O_F \oplus \mathfrak D^{-1})\ \middle|\ c \in \frak N \right\}, \\
	\Gamma_{1}(\frak N) & = \left\{\begin{pmatrix}
		a & b \\
		c & d
	\end{pmatrix} \in \GL(\O_F \oplus \mathfrak D^{-1})\ \middle|\ c, a-1 \in \frak N \right\},
\end{align*}
where $\mathfrak D$ is the different ideal of $F$. Note that $\Gamma_{1}(\mathfrak N) \subseteq \Gamma_0(\mathfrak N)$ and the quotient is isomorphic to $(\O/\mathfrak N)^\times$. We let
$$X_0(\mathfrak N), X_{1}(\mathfrak N)  = \text{Hilbert modular variety with }\Gamma_0(\mathfrak N), \Gamma_{1}(\mathfrak N)\text{-level structure, respectively}.$$
For $\mathfrak N$ large enough, both of these are schemes over $\Z[1/N_{F/\Q} \mathfrak N]$ (c.f.\ \cite{Dimitrov_Gamma_1_n}) and they have good reduction modulo primes~$p$ not dividing $N_{F/\Q}\mathfrak N$. The covering
$$X_{1}(\mathfrak N) \to X_0(\mathfrak N)$$
descends to a covering
$$ X_{1}(\mathfrak N)_k \to X_0(\mathfrak N)_k$$
with Galois group $(\O/\mathfrak N)^\times$.

Let $q > 5$ be a prime and $\mathfrak q$ be a prime of $F$ above it. Then
$$X_{1}(\mathfrak q) \to X_0(\mathfrak q)$$
is a $(\O/\mathfrak q)^\times$-covering. We may pass to the unique subcovering with Galois group $\Delta = (\O/\mathfrak q)^\times_{p^n}$:
$$X_{1}(\mathfrak q)^{\Delta} \to X_0(\mathfrak q).$$
This extends to an \'etale covering of schemes over $\Z[1/q]$, and hence induces an \'etale covering
$$X_{1}(\mathfrak q)^{\Delta}_k \to X_0(\mathfrak q)_k$$
(c.f.\  \cite[Corollary 2.3]{Mazur} for $[F:\Q] = 1$ and \cite[Prop.\ 3.4]{Dimitrov_arithmetic_aspects_of_HMF} for $[F:\Q] > 1$; the assumption that $q > 5$ is needed to avoid elliptic points). 

We hence get a class
\begin{equation}
	\mathfrak S_k \in H^1_{\text{\'et}}(X_0(\mathfrak q)_k, k \langle 1 \rangle),
\end{equation}
where we recall that $k \langle 1 \rangle \iso k \otimes \Delta$. Using the natural map $k \to \mathbb G_a$ of \'etale sheaves over~$X_0(\mathfrak q)_k$, we obtain a class:
\begin{equation}
	\mathfrak S_{\mathbb G_a} \in H^1_{\text{\'et}}(X_0(\mathfrak q)_k, \mathbb G_a \langle 1 \rangle).
\end{equation}
Finally, using Zariski--\'etale comparison, we have an isomorphism:
$$H^1(X_0(\mathfrak q)_k, \O\langle 1 \rangle) \to H^1(X_0(\mathfrak q)_k, \mathbb G_a\langle 1 \rangle) $$
and hence $\mathfrak S_{\mathbb G_a}$ defines a class
\begin{equation}\label{eqn:Shimura_class}
	\mathfrak S \in H^1(X_0(\mathfrak q)_k, \O \langle 1 \rangle).
\end{equation}

\begin{definition}
	The {\em Shimura class} is the cohomology class $\mathfrak S \in H^1(X_0(\mathfrak q)_k, \O \langle 1 \rangle)$ obtained above~(\ref{eqn:Shimura_class}).
\end{definition}

We will use it next to construct a mod $\p^n$ derived Hecke operator.

\subsection{Construction of derived Hecke operators}\label{subsec:derived_Hecke}

Let $\mathfrak N$ be the level of $f$ and recall that we consider $X = X_{1}(\mathfrak N)$ over $\Z[1/N_{F/\Q}\mathfrak N]$.

Write $X_{0, 1}(\mathfrak q, \mathfrak N)$ for $X$ with added $\Gamma_0(\mathfrak q)$-level structure at $\mathfrak q$. This is a Hilbert modular variety for the group $\Gamma_1(\mathfrak q, \mathfrak N)$ in the notation of \cite{Dimitrov_Gamma_1_n}, and hence also has a smooth, projective, integral model. 

Then the Shimura class $\mathfrak S$ pulls back to a class
$$\mathfrak S_{X} \in H^1(X_{0, 1}(\mathfrak q, \mathfrak N)_k, \O \langle 1 \rangle).$$
Cupping with this class gives a map
\begin{equation}
	H^0(X_{0,1}(\mathfrak q, \mathfrak N)_k, \omega) \overset{\cup \mathfrak S_{X}}{\longrightarrow}  H^1(X_{0, 1}(\mathfrak q, \mathfrak N)_k, \omega)\langle 1 \rangle.
\end{equation}

Classically, Hecke operators are defined as operators on cohomology induced by certain correspondences:
\begin{center}
	\begin{tikzcd}
		& X_{0, 1}(\mathfrak q, \mathfrak N) \ar[ld, swap, "\pi_1"] \ar[rd, "\pi_2"] \\
		X & & X
	\end{tikzcd}
\end{center}

We define the \textit{derived Hecke operator} by the same push-pull procedure but cupping with~$\mathfrak S_{X}$ in the middle:
\begin{center}
	\begin{tikzcd}
		H^0(X_k, \omega) \ar[r, "\pi_1^\ast"] & H^0(X_{0,1}(\mathfrak q, \mathfrak N)_k, \omega) \ar[r, "\cup \mathfrak S_{X}"] & H^1(X_{0, 1}(\mathfrak q, \mathfrak N)_k, \omega)\langle 1 \rangle \ar[r, "\pi_{2, \ast}"] & H^1(X_k, \omega)\langle 1 \rangle.
	\end{tikzcd}
\end{center}

Finally, for any $z \in k \langle -1 \rangle$, we define
\begin{equation}
	T_{\mathfrak q, z} \colon H^0(X_k, \omega) \to H^1(X_k, \omega)
\end{equation}
by composing the above map with multiplication by $z$.

More generally, for each $z \in k \langle -1 \rangle$, there is an operator
\begin{equation}\label{eqn:derived_Hecke}
	T_{\mathfrak q, z} \colon H^q(X_k, \omega) \to H^{q + 1}(X_k, \omega),
\end{equation}
defined analogously.

Recall that equation~\eqref{eqn:dual_units_mod_p} defines a map:
$$\theta_{\mathfrak q}^\vee \colon \bigoplus_{\sigma \in G/G'} k\langle -1 \rangle \to U_f^\vee \otimes k.$$
We may hence define an action of the codomain on coherent cohomology of the special fiber as follows.

\begin{definition}\label{def:action_mod_p}
	For each $\sigma \in G/G'$ and $z \in k\langle - 1 \rangle$, we define the action of $z$ in the $\sigma$-component of $\bigoplus\limits_{\sigma \in G/G'} k\langle -1 \rangle$ by:
	$$T_{\sigma \mathfrak q, z} \colon H^\ast(X_k, \omega)_f \to H^{\ast + 1}(X_k, \omega)_f.$$
	This naturally extends to an action of $\bigwedge^\ast \bigoplus\limits_{\sigma \in G/G'} k\langle -1 \rangle$ on $H^\ast(X_k, \omega)_f$.
\end{definition}

\subsection{The conjecture}\label{subsec:conj_mod_p}

We conjecture there is an action of $U_f^\vee$ on the $f$-isotypic component of the cohomology space $H^\ast(X, \omega)_f$ which reduces modulo $\p^n$ to the action of the operators~$T_{\mathfrak q, z}$.

For $h \in H^\ast(X_{\O[1/N(\mathfrak N)]}, \omega)$, we write $\overline h \in H^\ast(X_{k}, \omega)$ for its reduction. Equation~\eqref{eqn:dual_units_mod_p} defines a map:
$$\theta^\vee_{\mathfrak q} \colon \bigoplus_{\sigma \in G/G'} k\langle -1 \rangle \to U_f^\vee \otimes k.$$
associated to a Taylor--Wiles primes $\mathfrak q$ of $F$ and a prime $\mathfrak Q$ above it. In Definition~\ref{def:action_mod_p}, we defined an action of the domain by derived Hecke operators. We conjecture that the resulting action of $ U_f^\vee \otimes k$ on the special fiber lifts to an integral action of $U_f^\vee$.

\begin{conjecture}\label{conj:modulo}
	There is an action $\star$ of the exterior algebra $\bigwedge^\ast(U_f^\vee)$ on $H^\ast(X_{\O[1/N(\mathfrak N)]}, \omega)_f$ such that the induced action of $\bigwedge^\ast(U_f^\vee) \otimes k$ on  $H^\ast(X_{\O[1/N(\mathfrak N)]}, \omega)_f \otimes k$ is the one described above. More specifically, fix a quadruple $(\mathfrak p, n, \sigma, \mathfrak q)$ with
	\begin{itemize}
		\item $\mathfrak p$ a prime of ${E}$ satisfying the above conditions,
		\item $n \geq 1$ an integer,
		\item $\sigma \in G/G'$,
		\item $q > 5$ a prime and $\mathfrak q$ a Taylor--Wiles primes of level $n$ above it; in particular $N \mathfrak q \equiv 1\ (p^n)$.
	\end{itemize}
	For an element $u^\vee \in U_f^\vee$, consider its reduction $\overline{u^\vee} \in U_f^\vee \otimes k$, and suppose that
	$$\overline{u^\vee} = \sum_{\sigma \in G/G'} \theta^\vee_{\mathfrak q}( z_\sigma) \qquad\text{for some $z_\sigma \in k\langle -1 \rangle$}.$$
	Then:
	$$ \alpha \cdot \overline{u^\vee \star \omega_f} = \sum_{\sigma \in G/G'} T_{\sigma \mathfrak q, z_\sigma} \overline{\omega_f}$$
	for some constant $\alpha$.
\end{conjecture}

\begin{remark}
	Harris and Venkatesh~\cite{Harris_Venkatesh} and Marcil~\cite{Marcil} provide numerical evidence for this conjecture for $F = \Q$ and~$n =1$. To do that, they first perform an {\em explication} (\cite[Section~5]{Harris_Venkatesh}), putting the conjecture in a more computable form. They relate it to a question about a pairing considered by Mazur \cite{Mazur} and then rely on a computation of this pairing due to Merel~\cite{Merel}. While the initial steps of the explication can be performed in our case, putting Conjecture~\ref{conj:modulo} in a similar framework, the analogue of Merel's computation is currently not available in the literature.
	
	In dihedral cases, the conjecture of Harris and Venkatesh has since been proved by Darmon--Harris--Rotger--Venkatesh~\cite{Darmon_Harris_Rotger_Venkatesh}.
\end{remark}

When $F = \Q$ and $n = 1$, Harris--Venkatesh~\cite[Section 4]{Harris_Venkatesh} prove the following result:
$$\text{vanishing of }T_{q, z} \overline{f} \quad \Longrightarrow\quad \text{vanishing of the map $\theta_q^\vee \colon k \langle -1 \rangle \to U_f^\vee \otimes k$},$$
assuming an ``$R = T$'' theorem. It would be interesting to obtain a similar result in our case. We expect that the rank $r$ of the map
$$\theta^\vee_{\mathfrak q}  \colon \bigoplus_{\sigma \in G/G'} k\langle -1 \rangle \to U_f^\vee \otimes k$$
from equation~\eqref{eqn:dual_units_mod_p} can be any number $0 \leq r \leq d$. Hence the strongest analogue of the above result should be:
$$\rank \langle T_{\sigma \mathfrak q, z} \overline{f} \ | \ \sigma \in G/G' \rangle = \rank(\theta^\vee_{\mathfrak q}).$$
A weaker version simply states:
$$\text{vanishing of }T_{\sigma \mathfrak q, z} \overline{f} \text{ for all }\sigma \in G/G' \quad \Longrightarrow\quad \text{vanishing of the map $\theta^\vee_{\mathfrak q}$}.$$

Note that the proof in the case $F = \Q$ relies on the approach of Calegari--Geraghty~\cite{Calegari_Geraghty} to modularity lifting. Since their results apply to general $F$, one could hope to prove the above results in a similar way, but we have not explored this further yet.

Since we expect that the map $\theta^\vee_{\mathfrak q}$ may sometimes have rank $d$, we want to make sure that we can produce a rank $d$ group of operators $T_{\mathfrak q, z}$ in order to pin down the conjectural action.

\begin{lemma}
	For any $\mathfrak p$ and $n$, there is a prime $q \equiv  1\ (p^n)$ which splits completely in $F$ and the primes $\mathfrak q_1, \ldots, \mathfrak q_d$ above $q$ are Taylor--Wiles primes for $f$ of level $n$.
\end{lemma}
\begin{proof}
	We first show that there exists a positive density of primes $q$ of $\Q$ which split completely in $F$ such that $q \equiv 1\ (p^n)$. Consider the field $F(\zeta_{p^n})$ for a primitive $p^n$th root of unity and a prime $q$ of $\Q$ in the field diagram:
	\begin{center}
		\begin{tikzcd} 
			& \mathfrak Q & F(\zeta_{p^n}) \\
			\mathfrak q \ar[ru, no head] &	F \ar[ru, no head] & & \Q(\zeta_{p^n}) \ar[lu, no head] \\
			& q \ar[lu, no head] & \Q \ar[lu, no head] \ar[ru,no head]
		\end{tikzcd}
	\end{center}
	Since we assume that $p$ has good reduction in $F$, the fields $\Q(\zeta_{p^n})$ and $F$ have disjoint ramification, and hence we have isomorphisms:
	\begin{center}
		\begin{tikzcd}
			G_{F(\zeta_{p^n})/F} \ar[r, "\iso"] & G_{\Q(\zeta_{p^n})/\Q} \ar[r, "\iso"] & (\Z/p^n \Z)^\times \\
			D(\mathfrak Q/\mathfrak q) \ar[r, "\iso"] \ar[u, hook] & D(\mathfrak Q \cap \Q(\zeta_{p^n})/ q) \ar[u, hook] \ar[r, "\iso"] & \langle q \rangle \ar[u, hook]
		\end{tikzcd}
	\end{center}
	via the restriction map. By Cheboratev density theorem, there is a positive density of primes~$q$ of $\Q$ that split completely in $F(\zeta_{p^n})$. These $q$ also split completely in $F$ and in $\Q(\zeta_{p^n})$ which shows that
	$$q \equiv 1 \mod p^n$$
	using the above diagram.
	
	Since there is a positive density of primes $q$ with the above property, there exists a positive density for which $\mathfrak q_1, \ldots, \mathfrak q_d$ are Taylor--Wiles primes for $f$ of level $n$.
\end{proof}

In this case, we have $d$ derived Hecke operators $T_{\mathfrak q_1, z_1}, \ldots, T_{\mathfrak q_d, z_d}$ and we expect that if they are linearly independent, then the map $\theta^\vee_{\mathfrak q}$ is an isomorphism.

\section{Archimedean realization of the motivic action}\label{sec:arch_real}

We continue using the notation of Section~\ref{subsec:Stark_for_HMF}: the Artin representation $\varrho_f$ associated to $f$ factors through a finite Galois extension $L/F$ and has coefficients in a number field $E$, i.e.\ $\varrho_f \colon \Gal(L/F) \to \GL_2(E)$.

Fix embeddings $\tau \colon L \hookrightarrow \C$ and $\iota \colon E \hookrightarrow \C$. We will describe:
\begin{itemize}
	\item an {\em isomorphism} 
	$$\theta^\vee_\C \colon  \bigoplus_{j=1}^d U_{f, j}^\C \overset\iso\to U_f^\vee \otimes_{\iota} \C$$
	for some one-dimensional spaces $U_{f, j}^\C$ in Proposition~\ref{prop:dual_units_over_C};
	
	\item an action of the codomain via {\em partial complex conjugation operators}: 
	\begin{align*}
		H^q(X_\C, \omega)_f & \to H^{q+1}(X_\C, \omega)_f \\
		\omega_f & \mapsto \omega_f^{\sigma_j}
	\end{align*}
	for a chosen element of $U_{f, j}^\C$ in Sections~\ref{subsec:complex_conj} and~\ref{subsec:action_over_C} (Definition~\ref{def:action});
\end{itemize}
and conjecture that the resulting action of $U_f^\vee \otimes E \subseteq U_f^\vee \otimes_{\iota} \C$ preserves the rational structure on coherent cohomology in Section~\ref{subsec:conj_C} (Conjecture~\ref{conj:motivic_action}).

\subsection{Partial complex conjugation and Harris'  periods}\label{subsec:complex_conj}

Following \cite{Shimura_HMF, Harris_periods_I}, we briefly recall the definition of {\em partial complex conjugation} operators on Hilbert modular forms. We encourage the reader to consult~\cite{Harris_periods_I, Harris_delbar, Su} for details.

Let $Y$ be an open Hilbert modular variety of level $\Gamma_1(\mathfrak N)$ and write $X$ for a smooth toroidal compactification of $Y$ defined over $\Q$. Associated to a weight $(\underline k, r)$ where $\underline k \in \Z^d$ and $k_j \equiv r \bmod 2$ for $1 \leq j \leq d$ is an automorphic sheaf $\mathcal E_{\underline k, r}$ over $Y$ whose sections are weight $(\underline k, r)$ Hilbert modular forms. We normalize the isomorphism between the sections $H^0(Y_{\C}, \mathcal E_{\underline k, r})$ and Hilbert modular forms of weight $(\underline k, r)$ so that Hilbert modular forms with Fourier coefficients in $E$ give sections of $H^0(Y, \mathcal E_{\underline k, r}) \otimes E$. In particular, this differs from Harris' normalization by a factor of $(2 \pi i)^{\frac{1}{2} (dr + \sum\limits_j k_j)}$; see~\cite[(1.6.4)]{Harris_periods_I}. For simplicity, we assume that $r \in \{0,1\}$ according to the parity of $k_j$. 

The automorphic sheaf $\mathcal E_{\underline k, r}$ can be extended to $X$ in two ways, denoted $\mathcal E_{\underline k, r}^{\mathrm{can}}$ and $\mathcal E_{\underline k, r}^{\mathrm{sub}}$. The cohomology of these sheaves is independent of the choice of toroidal compactification. Following Harris, we will be interested in the space:
$$H^q(X, \mathcal E_{\underline k, r}) =  \im(H^q(X, \mathcal E_{\underline k, r}^{\mathrm{sub}}) \to H^q(X, \mathcal E_{\underline k ,r}^{\mathrm{can}}))$$
which is a vector space over $F(\underline k) = F^{\Gamma(\underline k)}$ where $\Gamma(\underline k) = \{\sigma \in G_\Q \ | \ \underline k^\sigma = \underline k\}$.

Let $f$ be a normalized Hilbert modular eigenform $f$ of weight $(\underline k, r)$ and level $\Gamma_1(\mathfrak N)$ such that $T(\mathfrak p) f = a_{\mathfrak p} f$ and $a_{\mathfrak p} \in E_f$. Hecke operators act on the higher cohomology groups and we write:
\begin{equation}
	H^q(X, \mathcal E_{\underline k, r})_f = \{\omega \in H^q(X, \mathcal E_{\underline k, r}) \otimes E_f \ | \ T(\mathfrak p) \omega = a_{\mathfrak p}  \omega\}
\end{equation}
for the $f$-isotypic component under the action of the Hecke algebra.

For any subset $J$ of the infinite places $\Sigma_\infty = \{\sigma_1, \ldots, \sigma_d\}$ of $F$, we assume that there exists a unit~$\epsilon_J \in \O_F^\times$ such that
\begin{equation}
	\begin{cases}
		\sigma(\epsilon_J) > 0 & \text{if }\sigma \not\in J, \\
		\sigma(\epsilon_J) < 0 & \text{if }\sigma \in J.
	\end{cases}
\end{equation}
When $d = 2$, this amounts to the standard assumption (e.g., \cite{Oda}) that $\O_F$ has a fundamental unit of negative norm.

Given $f$ and a subset $J$ of $\Sigma_\infty$, we can apply complex conjugation to variables corresponding to places in~$J$:
\begin{equation}\label{eqn:f^J}
	f^J(\underline z) = f(\underline z^J) \cdot  \prod_{j \in J} \mathrm{Im}(z_j)^{k_j}
\end{equation}
where
\begin{equation}
	(\underline z^J)_j = \begin{cases}
		(\epsilon_J)_j z_j & \text{if }\sigma_j \not\in J, \\
		(\epsilon_J)_j \overline z_j & \text{if }\sigma_j \in J.
	\end{cases}
\end{equation}
This defines a $C^\infty$-function on $\H^d$ which has weight $-k_j$ at places $\sigma_j \in J$ and $k_j$ at places $\sigma_j \not\in J$. We can then define a Dolbeault class associated to $f$ and $J$:
\begin{equation}
	\omega_f^J = \left[f^J(\underline z) \cdot \bigwedge_{j \in J} \frac{dz_j \wedge d\overline{z_j}}{y_j^2}\right] \in H^{0, |J|}(Y_\C, \mathcal E_{\underline k(J), r}) \iso H^{|J|}(Y_\C, \mathcal E_{\underline k(J), r}),
\end{equation}
where:
$$\underline k(J) = \begin{cases}
	k_j & \sigma_j \not\in J, \\
	2 - k_j & \sigma_j \in J.
\end{cases}$$

\begin{remark}\label{rmk:cc_automorphically}
	This cohomology classes $\omega_f^J$ are independent of the choice of the unit $\epsilon_J$ above. Moreover, if $f$ corresponds to the automorphic function $\varphi$, the function $f^J$ defined in equation~\eqref{eqn:f^J} corresponds to the automorphic function:
	$$\varphi^J(g) = \varphi(g \cdot g_J)$$
	where $g_J \in \GL_2(F \otimes \R)$ is:
	$$(g_J)_j = \begin{cases}
		\begin{pmatrix}
			1 & 0 \\
			0 & -1
		\end{pmatrix} & \sigma_j \in J, \\
		\begin{pmatrix}
			1 & 0 \\
			0 & 1
		\end{pmatrix} & \sigma_j \in J.
	\end{cases}$$
	This gives a definition of partial complex conjugation even in cases where the unit $\epsilon_J$ does not exist. See~\cite[Section 1.4]{Harris_periods_I} for details.
\end{remark}

\begin{theorem}[Harris, Su]\label{thm:Harris_Su}
	\leavevmode
	\begin{enumerate}
		\item The cohomology classes $\omega_f^J$ extend to toroidal compactifications:
		$$\omega_f^J \in H^{|J|}(X_\C, \mathcal E_{\underline k(J), r})_f.$$
		
		\item Let $J \subseteq \Sigma_\infty$ be any subset. Then a basis of $H^{|J|}(X_\C, \mathcal E_{\underline k(J), r})_f $ is given by
		$$\{\omega_f^I \ | \ |I| = |J| \text{ and }\underline k(I) = \underline k(J) \}.$$
		In particular, if we write $J_1 = \{\sigma_j \in \Sigma_\infty \ | \ k_j = 1\}$, then:
		$$\dim  H^{|J|}(X_\C, \mathcal E_{\underline k(J), r})_f = \binom{|J_1|}{|J \cap J_1|}.$$
	\end{enumerate}
\end{theorem}
\begin{proof}
	For $k_j \geq 2$, see \cite[Lemmas 1.4.3, 2.4.5]{Harris_periods_I}. When $k_j = 1$ for some $j$, this follows from the main theorem of \cite{Su} and an analogous computation of $(\mathfrak P, K)$-cohomology.
\end{proof}

We are particularly interested in the case $(\underline k, r) = (\underline 1, 1)$. In this case, $\mathcal E_{\underline 1, 1}$ is identified with the Hodge bundle $\omega$ used in the previous section.

\begin{corollary}\label{cor:dim_of_coh_group}
	Suppose $(\underline k, r) = (\underline 1, 1)$. Then a basis of $H^j(X_\C, \mathcal E_{\underline 1, 1})_f$ is given by
	$$\{\omega_{f}^J \ | \ J\subseteq \Sigma_\infty \text{ and } |J| = j\}.$$
	In particular,
	$$\dim H^{j}(X_\C, \mathcal E_{\underline 1, 1})_f = \binom{d}{j}.$$
\end{corollary}

It is also important to note when the cohomology spaces are one-dimensional.

\begin{corollary}\label{cor:1-dimensional}
	For any $J \subseteq \Sigma_\infty$, $\dim H^{|J|}(X_\C, \mathcal E_{\underline k(J), r})_f > 1$ if and only if both $J$ and $\Sigma_\infty \setminus J$ contain a place at which $f$ has weight one.
\end{corollary}
\begin{proof}
	For the `if' implication, take $\sigma \in J \cap J_1$ and $\sigma' \in  (\Sigma_\infty \setminus J) \cap J_1$, and define
	$$J' = (J \setminus \{\sigma\}) \cup \{\sigma'\}.$$
	Then $|J'| = |J|$ and $\underline k(J') = \underline k(J)$, so $\omega_f^{J}, \omega_f^{J'} \in  H^{|J|}(X_\C, \mathcal E_{\underline k(J), r})_f$ are linearly independent. 
	
	Conversely, suppose $\dim H^{|J|}(X_\C, \mathcal E_{\underline k(J), r})_f > 1$. Then there exists $J' \neq J$ such that $\omega_f^{J'} \in H^{|J|}(X_\C, \mathcal E_{\underline k(J), r})_f$, i.e.\  $|J'| = |J|$ and $(J \cup J') \setminus (J \cap J') \subseteq J_1$. Then $\sigma \in J \setminus J'$ belongs to $J \cap J_1$ and $\sigma' \in J' \setminus J$ belongs to $(\Sigma_\infty \setminus J) \cap J_1$.
\end{proof}

This leads to the definition of Harris' period invariants when the cohomology space is one-dimensional.

\begin{lemma}[{\cite[Lemma~1.4.5]{Harris_periods_I}}]\label{lemma:def_of_nu}
	Let $J$ be a set of infinite places which contains either all or none of the weight one places of $f$. Then there is a number $\nu^J(f) \in \C^\times$, well-defined up to multiplication by elements in $E_f(J)^\times$ where $E_f(J) = E_f F(\underline k(J))$, such that
	$$ \frac{\omega_f^J}{\nu^J(f)} \in H^{|J|}(X, \mathcal E_{\underline k, r})_f \subseteq H^{|J|}(X_\C, \mathcal E_{\underline k, r})_f.$$
	Clearly, when $J = \emptyset$, we may take $\nu^J(f) = 1$.
\end{lemma}

\begin{definition}\label{def:nu_Harris}
	Let $J$ be a set of infinite places which contains either all or none of the weight one places of $f$. Then the complex number $\nu^J(f)$ defined by Lemma~\ref{lemma:def_of_nu} is the {\em period} or {\em period invariant} associated to~$f$ and $J$. It is well-defined up to $E_f(J)^\times$.
\end{definition}

\begin{remark}
	Despite of the difference in trivializations of the line bundles, the above period invariants $\nu^J(f)$ agree with Harris' period invariants $\nu^J(\pi_f)$, where $\pi_f$ is the automorphic representation associated to $f$. For example, note that both of the normalizations result in $\nu^\emptyset(f) = \nu^\emptyset(\pi_f) = 1$. 
\end{remark}

Shimura defines periods by considering Petersson inner products on Shimura varieties associated to quaternion algebras over $F$. Harris' definition is much less explicit, but it is related to Petersson inner products as follows.

\begin{proposition}[{\cite[Prop.~1.5.6]{Harris_periods_I}}]\label{prop:fact_of_Pet}
	For any $J \subseteq \Sigma_\infty$, we have that:
	\begin{equation}
		\nu^J(f) \cdot \nu^{\Sigma_\infty \setminus J}(f^\varrho) \sim_{ E(J)^\times} \langle f, f \rangle
	\end{equation}
	where $f^\varrho(z) = \overline{f(-\overline z)}$ is Shimura's complex conjugation, and
	\begin{equation}
		\langle f, g \rangle = \int\limits_{\Gamma \backslash \H^d} f(\underline z) \overline{g(\underline z)} \prod_{j=1}^d y_j^{k_j} \frac{dz_j \wedge d\overline{z_j}}{y_j^2}.
	\end{equation}
\end{proposition}

Therefore, we may think of $\nu^J(f)$ as a certain factor of the Petersson inner product $\langle f, f \rangle$.

\begin{remark}\label{rmk:Shimura_norm}
	Here and elsewhere we use the above normalization of Petersson inner products. This is consistent with  \cite{Hida_GL2xGL2, Hida_Tilouine}, which we refer to later. This differs from Shimura's normalization of Petersson inner products~\cite[(2.27, 2.28)]{Shimura_HMF}:
	$$\langle f, g \rangle_{\text{Shimura}} = \frac{1}{\mu(\Gamma \backslash \H^d)} \langle f, g \rangle,$$
	where $\mu(\Gamma \backslash \H^d)$ is the volume of the fundamental domain. It also differs from Harris' normalization, since:
	\begin{equation}
		\langle f, f \rangle_{\mathrm{Harris}} \sim_{\Q^\times} (2\pi i)^{-dr} \langle f, f\rangle_{\text{Shimura}}.
	\end{equation}
	\cite[(1.6.3)]{Harris_periods_I}.
\end{remark}

\begin{remark}
	The proof in loc.\ cit.\ is based on the rationality of (a Tate twist of) the Serre duality pairing \cite[(1.5.4)]{Harris_periods_I}: 
	\begin{equation}\label{eqn:SD}
		\cup \colon H^{|J|}(X, \mathcal E_{\underline k(J), r})_f \times H^{|\Sigma_\infty \setminus J|}(X, \mathcal E_{\underline k(\Sigma_\infty \setminus J), r})_{f^\varrho} \to E(J)
	\end{equation}
	induced by the cup product, and the identity~\cite[(1.5.5.2)]{Harris_periods_I}:
	\begin{equation}\label{eqn:Petersson}
		\omega_f^{J} \cup \omega_{f^\varrho}^{\Sigma_\infty \setminus J} = \pm \langle f, f \rangle.
	\end{equation}
\end{remark}


\begin{remark}\label{rmk:rel_to_other_periods}
	In this extended remark, we discuss the relation of Harris' periods to other periods attached to Hilbert modular forms. 
	The study of period invariants was initiated by Shimura \cite{Shimura_alg_relations_between, Shimura:on_crit_val}, who studied the case when the weights at all places are at least two. 
	In this case, Shimura conjectured the existence of a set of period invariants $c_\sigma$, one attached to each infinite place $\sigma$ of $F$;
	moreover, he conjectured that if $B$ is any quaternion algebra over $F$ such that $f$ transfers to a form $f_B$ on $B^\times$, then the Petersson norm of $f_B$ (if $f_B$ is chosen to be algebraic) is 	
	essentially a product of some of the $c_\sigma$ up to algebraic factors. More precisely,
	defining 
	\[
	q_B (f) :=  \langle f_B, f_B \rangle,
	\]
	Shimura conjectured that 
	\begin{equation}\label{eqn:Shimura_fact}
		q_B (f) \sim_{\overline \Q^\times} \prod_{\sigma \in \Sigma_{B,\infty}} c_\sigma,
	\end{equation}
	where $\Sigma_{B,\infty}$ is the set of infinite places where $B$ is split. 
	This conjecture was proved by Harris \cite{Harris:JAMS}, using the theta correspondence for unitary groups. 
	In this work, the periods~$c_\sigma$ are essentially {\it defined} as suitable ratios of periods on 
	quaternion algebras. The fact that the definition of the periods does not depend on choices of quaternion algebras 
	boils down to proving relations between periods on different quaternion algebras, which provides the 
	main thread of Harris' argument. This work admits an integral refinement which is studied in the ongoing work of Ichino--Prasanna (for example, \cite{Ichino_Prasanna}).
	
	In related work~\cite{Harris_periods_I, Harris_periods_II}, Harris gave another definition of such period invariants using rational structures on coherent cohomology. This is what was recalled in Definition~\ref{def:nu_Harris}. The advantage of this definition is that it does not require working with quaternion algebras;
	rather everything happens on the Hilbert modular variety attached to the group $\GL_{2,F}$. This also makes it easy 
	to see the relations between these periods and the transcendental factors of Rankin--Selberg and triple 
	product $L$-functions attached to two (respectively, three) Hilbert modular forms. 
	
	The point of our work is to define periods attached to parallel weight one forms, and relate them to rational structures on coherent cohomology. 
	For dimension reasons, one cannot simply use these rational structures directly to define periods. Indeed, the proof of Lemma~\ref{lemma:def_of_nu} relies on higher cohomology groups being one-dimensional whereas the dimensions are greater than one for weight one forms (c.f.\ Corollary~\ref{cor:dim_of_coh_group}). Instead, we give an ad hoc definition using logarithms of units, and conjecture (Conjecture~\ref{conj:Harris_periods}) a relationship to rational structures.
\end{remark}

\subsection{The action}\label{subsec:action_over_C}

To define the action of $U_f^\vee \otimes \C$ on coherent cohomology via partial complex conjugation operators, we first give an identification of this group with the trace zero adjoint representation of $f$.

\begin{lemma}\label{lemma:pairing}
	For any Artin representation $\varrho_0 \colon G_{L/\Q}\to \GL(M_0)$ where $M_0$ is an ${E}$-vector space, there is a natural perfect pairing
	\begin{align*}
		(U_L[\varrho_0] \otimes_{\iota} \C) \times (M_0^{c_0} \otimes_{\iota} \C) & \to \C \\
		(\varphi, m) & \mapsto \log(|(\tau \otimes \iota)(\varphi(m))|)
	\end{align*}
	which induces an isomorphism
	$$U_L[\varrho_0]^\vee \otimes \C \overset \iso \to M_0^{c_0} \otimes \C.$$
\end{lemma}
\begin{proof}
	This is a paraphrase of Proposition~\ref{prop:Stark_basis}.
\end{proof}

\begin{proposition}\label{prop:dual_units_over_C}
	Let $\varrho \colon G' = G_{L/F} \to \GL(M)$ be the Artin representation associated to a Hilbert modular newform of parallel weight one. We then have an isomorphism:
	\begin{equation*}
		\theta^\vee_\C \colon \bigoplus_{j=1}^d (\Ad^0 M \otimes_{\iota} \C)^{c_j} \overset\iso\to U_f^\vee \otimes_{\iota} \C.
	\end{equation*}
	For each $j$, consider the element $m_{1,j}$ in $(\Ad^0 M)^{c_j}$ as in Corollary~\ref{cor:Stark_basis_HMF} and let $\{\varphi_j\}$ be the corresponding basis of $U_f \otimes E$. Finally, let $\{u_j^\vee\}$ be the dual basis of $U_f^\vee \otimes E$. Then the matrix of the map $\theta_\C^\vee$ in these bases is the Stark regulator matrix $R_f = ( \log|u_{jk}  |)_{j, k}$ (c.f.\ Definition~\ref{def:Stark_reg_HMF}).
\end{proposition}


\begin{proof}
	The result is obtained by applying Lemma~\ref{lemma:pairing} to $\varrho_0 = \Ind_{G'}^G \Ad^0  \varrho$ and recalling that $M_0^{c_0} \iso \bigoplus\limits_{\sigma \in G/G'} \left(\Ad^0 M   \right)^{\sigma c_0\sigma^{-1}}$ by the proof of Corollary~\ref{cor:Stark_basis_HMF}. The explicit description of the map is given by the second part of Corollary~\ref{cor:Stark_basis_HMF}.
\end{proof}

\begin{remark}
	Note that both $U_f^\vee \otimes_{\iota} \C$ and $(\Ad^0 M \otimes_{\iota} \C)^{\sigma c_0 \sigma^{-1}}$ have natural ${E}$-rational structures $U_f^\vee \otimes {E_f}$ and $(\Ad^0 M \otimes {E_f})^{\sigma c_0 \sigma^{-1}}$ but the above isomorphism does not respect them. The rational structures differ by the Stark regulator matrix.
\end{remark}

\begin{definition}\label{def:action}
	We define the action of $\bigoplus\limits_{j=1}^d (\Ad^0 M \otimes_{\iota} \C)^{c_j}$ on $H^\ast(X_\C, \mathcal E_{\underline 1, 1})_f$ by letting $m_{1,j}$ act by 
	\begin{align*}
		H^j(X_\C, \mathcal E_{\underline 1,1})_f & \to H^{j+1}(X_\C, \mathcal E_{\underline 1,1})_f \\
		\omega_f^J & \mapsto \begin{cases}
			\omega^{J \cup \{\sigma_j\}} & \sigma_j \not\in J \\
			0 & \sigma_j \in J.
		\end{cases}
	\end{align*}
	This defines a graded action of $\bigwedge^\ast \bigoplus\limits_{j=1}^d (\Ad^0 M \otimes_{\iota} \C)^{\sigma_j c_0 \sigma_j^{-1}}$ on $H^\ast(X_\C, \mathcal E_{\underline 1, 1})_f$ such that $H^\ast(X_\C, \mathcal E_{\underline 1, 1})_f$  is generated in degree 0 by $f \in H^0(X, \omega)_f$.
\end{definition}

\subsection{The conjectures}\label{subsec:conj_C}

Recall that Proposition~\ref{prop:dual_units_over_C} defined an isomorphism:
\begin{equation}\label{eqn:dual_units_map}
	\theta^\vee_\C \colon \bigoplus_{j=1}^d (\Ad^0 M \otimes_{\iota} \C)^{c_j}  \overset\iso\to  U_f^\vee \otimes_{\iota} \C
\end{equation}
and Definition~\ref{def:action} described an action of the latter group on coherent cohomology. We conjecture that the resulting action of $U_f^\vee \otimes E$ is rational.

\begin{conjecture}\label{conj:motivic_action}
	Fix embeddings $\tau \colon L \to \C$ and $\iota \colon E \to \C$. Then the action of $U_f^\vee \otimes E \subseteq U_f^\vee \otimes_{\iota} \C$ on $H^\ast(X_\C, \mathcal E_{\underline 1, 1})_f$ via equation~\eqref{eqn:dual_units_map} and Definition~\ref{def:action} preserves the rational structure $H^\ast(X, \mathcal E_{\underline 1, 1})_f \otimes_{E_f} E$.
\end{conjecture}

This is the analogue of the main conjecture of Prasanna--Venkatesh~\cite{Prasanna_Venkatesh}. In Appendix~\ref{app:Prasanna-Venkatesh}, we discuss the specific relation to their conjecture and justify why the definition of the action is natural.

Next, we give a more explicit statement of rationality of cohomology classes, via Propositon~\ref{prop:dual_units_over_C}.

\begin{conjecture}\label{conj:Harris_periods}
	Let $A = (a_{ij}) = R_f^{-1}$ be the inverse of the Stark regulator matrix. Then for $j = 1, \ldots, d$, the cohomology classes
	$$u_i^\vee \star f = \sum_{i=1}^n a_{ij} \omega_f^{\sigma_i} \in H^1(X_\C, \mathcal E_{\underline 1, 1})_f$$
	belong to the rational subspace $ H^1(X, \mathcal E_{\underline 1, 1})_f \otimes E$. More generally, the rational cohomology classes in $H^j(X_\C, \mathcal E_{\underline 1,1})_f$ are given by the entries of the vector:
	$$\left( \textstyle
	\bigwedge^j A \right) \begin{pmatrix}
		\omega_f^{J_1} \\
		\vdots\\
		\omega_f^{J_{\binom{d}{j}}}
	\end{pmatrix} $$
	where $J_1, \ldots, J_{\binom{d}{j}}$ are the subsets of $\Sigma_\infty$ of order $j$. In particular, the cohomology class
	$$\frac{\omega_f^{\Sigma_\infty}}{\det R_f} \in H^d(X_\C, \mathcal E_{\underline 1,1})_f$$
	is rational.
\end{conjecture}

The final statement is equivalent to Stark's conjecture~\ref{conj:Stark} for $\Ad^0 \varrho_f$ (Theorem~\ref{thm:Stark_implies_conj}). Therefore, this conjecture may be interpreted as a refinement of Stark's conjecture in this case.

\begin{remark}
	A previous version of this manuscript incorrectly assumed that the Stark regulator matrix $R_f$ is diagonal, which lead to a different rationality statement. 
\end{remark}

\begin{example}[$d=1$]
	Suppose $d = 1$, i.e.\ $f$ is a modular form of weight one. Then the conjecture simply asserts that:
	\begin{equation}
		\frac{\omega_f^{\infty}}{\log|\tau (u_f)|} \in H^1(X, \mathcal E_{\underline 1, 1}) \otimes E
	\end{equation}
	where $u_f \in U_L$ is a unit associated to $f$. As far as we know, this conjecture is new in this case. It gives an archimedean analogue of the main conjecture of Harris--Venkatesh~\cite{Harris_Venkatesh}. As we will see (Corollary~\ref{cor:conj_for_F=Q}), it is equivalent to Stark's conjecture~\ref{conj:Stark} for $\Ad^0 \varrho_f$, and hence is true when the Fourier coefficients of $f$ are rational or when $f$ has CM.
\end{example}

\begin{example}[$d=2$]\label{ex:d=2}
	Suppose $d = 2$, i.e.\ $f$ is a Hilbert modular form of parallel weight one for a real quadratic field $F$. Then there are four units $u_{11}, u_{12}, u_{21}, u_{22} \in U_L$ associated to $f$ and
	\begin{equation}
		R_f = \begin{pmatrix}
			\log|\tau (u_{11}) | & \log|\tau(u_{12})| \\
			\log|\tau(u_{21})| & \log|\tau(u_{22}) |
		\end{pmatrix}.
	\end{equation}
	Its inverse is:
	\begin{equation}
		A = \frac{1}{\det R_f}
		\begin{pmatrix}
			\log|\tau (u_{22}) | & -\log|\tau(u_{12})| \\
			-\log|\tau(u_{21})| & \log|\tau(u_{11}) |
		\end{pmatrix}.
	\end{equation}
	Therefore, the rational classes in $H^1(X_\C, \mathcal E_{\underline 1, 1})_f$ should be:
	\begin{align}
		u_1^\vee \star f & = \frac{ \log|\tau(u_{22}) | \cdot \omega_f^1 - \log|\tau (u_{21}) | \cdot \omega_f^2}{\det R_f} \in H^1(X, \mathcal E_{\underline 1,1})_f \otimes E, \\
		u_2^\vee \star f & = \frac{-\log|\tau (u_{12}) | \cdot \omega_f^1 +\log|\tau (u_{11}) |\cdot \omega_f^2}{\det R_f} \in H^1(X, \mathcal E_{\underline 1,1})_f \otimes E.
	\end{align}
	We will give the following evidence for this:
	\begin{enumerate}
		\item the determinant of this basis of $H^1(X_\C, \mathcal E_{\underline 1, 1})_f$ is rational, assuming Stark's conjecture~\ref{conj:Stark} (Section~\ref{sec:Stark_conj}),
		\item in base change cases, we give numerical evidence that the restrictions of these cohomology classes to an embedded modular curve is rational (Section~\ref{sec:BC}).
	\end{enumerate}
	Finally, we expect the following class in $H^2(X_\C, \mathcal E_{\underline 1, 1})_f$ to be rational:
	\begin{equation}
		 (u_1^\vee \wedge u_2^\vee) \star f =  \frac{\omega_f^{\sigma_1, \sigma_2}}{\det R_f} \in H^2(X, \mathcal E_{\underline 1, 1})_f \otimes E.
	\end{equation}
	We prove this assertion in Corollary~\ref{cor:Stark_implies_conj}.
\end{example}

The goal of the next two sections is to present our evidence for Conjecture~\ref{conj:motivic_action}.

\section{Evidence: Stark conjecture }\label{sec:Stark_conj}

In this section, we present the theoretical evidence for Conjecture~\ref{conj:Harris_periods}. These follow from results of Stark and Tate presented in Section~\ref{subsec:Stark_conj}.

\subsection{Action of top degree elements}

We show that Stark's conjecture~\ref{conj:Stark} for $\Ad^0 \varrho_f$ is equivalent to the following consequence of Conjecture~\ref{conj:Harris_periods}. In particular, Theorem~\ref{thm:Stark} implies this consequence when $f$ has rational Fourier coefficients.

\begin{theorem}\label{thm:Stark_implies_conj}
	Let $f$ be a parallel weight one Hilbert modular form and $\varrho_f$ be the associated Artin representation. Stark's conjecture~\ref{conj:Stark} for $\Ad^0 \varrho_f$ is equivalent to the statement:
	\begin{equation}\label{eqn:Petersson_is_detR_f}
		\langle f, f \rangle \sim_{E^\times} f_{\varrho, 2}^{1/2} \det R_f,
	\end{equation}
	where $f_{\varrho, 2} = 2^{a(\varrho, 2)}$ is the Artin conductor at $p = 2$ of the trace 0 adjoint representation. In particular, equation~\eqref{eqn:Petersson_is_detR_f}  is true unconditionally if $f$ has rational Fourier coefficients.
\end{theorem}

\begin{remark}
	We expect that the factor $f_{\varrho,2}^{1/2}$ is rational; see Remark~\ref{rmk:adjoint_conductor} for more details. If we could prove this, we could remove ``up to a possible factor of $\sqrt{2}$'' in the corollaries below.
\end{remark}

Before presenting the proof of Theorem~\ref{thm:Stark_implies_conj}, we give two corollaries.

\begin{corollary}\label{cor:Stark_implies_conj}
	Stark's conjecture~\ref{conj:Stark} for the Artin representation $\Ad^0 \varrho_f$ is equivalent to the assertion that top degree elements, i.e.\ elements in $\bigwedge^d U_f^\vee \otimes E$, act rationally, up to a possible factor of~$\sqrt{2}$.  In particular, the latter is true if $f$ has rational Fourier coefficients.
\end{corollary}
\begin{proof}
	Recall from Conjecture~\ref{conj:Harris_periods} that top degree elements act by
	$$f \mapsto \frac{\omega_f^{\Sigma_\infty}}{\det R_f}.$$
	Then:
	$$ \left\langle f^\varrho, \frac{\omega_f^{\Sigma_\infty}}{\det R_f} \right\rangle_{\mathrm{SD}} = \frac{\langle f, f \rangle }{\det R_f}.$$
	Since $H^d(X, \mathcal E_{\underline 1, 1})_f$ is one-dimensional and the Serre duality pairing is rational, the rationality of~$\frac{\omega_f^{\Sigma_\infty}}{\det R_f}$ is equivalent to equation~\eqref{eqn:Petersson_is_detR_f}.
\end{proof}

\begin{corollary}\label{cor:conj_for_F=Q}
	Conjecture~\ref{conj:motivic_action} is equivalent to Stark's conjecture~\ref{conj:Stark} for $\Ad^0 \varrho_f$ when $F = \Q$, up to a possible factor of~$\sqrt{2}$. Hence Conjecture~\ref{conj:motivic_action} is true unconditionally when $f$ has rational Fourier coefficients or complex multiplication.
\end{corollary}

\begin{remark}
	We checked computationally (using the method of Collins~\cite{Collins}) that for a few modular forms $f$ of weight one from Example~\ref{ex:cubic_fields}, we have that $\langle f, f \rangle = 3 \log(|\iota(u_f)|)$. This was already observed by Stark~\cite[pp.\ 91]{StarkII}.
\end{remark}

The proof of Theorem~\ref{thm:Stark_implies_conj} requires 2 steps:
\begin{enumerate}
	\item relating $L(1, \Ad^0 \varrho_f)$ to $\langle f, f \rangle$,
	\item showing that $f_\varrho$ is a square when $\varrho = \Ind_{G_F}^{G_\Q} \Ad^0(\varrho_f)$, so that $f_\varrho^{1/2} \in \Q^\times$ (away from 2).
\end{enumerate}
We will then conclude Theorem~\ref{thm:Stark_implies_conj} from Proposition~\ref{prop:Stark_for_HMF}.

The relation of the adjoint $L$-value to the Petersson inner product was first observed by Hida, based on the work of Shimura~\cite{Shimura:special_values}. He also related the prime factors of the quotient $\frac{L(1, \Ad(f))}{\langle f,  f \rangle}$ to congruence primes of the modular form $f$~\cite{Hida:Cong_of_cusp_forms, Hida:Cong_of_cusp_forms, Hida:modules_of_congruence}. This work was later generalized to Hilbert modular forms~\cite{Hida_GL2xGL2, Hida_Tilouine, Ghate_adjoint_L-values}. An integral refinement of Conjecture~\ref{conj:Harris_periods} would hence have to account for congruence primes.

\begin{theorem}[{\cite[Theorem~7.1]{Hida_Tilouine}}]\label{thm:adj_and_pet}
	Let $f$ is a primitive Hilbert modular form of weight $(\underline k, r)$, level $\mathfrak N$. Then
	$$\langle f, f \rangle = |D_F|^{m -1 } \Gamma_F(k) N_{F/\Q} (\mathfrak N) 2^{-2\{\underline k\} +1} \pi^{-d - \{k\}} L_S(1,  f, \mathrm{Ad}),$$
	where
	$$L_S(s, f, \mathrm{Ad}) = \prod_{\mathfrak q \in S} L_{\mathfrak q} (N_{F/\Q}(\mathfrak q)^{-s}) L(s, f, \Ad),$$
	$S$ is a set of bad places, $L_{\mathfrak q}(N_{F/\Q}(\mathfrak q)^{-s})$ are bad local factors, $\{k\} = \sum\limits_j k_j$, and $m$ is an explicit integer which accounts for Hida's unitarization~\cite[(4.2a), (7.1)]{Hida_GL2xGL2}.
\end{theorem}

For an automorphic proof relating $L(1, \Ad(f))$ to $\langle f, f \rangle$, see~\cite[Prop.~6.6]{Ichino_Prasanna}. 

For parallel weight one Hilbert modular forms, this specializes to the following result we will use.

\begin{corollary}\label{prop:adjoint_and_Pet}
	Suppose $(\underline k, r) = (\underline 1, 1)$. Then:
	$$\langle f, f \rangle \sim_{E^\times}  \pi^{-2d} L(1, f, \Ad).$$
\end{corollary}

To finish the proof of Theorem~\ref{thm:Stark_implies_conj}, we need to check that $f_\varrho$ is a square (away from $p = 2$).

\begin{proposition}\label{prop:adjoint_conductor}
	Let $\pi_v$ be the local representation of $\GL_2(F_v)$ associated to $f$ at a finite place~$v$ of $F$. When $v$ lies above $2$, assume that $\pi_v$ is not a theta lift from a ramified quadratic extension. Then the adjoint conductor of $\pi_v$ is a square.
\end{proposition}
\begin{proof}
	It is enough to prove that the analytic conductors of the Rankin--Selberg $L$-functions $L(\pi_v \otimes \pi_v^\vee, s)$ are squares. When $\pi_v$ is not supercuspidal, Jacquet's results~\cite{Jacquet} give explicit formulas for the local conductors (see, for example, \cite[Section 4.2]{Collins}) and they are visibly squares.
	
	We hence just need to show the conductor is a square at places $v$ where $\pi_v$ is supercuspidal. Suppose throughout the rest of the proof that $F$ is a finite extension of $\Q_p$ and $\pi$ is a supercuspidal representation of $\GL(2, F)$.  We write $a(-)$ for the valuation of the conductor of a representation and prove that $a(\pi \times \pi^\vee)$ is even. 
	
	Since $\pi$ is supercuspidal, it is a theta lift of a character $\xi$ of a quadratic extension $K/F$ \cite[Theorem 7.4]{Gelbart}. Then:
	\begin{equation}\label{eqn:adjoint_conductor}
		a(\pi \times \pi^\vee) = 2 v_F(d_{K/F}) + f_{K/F} \cdot a(\xi (\xi^\varrho)^{-1})
	\end{equation}
	where $d_{K/F}$ is the discriminant of $K/F$, $f_{K/F}$ is the residue degree of $K/F$, and $\varrho$ is the non-trivial element of $\Gal(K/F)$. Indeed, if $\varrho$ is the Galois representation corresponding to $\pi$ via the local Langlands correspondence, then $\varrho = \Ind_{K}^F(\chi)$ where $\chi$ corresponds to $\xi$ via class field theory, and hence
	\begin{align*}
		a(\pi \times \pi^\vee) & = a(\varrho \otimes \varrho^\vee) \\
		& = a(\Ind_{K}^F \chi \otimes \Ind_K^F \chi^{-1}) \\
		& = a(\Ind_{K}^F \mathbbm{1} \oplus \Ind_K^F \chi (\chi^\varrho)^{-1}) \\
		& = a(\Ind_{K}^F \mathbbm{1}) + a(\Ind_K^F \chi (\chi^\varrho)^{-1}) \\
		& = 2v_F(d_{K/F}) + f_{K/F} \cdot a(\chi (\chi^\varrho)^{-1}) & \textnormal{\cite[pp.\ 101]{Serre:Local_fields}} \\
		& = 2v_F(d_{K/F}) + f_{K/F} \cdot a(\xi (\xi^\varrho)^{-1}).
	\end{align*} 
	
	When $K/F$ is unramified, $f_{K/F} = 2$, so $a(\pi \times \pi^\vee)$ is even by equation~(\ref{eqn:adjoint_conductor}). Suppose that $K/F$ is ramified and has residue characteristic different than 2. Let $\varpi = \varpi_K$, $\varpi_F$ be uniformizers of $K$, $F$, respectively. Then $\varpi_K^\varrho = -\varpi_K$. Also, since $f_{K/F} = 1$, $\O_K/\varpi_K \iso \O_F/\varpi_F$. There is a filtration on the unit group $U_K$
	$$U_K^{0} = U_K, \quad U_K^{i} = 1 + \varpi_K^i \O_K\quad  \text{for $i \geq 1$}$$
	with quotients:
	\begin{equation}\label{eqn:quotients_of_units}
		U_K^{0}/U_K^{1} \iso (\O_K/\varpi_K)^\times, \quad U_K^i/U_K^{i+1} \iso \O_K/\varpi_K.
	\end{equation}
	
	We show that if $\xi(\xi^\varrho)^{-1}|_{U_K^i} = \mathbbm 1$ for $i$ odd, then $\xi(\xi^\varrho)^{-1}|_{U_K^{i-1}} = \mathbbm 1$. 
	
	For $i = 1$, if $\xi(\xi^\varrho)^{-1}|_{U_K^1} = \mathbbm 1$, then $\xi(\xi^\varrho)^{-1}(x)$ for $x \in U_K$ depends only on the residue class of $x$ (equation~(\ref{eqn:quotients_of_units})). We may hence assume $x \in \O_F$ since $\O_K/\varpi_K \iso \O_F /\varpi_F$. Then
	$$\xi(\xi^\varrho)^{-1}(x) = \xi(x) \xi(x^\varrho)^{-1} = 1.$$
	
	Similarly, for $i > 1$ odd, if $(\xi(\xi^\varrho)^{-1})|_{U_K^i} = \mathbbm 1$, then $\xi(\xi^\varrho)^{-1} (1 + \varpi^{i-1} x)$ for $x \in \O_K$ depends only on the residue class of $x$ (equation~(\ref{eqn:quotients_of_units})). We may hence assume $x \in \O_F$ since $\O_K/\varpi_K \iso \O_F /\varpi_F$. Then
	$$\xi(\xi^\varrho)^{-1}(1 + \omega_K^{i-1} x) = \xi(1 + \omega_K^{i-1} x) \xi(1 + (-\omega_K)^{i-1} x^\varrho)^{-1} = 1.$$
	
	Therefore, $a(\xi(\xi^\varrho)^{-1})$ is even, which completes the proof.
\end{proof}

\begin{remark}
	The strategy in the proof of Proposition~\ref{prop:adjoint_conductor} gives an explicit formula for $a(\pi \times \pi^\vee)$ in terms of $a(\xi)$ when $p \neq 2$. For example, when $K/F$ is ramified:
	$$a(\pi \times \pi^\vee) = \begin{cases}
		a(\xi) + 2 & \text{if $a(\xi)$ is even} \\
		a(\xi) + 1 & \text{if $a(\xi)$ is odd}.
	\end{cases}$$
	A similar result was obtained by Nelson--Pitale--Saha~\cite[Proposition~2.5]{Nelson_Pitale_Saha} when $F = \Q$ and the central character of $\pi_v$ is trivial.
	
	It would be interesting to compare these formulas with the ones given in~\cite{Bushnell_Henniart-Kutzko}, but we have not attempted to do this.

\end{remark}

\begin{remark}\label{rmk:adjoint_conductor}
	In fact, Nelson--Pitale--Saha~\cite{Nelson_Pitale_Saha} prove that the adjoint conductor is always a square when $F = \Q$ and $f$ has trivial Nebentypus. We expect that the adjoint conductor is a square also in our more general setting. However, proving this would require a careful analysis of dyadic representations~\cite[Chapter 12]{Bushnell_Henniart} and we decided not to pursue it here.
\end{remark}

We are finally ready to prove Theorem~\ref{thm:Stark_implies_conj}.

\begin{proof}[Proof of Theorem~\ref{thm:Stark_implies_conj}]
	By construction of $\varrho_f$~\cite{Rogawski_Tunnell},
	$$L(1, f, \Ad) = L(1, \Ad^0\varrho_f).$$
	Then, by Corollary~\ref{prop:adjoint_and_Pet}, we have that
	$$\langle f, f \rangle = \langle f^\varrho, f^\varrho \rangle \sim_{E^\times} \pi^{-2d} L(1, \Ad(f), \overline{\iota}) =  \pi^{-2d} L(1, \Ad^0 \varrho_f, \overline{\iota}).$$
	By Proposition~\ref{prop:Stark_for_HMF}, Stark's conjecture for $\Ad^0 \varrho_f$ is equivalent to the statement:
	$$L(1, \Ad^0 \varrho_f, \overline{\iota}) \sim_{E^\times} \frac{\pi^{2d}}{f_\varrho^{1/2}} \cdot \det R_f.$$
	Putting these together and noting that $W(\varrho) = \pm 1$ and $f_\varrho$ is a square away from $p=2$ (Proposition~\ref{prop:adjoint_conductor}) gives the result.
\end{proof}

\subsection{Further evidence}\label{subsec:further_ev}

We now present further evidence for the conjecture which may be deduced from Stark's conjecture~\ref{conj:Stark}.

We first observe that we have an algebraic operation given by complex conjugation. Recall that the vector space $H^{|J|}(X, \mathcal E_{\underline k(J), r})$ is defined over the field $F(\underline k(J)) \subseteq F$ which is totally real, and hence
$H^{|J|}(X, \mathcal E_{\underline k(J), r}) \otimes_{F(J)} \C \iso H^{|J|}(X_\C, \mathcal E_{k(J), r})$ has an action of complex conjugation~$F_\infty$. By definition, it preserves the rational structure $H^{|J|}(X, \mathcal E_{\underline k(J), r})$.

\begin{lemma}\label{lemma:cc}
	The complex conjugation $F_\infty \colon H^{|J|}(X_\C, \mathcal E_{k(J), r}) \to H^{|J|}(X_\C, \mathcal E_{k(J), r})$ is given on the basis $\omega_f^I$ where $|I| = |J|$ and $k(I) = k(J)$ by
	$$\omega_f^I \mapsto \omega_{f^\varrho}^I,$$
	where $f^\varrho(z) = \overline{f(-\overline z)}$ is Shimura's complex conjugation. In particular, on $f$-isotypic subspaces, it defines a map:
	$$F_\infty \colon H^{|J|}(X, \mathcal E_{k(J), r})_{f} \to H^{|J|}(X, \mathcal E_{k(J), r})_{f^\varrho}.$$
\end{lemma}
\begin{proof}
	This is a paraphrase of an observation of Harris~\cite[pp.\ 164]{Harris_periods_I}.
\end{proof}

\begin{proposition}
	There is an $E$-linear isomorphism $U_f \iso U_{f^\varrho}$. In particular, Conjecture~\ref{conj:Harris_periods} for $f$ is equivalent to Conjecture~\ref{conj:Harris_periods} for $f^\varrho$.
\end{proposition}
\begin{proof}
	The first assertion follows from the observation that $\varrho_f^\vee \iso \overline{\varrho_f} = \varrho_{f^\varrho}$, so we can realize $\Ad^0 \varrho_f \subseteq \varrho_f \otimes \varrho_{f^\varrho}$. Since $\varrho_f \otimes \varrho_{f^\varrho} \iso \varrho_{f^\varrho} \otimes \varrho_f$, we have that $\Ad^0 \varrho_f \iso \Ad^0 \varrho_{f^\varrho}$. This induces an isomorphism $U_f \iso U_{f^\varrho}$.
\end{proof}

Next, recall that we have a Serre duality pairing~(\ref{eqn:SD}):
\begin{equation}\label{eqn:SD_again}
	\langle -,- \rangle_{\mathrm{SD}} \colon H^{|J|}(X, \mathcal E_{\underline k(J), r})_f \otimes H^{|\Sigma_\infty \setminus J|}(X, \mathcal E_{\underline k(\Sigma_\infty \setminus J), r})_{f^\varrho} \to E(J)
\end{equation}
which is $E(J)$-rational. We modify it slightly to replace $f^\varrho$ with $f$ via Lemma~\ref{lemma:cc}.

\begin{definition}\label{def:pairing}
	We define a pairing
	$$\langle - , - \rangle \colon H^j(X, \E_{\underline 1, 1})_f \times H^{d- j}(X, \E_{\underline 1, 1})_f \to E^\times$$
	by $\langle -, - \rangle = \langle -, F_\infty(-) \rangle_{\mathrm{SD}}$.
\end{definition}

\begin{proposition}\label{prop:qiffd-q}
	Assume Stark's conjecture~\ref{conj:Stark}. Conjecture~\ref{conj:Harris_periods} in cohomological degree~$j$ if equivalent to Conjecture~\ref{conj:Harris_periods} in cohomological degree $d-j$ (up to a factor of $\sqrt{2}$).
\end{proposition}
\begin{proof}
	Recall that Conjecture~\ref{conj:Harris_periods} in cohomological degree $j$ states that the elements:
	$$\left(\textstyle \bigwedge^j A \right) \begin{pmatrix}
		\omega_f^{J_1} \\
		\vdots\\
		\omega_f^{J_{\binom{d}{j}}}
	\end{pmatrix} $$
	give a rational basis of $H^j(X, \mathcal E_{\underline 1, 1})_f$. Let us assume that this is true and prove that the elements
	$$\left(\textstyle \bigwedge^{d-j} A \right) \begin{pmatrix}
		\omega_f^{J_1} \\
		\vdots\\
		\omega_f^{J_{\binom{d}{d-j}}}
	\end{pmatrix} $$
	are rational in $H^{d-j}(X, \mathcal E_{\underline 1, 1})_f$. It is enough to check that each of these classes pair rationally with the classes in $H^j(X, \mathcal E_{\underline 1, 1})_f$ using the pairing $\langle -, - \rangle$. Note that the pairing $\langle -, - \rangle$ is induced by cup product and
	$$\langle \omega_f^J, \omega_f^{J'} \rangle = \begin{cases}
		\pm \langle f, f \rangle & \text{if }J' = \Sigma_\infty \setminus J, \\
		0 & \text{otherwise}.
	\end{cases}$$
	Since $A = R_f^{-1}$ and $\langle f, f \rangle \sim_{E^\times} f_{\varrho, 2}^{1/2} \det R_f$ by Theorem~\ref{thm:Stark_implies_conj}, this completes the proof.
\end{proof}

%

Now, suppose that $j = d-j$, i.e.\ $d = 2j$ is even and we consider the middle degree sheaf cohomology. Definition~\ref{def:pairing} then gives a non-degenerate bilinear pairing
$$\langle - , - \rangle \colon H^j(X, \mathcal E_{\underline 1, 1})_f \otimes  H^j(X, \mathcal E_{\underline 1, 1})_f \to E$$
which satisfies:
$$\langle \omega_1, \omega_2 \rangle = (-1)^j \langle \omega_2, \omega_2 \rangle.$$

\begin{proposition}\label{prop:matrix_of_pairing}
	Suppose $d = 2j$ is even. Consider the basis of $H^j(X_\C, \mathcal E_{\underline 1, 1})_f$ given by the entries of the vector
	$$\left( \textstyle \bigwedge^j A \right) \begin{pmatrix}
		\omega_f^{J_1} \\
		\vdots\\
		\omega_f^{J_{\binom{d}{j}}}
	\end{pmatrix},$$
	ordered so that the pairs $\omega_f^J$ and $\omega_f^{\Sigma_\infty \setminus J}$ are consecutive. Then the of the pairing $\langle -, - \rangle \otimes \C$ is block-diagonal with $2 \times 2$ blocks given by
	$$\begin{pmatrix}
		0 & \ast \\
		(-1)^j\ast & 0
	\end{pmatrix}.$$
	Moreover, assuming Stark's conjecture~\ref{conj:Stark}, we have that $\ast \in {E}[\sqrt{2}]^\times$.
\end{proposition}
\begin{proof}
	This follows from the same arugment as the proof of Proposition~\ref{prop:qiffd-q}.
\end{proof}

\begin{corollary}\label{cor:d=2det}
	When $d = 2$, we showed in Example~\ref{ex:d=2} that Conjecture~\ref{conj:Harris_periods} predicts that
	\begin{align*}
		u_1^\vee \star f & = \frac{ \log|\tau(u_{22}) | \cdot \omega_f^1 - \log|\tau (u_{21}) | \cdot \omega_f^2}{\det R_f} \in H^1(X, \mathcal E_{\underline 1,1})_f \otimes E, \\
		u_2^\vee \star f & = \frac{-\log|\tau (u_{12}) | \cdot \omega_f^1 +\log|\tau (u_{11}) |\cdot \omega_f^2}{\det R_f} \in H^1(X, \mathcal E_{\underline 1,1})_f \otimes E.
	\end{align*}
	Assuming Stark's conjecture~\ref{conj:Stark}, the determinant of this basis lies in $E[\sqrt{2}]^\times$.
\end{corollary}
\begin{proof}
	Suppose that $\omega_1, \omega_2$ is a rational basis of $H^1(X, \ \mathcal E_{\underline 1,1})_f \otimes E$ and
	\begin{align*}
		\omega_1 & = a u_1^\vee \star f  + bu_2^\vee \star f \\
		\omega_2 & = c u_1^\vee \star f  + d u_2^\vee \star f
	\end{align*}
	for some $a, b, c, d \in \C$. Then:
	\begin{align*}
		\langle \omega_1, \omega_2 \rangle & = \langle  a u_1^\vee \star f  + bu_2^\vee \star f , c u_1^\vee \star f  + d u_2^\vee \star f \rangle \\
		& = (ad - bc) \langle u_1^\vee \star f, u_2^\vee \star f \rangle,
	\end{align*}
	showing that $ad - bc \in E[\sqrt{2}]^\times$ by Proposition~\ref{prop:matrix_of_pairing}, assuming Stark's conjecture~\ref{conj:Stark}. Finally, this shows that $(u_1^\vee \star f) \wedge (u_2^\vee \star f) = \frac{\omega_1 \wedge \omega_2 }{ad - bc}$ is $E[\sqrt{2}]^\times$-rational.
\end{proof}

\section{Evidence: base change forms}\label{sec:BC}

Let $F_0$ be a totally real number field and consider a totally real extension $F$ of $F_0$. Any Galois representation of $G_{\overline{\Q}/F_0}$ may be restricted to a Galois representation  $G_{\overline{\Q}/F}$. Hence, according to Langlands' functoriality conjecture, for any automorphic representation~$\pi_0$ of $\Res_{F_0/\Q} \GL_{2, F_0}$, there exists an associated {\em base change} representation $\pi$ of $\Res_{F/\Q} \GL_{2, F}$, written $\pi = \mathrm{BC}_{F_0}^F \pi_0$. This is discussed in detail and proved when $F/F_0$ is a cyclic Galois extension in \cite{Langlands:BC}. See also \cite{Arthur_Clozel}.

We now make the following definition.

\begin{definition}
	A Hilbert modular form $f$ for $F$ is a {\em base change form} from $F_0$, if the associated automorphic representation $\pi$ is equal to $\mathrm{BC}_{F_0}^F \pi_0$ for some automorphic representation~$\pi_0$.
\end{definition}

Of course, this leaves the following question: given a Hilbert modular form $f_0\in \pi_0$, how to choose an explicit Hilbert modular form $f \in \pi = \mathrm{BC}_{F_0}^F \pi_0$? As far as we know, there is no canonical choice of $f$ in this generality.

When $F$ is a real quadratic extension of $F_0 = \Q$ and the weight of $f_0$ is at least two, one can define~$f$ as a theta lift of $f_0$, called the {\em Doi--Naganuma lift}. The reader can consult~\cite{DoiNaganuma, Naganuma, Zagier} for the original results and \cite[Ch.~III]{Oda} or \cite[Ch.~VI.4]{van_der_Geer} for an overview. In examples below, we will primarily be interested in cases where the level of $f_0$ is coprime to the discriminant of $F$; such cases were treated by Kumar--Manickam~\cite{Kumar_Manickam}. When $f_0$ has weight one, we are not aware of an explicit construction of the base change of $f_0$ to a real quadratic extension in the literature. We expect these forms can be constructed using the theta correspondence as above. 

We will instead satisfy ourselves with the fact that these forms exist according to the Strong Artin Conjecture, which is known in several relevant cases~\cite{Kassaei, Kassaei_Sasaki_Tian}.

\begin{definition}\label{def:BC_wt1}
	Let $f_0$ be a normalized parallel weight one Hilbert modular eigenform for $F_0$ and $\varrho_0$ be the associated Artin representation. The {\em base change of $f_0$} to $F$ is the normalized parallel weight one Hilbert modular eigenform~$f$ whose associated Galois representation is $\varrho_f = \Res_{G_{\overline \Q/ F}} \varrho_0$.
\end{definition}

The goal of this section is to consider Conjecture~\ref{conj:Harris_periods} for base change forms. We compute Stark units for base change forms, give a more explicit from of the conjecture in this case, and provide numerical evidence for it in the case of real quadratic extensions.

\subsection{Stark units for base change forms}

For a Hilbert modular form $f$ which is the base change of $f_0$, we want to relate the unit groups $U_f$ and $U_{f_0}$. We fix a common splitting field~$L$ which is Galois over $\Q$. We denote the three Galois groups by:
$$G = G_{L/\Q} \supseteq G_0' = G_{L/F_0} \supseteq G' = G_{L/F}.$$
If $\varrho_0 \colon G_0' \to \GL_2(E)$ is the Artin representation associated to $f_0$, then the Artin representation $\varrho$ associated to $f$ is $\varrho = \Res^{G_0'}_{G'} \varrho_0$ by our definition of base change forms (Definition~\ref{def:BC_wt1}).

The goal of this section is to discuss the relation between the Stark unit groups and regulators for~$f$ and~$f_0$.

\begin{proposition}\label{prop:units_for_BC}
	\leavevmode
	\begin{enumerate}
		\item We have a natural isomorphism:
		$$U_{f}\iso U_L[\Ad^0 \varrho_0 \otimes P],$$
		where $P$ is the permutation representation of $G_0'$ on the cosets $G_0'/G'$. 
		
		\item In particular, if we consider the $G_0'$-invariant subrepresentation
		$$P_0 = \mathrm{span} \left\{ \sum_{\sigma G' \in G_0'/G'} \sigma G'  \right\} \subseteq P,$$
		then
		$$U_{f_0} \iso U_L[ \Ad^0 \varrho_0 \otimes  P_0] \subseteq U_f.$$
	\end{enumerate}
\end{proposition}
\begin{proof}
	Part (2) clearly follows from part (1), so we just prove part (1). We have that:
	\begin{align*}
		U_{f} & = \Hom_{G'}(\Ad^0 \varrho, \Res^G_{G'} U_L) \\
		& = \Hom_{G'}(\Res_{G'}^{G_0'} \Ad^0 \varrho_0,  \Res^G_{G'} U_L) \\
		& = \Hom_{G}( \Ind^G_{G'} \Res_{G'}^{G_0'} \Ad^0 \varrho_0, U_L) \\
		& = \Hom_{G}( \Ind_{G_0'}^G( \Ad^0 \varrho_{0} \otimes P ), U_L), \\
		& = \Hom_{G_0'}(\Ad^0 \varrho_0 \otimes P, U_L)
	\end{align*}
	as claimed. The penultimate equality follows form the following fact from representation theory: if $K \subseteq H \subseteq G$ and $V$ is a representation of $H$, then
	\begin{align*}
		\Ind_K^G \Res^{H}_K V \iso \bigoplus_{g \in G/K} g (\Res^H_K V) \iso \bigoplus_{g \in G/H} \bigoplus_{h \in H/K} gh (\Res^H_K V)  & \iso \bigoplus_{g \in G/H} g(V \otimes P) \iso \Ind_{H}^G (V \otimes P) \\
		g (h v) & \mapsto g  ( h \cdot v \otimes hK) 
	\end{align*}
	where $P$ is a permutation representation of $H$ on the cosets $H/K$.
\end{proof}

Suppose now that $F_0 = \Q$ for simplicity.

\begin{proposition}\label{prop:units_for_BC_explicit}
	\leavevmode
	\begin{enumerate}
		\item Let $f$ be the base change of a modular form $f_0$ of weight one. Then the units $u_{jk}^f$ associated to $f$ as in Definition~\ref{def:Stark_reg_HMF} are given by:
		$$u_{jk}^f = \prod_{\sigma' \in G'} (\epsilon^{(\sigma_k \sigma' \sigma_j^{-1})^{-1}})^{a^0(\sigma_k \sigma' \sigma_j^{-1})_{11}},$$
		where $a^0(\sigma)$ is the matrix of $\Ad^0 \varrho_0(\sigma)$ in the basis $m_{i, 0}$.
		\item For any $j$, we have that
		$$\prod_{k=1}^d u_{jk}^f = u_{f_0}.$$
		In particular,
		$$R_f \begin{pmatrix}
			1 \\
			\vdots \\
			1
		\end{pmatrix} =\log|u_{f_0}| \begin{pmatrix}
			1 \\
			\vdots \\
			1
		\end{pmatrix}.$$
	\end{enumerate}
\end{proposition}
\begin{proof}
	For part (1), we may take $M_j  = \Ad^0 \varrho_0(\sigma_j)$ for $j=1,\ldots, d$ in Corollary~\ref{cor:Stark_basis_HMF} to get this expression for $u_{jk}^f$. Part (2) then follows from Proposition~\ref{prop:units_for_BC}~(2).
\end{proof}

%
%

\begin{corollary}\label{cor:units_for_real_quadratic_BC}
	Suppose  $[F : \Q] = 2$. Let $u_{f_0}$ be the unit associated to $f_0$ and $u_{f_0}^F$ be the unit associated to the Artin representation $\Ad^0 \varrho_{0} \otimes \omega_{F/\Q}$, where $\omega_{F/\Q}$ is the quadratic character associated to the extension $F/\Q$. Then:
	\begin{align}
		u_{11} \cdot u_{12} & = u_{f_0},\label{eqn:u11u12=uf0} \\ 
		u_{21} \cdot u_{22} & = u_{f_0}, \\
		u_{11} \cdot u_{12}^{-1} & = u_{f_0}^F, \\
		u_{21}^{-1} \cdot u_{22} & = u_{f_0}^F. \label{eqn:u21u22-1=uf0F}
	\end{align}
	In particular,
	\begin{equation}\label{eqn:R_f_BC}
		R_f = \begin{pmatrix}
			\log|u_{11}| & \log|u_{12}| \\
			\log|u_{21}| & \log|u_{22}|
		\end{pmatrix} = \begin{pmatrix}
			1 & -1 \\
			1 & 1
		\end{pmatrix}  \begin{pmatrix}
			\log|u_{f_0}| & 0 \\
			0 & \log|u_{f_0}^F|
		\end{pmatrix}  \begin{pmatrix}
			1 & -1 \\
			1 & 1
		\end{pmatrix} ^{-1}.
	\end{equation}
\end{corollary}
\begin{proof}
	Fix representatives $\sigma_1, \sigma_2$ of $G/G'$ and assume that $\sigma_1 \in G'$. Then the permutation representation $P$ of $G$ on $G'/G$ decomposes as
	$$P \iso \Q(\sigma_1 + \sigma_2) \oplus \Q(\sigma_1 - \sigma_2).$$
	Therefore,
	$$U_f \iso U_{f_0} \oplus U_{L}[\Ad^0 \varrho_0 \otimes \omega_{F/\Q}]$$
	by Proposition~\ref{prop:units_for_BC}. Tracing through this isomorphism under the chosen bases, we obtain equations~\eqref{eqn:u11u12=uf0}--\eqref{eqn:u21u22-1=uf0F} and the resulting equation~\eqref{eqn:R_f_BC}.
\end{proof}

\subsection{Consequences of Conjecture~\ref{conj:Harris_periods}}

Recall that we can use the matrix $R_f^{-1}$ to predict which cohomology classes in $H^1(X_\C, \mathcal E_{\underline 1,1})$ are rational. When $f$ is the base change of a modular form $f_0$, Proposition~\ref{prop:units_for_BC_explicit}~(2) implies that:
$$R_f^{-1} \begin{pmatrix}
	1 \\
	\vdots \\
	1
\end{pmatrix} = \frac{1}{ \log|u_{f_0}|} \begin{pmatrix}
	1 \\
	\vdots \\
	1
\end{pmatrix}.$$
Therefore, the following is a consequence of Conjecture~\ref{conj:Harris_periods}.

\begin{conjecture}
	Suppose $f$ is the base change of a modular form $f_0$ of weight one. Then the cohomology class
	$$\frac{1}{\log|u_{f_0}|}  \sum_{j=1}^d \omega_f^{\sigma_j} \in H^1(X, \mathcal E_{\underline 1,1})_f$$
	is rational.
\end{conjecture}

When $[F:\Q] = 2$, Corollary~\ref{cor:units_for_real_quadratic_BC} gives the following stronger rationality statement.

\begin{conjecture}\label{conj:real_quad}
	Suppose $[F:\Q] = 2$ and $f$ is the base change to $F$ of a modular form $f_0$ of weight one. Then a rational basis for $H^1(X, \mathcal E_{\underline 1, 1})_f$ is given by:
	$$\frac{\omega_f^{\sigma_1} + \omega_f^{\sigma_2}}{\log|u_{f_0}|}, \quad \frac{\omega_f^{\sigma_1} - \omega_f^{\sigma_2}}{\log|u_{f_0}^F|}.$$
\end{conjecture}

In light of Corollary~\ref{cor:Stark_implies_conj}, this is equivalent to Conjecture~\ref{conj:motivic_action}.

\subsection{Embedded Hilbert modular varieties}

To check if Conjecture~\ref{conj:Harris_periods} is compatible with base change, we consider the Hilbert modular variety for $F_0$ embedded in the Hilbert modular variety for $F$.

We will write $d = [F : F_0]$ and $d' = [F_0 : \Q]$. Let $\tau_1, \ldots, \tau_{d'}$ be the infinite places of $F_0$. Above each place $\tau_i$, there are $d$ infinite places $\sigma_{i,j}$ for $j = 1, \ldots, d$ of $F$.  We write $\zeta_i$, $i = 1, \ldots, d'$, for the variables on $\H \otimes F_0$ and $z_{i,j}$, $i = 1, \ldots, d'$, $j = 1, \ldots, d$ for the variables on $\H \otimes F$. Here $\zeta_i$ corresponds to $\tau_i$ and $z_{i,j}$ corresponds to $\sigma_{i,j}$. 

We write $X_0$ and $X$ for the Hilbert modular varieties associated to $F_0$ and $F$, respectively. There is a natural embedding 
$$\iota \colon X_0 \hookrightarrow X.$$ 
Over $\C$, it descends from the map
\begin{align*}
	\H \otimes F_0 & \hookrightarrow \H \otimes F \\
	(\zeta_1, \ldots, \zeta_{d'}) & \mapsto (\zeta_1, \ldots, \zeta_1, \zeta_2, \ldots, \zeta_2, \zeta_{d'}, \ldots, \zeta_{d'}),
\end{align*}
i.e.\ the subvariety is given by the equation $z_{i,j} = \zeta_i$ for all $i, j$.

We are interested in the restriction map
$$H^i(X, \mathcal E_{\underline 1, 1}) \overset{\iota^\ast}\to H^i(X_0, \mathcal E_{\underline d, d}).$$
Particularly, we defined a class $\omega_f^J \in H^I(X, \mathcal E_{\underline 1, 1})$ associated to $f \in H^0(X, \mathcal E_{\underline 1, 1})$ which is represented by
\begin{equation}\label{eqn:f^I}
	\omega_f^J(\underline z) = f(\underline z^J) \cdot y^J \cdot \bigwedge_{\sigma_{i,j} \in J} \frac{d z_{i,j} \wedge d\overline{z_{i,j}}}{y_{i,j}^2}
\end{equation}
as a Dolbeault class, and we consider $\iota^\ast(\omega_f^J)$.

\begin{lemma}
	If $J$ contains $\sigma_{i,j}$ and $\sigma_{i,j'}$ for $j \neq j'$, 
	$$\iota^\ast(\omega_f^J) = 0.$$
\end{lemma}
\begin{proof}
	This follows immediately from the expression~(\ref{eqn:f^I}) and the identity $z_{i,j} = \zeta_i$ on $X_0$.
\end{proof}

Let us assume that $J$ only contains at most one $\sigma_{i,j}$ for each $i$, so that it is possible that $\iota^\ast(\omega_f^J)$ is non-zero.

The following conjecture is a consequence of Conjecture~\ref{conj:Harris_periods}.

\begin{conjecture}\label{conj:cons_of_conj}
	Let $A = (a_{ij}) = R_f^{-1}$ be the inverse of the Stark regulator matrix. Then for all $j = 1\ldots, d$:
	$$\sum_{i=1}^n a_{ij} \iota^\ast( \omega_f^{\sigma_i} ) \in H^{|I|}(X_0, \mathcal E_{\underline{d}, d}) \otimes {E} \subseteq H^{|I|}((X_0)_\C, \mathcal E_{\underline{d}, d}) \otimes {E}.$$
\end{conjecture}

Note that it is possible that $\iota^\ast(\omega_f^{\sigma_j}) = 0$ for all $j$ in which case this conjecture is void. In fact, we expect that $\iota^\ast(\omega_f^{\sigma_j}) = 0$ if $f$ is not a base change form from $F_0$ (see Proposition~\ref{prop:Oda} for an example of this phenomenon).

\subsection{The case of real quadratic extensions}\label{subsec:quad_ext}

We finally restrict our attention to real quadratic extensions $F/\Q$. In the previous notation, $F_0 = \Q$ and $d = 2$. We denote by $z_1, z_2$ (instead of $z_{1,1}, z_{1,2}$)  the variables on $X_\C$ and by $z$ (instead of $\zeta_1$) the variable on $(X_0)_\C$.

Let $f$ be a holomorphic Hilbert modular form of parallel weight $(k, k)$ and consider $\omega_f^{\sigma_1} \in H^1(X_\C^\an, \mathcal E_{(2-k,k)}^\an)$, given by:
\begin{equation}
	\omega_f^{\sigma_1}(z_1, z_2) = f(\epsilon_1 \overline{z_1}, \epsilon_2  z_2) y_1^{k}\, \frac{dz_1 \wedge d\overline{z_1}}{y_1^2}.
\end{equation}

There are embedded modular curves $\iota \colon C \hookrightarrow X$ in the Hilbert modular surface, studied extensively by Hirzebruch--Zagier~\cite{HirzebruchZagier}.  We only consider the simplest example which is obtained by considering the map:
\begin{align*}
	\iota \colon C_\C^\an & \hookrightarrow X_\C^\an \\
	z & \mapsto (z, z)
\end{align*}
over $\C$ which descends to varieties over $\Q$. Via this map,
$$\iota^\ast(\mathcal E_{(2-k, k)}^\an) \iso \mathcal E_{2}^\an \iso \Omega_{C}^{1, \an}(\infty)$$
by the Kodaira--Spencer isomorphism, where $(\infty)$ indicates that differentials are allowed to have poles of orders at most one at the cusps. Hence:
\begin{align*}
	\iota^\ast(\omega_f^{\sigma_1})(z) & = f(\epsilon_1 \overline z, \epsilon_2  z) y^k \frac{d z \wedge d\overline z}{y^2}
\end{align*}
defines a class in $H^1(C_\C^\an, \Omega_C^{1, \an}(\infty))$. Via the trace map, we have:
\begin{align*}
	\mathrm{Tr} \colon H^1(C_\C^\an, \Omega_C^{1, \an}(\infty)) & \overset \iso\to \C, \\
	\iota^\ast(\omega_f^{\sigma_1})(z) & \mapsto \int\limits_{C_\C^\an} f(\epsilon_1   \overline z, \epsilon_2  z) y^k \frac{d z \wedge d\overline z}{y^2},
\end{align*}
and the isomorphism respects rational structures.

\begin{lemma}\label{lemma:restrict_classes}
	For a Hilbert modular form of weight $(k, k)$, $\iota^\ast(\omega_f^{\sigma_1}) = (-1)^{k+1}\iota^\ast( \omega_f^{\sigma_2} )$.
\end{lemma}
\begin{proof}
	It suffices to check that $\mathrm{Tr}(\iota^\ast(\omega_f^{\sigma_1})) = \mathrm{Tr}( \iota^\ast( \omega_f^{\sigma_2} ))$. This follows by a change of variables:
	\begin{align*}
		\Tr (\iota^\ast(\omega_f^{\sigma_1})) &= \int\limits_{C_\C^\an} f(\epsilon_1 \overline z, \epsilon_2 z) y^k \, \frac{dz \wedge d\overline z}{y^2} \\
		& = - \int\limits_{C_\C^\an} f(- \epsilon_1 z, -\epsilon_2 \overline z) y^k \, \frac{d z \wedge d \overline z}{y^2}\\
		& = (-1)^{k+1}  \int\limits_{C_\C^\an} f( (-\epsilon_1^{-1}) z, (-\epsilon_2^{-1}) \overline z) y^k \, \frac{dz \wedge d\overline z}{y^2} & \begin{pmatrix}
			1 \\
			& \epsilon^2
		\end{pmatrix} \in \Gamma, N(\epsilon) = -1 \\
		& = (-1)^{k+1} \Tr(\iota^\ast(\omega_f^{\sigma_2})),
	\end{align*}
	as claimed.
\end{proof}


Putting this together with Conjectures~\ref{conj:real_quad} and~\ref{conj:cons_of_conj}, we get the following conjecture. 

\begin{conjecture}\label{conj:cons_of_conj_quadratic}
	Let $f$ be the base change of a weight one modular form~$f_0$. Then:
	$$\int\limits_{C_\C^\an} f(\epsilon_1 \overline z, \epsilon_2 z) y^k \frac{dz \wedge d\overline z}{y^2} \sim_{E^\times} \log |u_{f_0}| .$$
\end{conjecture}


For $k \geq 2$ and full level, these integrals were considered by Asai~\cite{Asai}. The following result was also obtained~by Oda~\cite{Oda}. See also~\cite[Proposition (VI.7.9)]{van_der_Geer}.

\begin{proposition}[{\cite[Theorem 16.5]{Oda}}]\label{prop:Oda}
	Suppose $f$ is a Hilbert modular form of parallel weight $k \geq 2$ and level one.  If $f$ is not a base change form, then
	$$\int\limits_{C_\C^\an} f(\epsilon_1   \overline z, \epsilon_2 z) y^k \frac{d z \wedge d\overline z}{y^2} = 0.$$
	Otherwise, if $f$ is the Doi--Naganuma lift of a modular form $g$ of weight $k \geq 2$, level $D = \mathrm{disc}(F/\Q)$, and character~$\omega_{F/\Q}$, then there is a constant $c \in \Q^\times$ such that 
	\begin{equation}\label{eqn:Oda}
		\int\limits_{C_\C^\an} f(\epsilon_1  \overline z, \epsilon_2  z) y^k \frac{d z \wedge d\overline z}{y^2} = c \frac{\langle f, f \rangle}{\langle g, g \rangle}.
	\end{equation}
\end{proposition}

\begin{remark}\label{rmk:BC_proof_d=2}
	The proof of Proposition~\ref{prop:Oda} in loc.\ cit.\ uses the explicit realization of $f$ as a Doi--Naganuma lift of a modular form $g$, which is currently not available in the literature for weight one forms. If an appropriate analogue of Proposition~\ref{prop:Oda} holds for a weight one forms $f_0$ of arbitrary weight, level, and character, then we expect that Stark's conjecture~\ref{conj:Stark} implies Conjecture~\ref{conj:cons_of_conj_quadratic} for base change forms of~$f_0$ to a real quadratic fields.
	
	%
	
	Verifying the details of this would take us too far afield, so we will pursue this elsewhere. Instead, in the next section we describe some explicit numerical computations that support Conjecture~\ref{conj:cons_of_conj_quadratic}.
\end{remark}

We end this section by proving that Conjectures~\ref{conj:real_quad} and~\ref{conj:cons_of_conj_quadratic} are equivalent for base change forms, as long as $\iota^\ast(\omega_f^{\sigma_1}) \neq 0$.

\begin{proposition}\label{prop:d=2_equivalent}
	Let $f$ be the base change of a weight one modular form $f_0$. Assume:
	\begin{enumerate}
		\item Stark's conjecture for the adjoint representation associated to $f$,
		\item $\iota^\ast(\omega_f^{\sigma_1}) \neq 0$,
	\end{enumerate}
	Then Conjecture~\ref{conj:real_quad} for $f$ is equivalent to Conjecture~\ref{conj:cons_of_conj_quadratic} for $f$, up to a potential factor of $\sqrt{2}$.
\end{proposition}
\begin{proof}
	Clearly, Conjecture~\ref{conj:real_quad} implies Conjecture~\ref{conj:cons_of_conj_quadratic}. We will prove the converse.
	
	Consider the algebraic map $\varphi \colon X \to X$ given on $X_\C \to X_\C$ by $(z_1, z_2) \mapsto (z_2, z_1)$. By examining the proof of Theorem~\ref{thm:Harris_Su}, one can deduce that if $f$ is a base change form, then $\varphi$ preserves $f$-isotypic components of coherent cohomology and hence induces a map:
	$$\varphi^\ast \colon H^1(X, \mathcal E_{\underline 1, 1})_f \to H^1(X, \mathcal E_{\underline 1,1})_f.$$
	Clearly, $\varphi^\ast_\C(\omega_f^{\sigma_1}) = \omega_f^{\sigma_2}$ and $\varphi^\ast_\C(\omega_f^{\sigma_2}) = \omega_f^{\sigma_1}$. Letting
	$$\omega^{\pm} = \omega_f^{\sigma_1} \pm \omega_f^{\sigma_2},$$
	we see that $\varphi^\ast_\C(\omega_f^{\pm}) = \pm \omega_f^{\pm}$. Hence $\omega_f^\pm$ are eigenvectors for the linear map $\varphi_\C^\ast$ with distinct eigenvalues, and so there exist $\lambda^\pm \in \C$ such that:
	$$\lambda^\pm \omega^\pm \in H^1(X, \mathcal E_{\underline 1,1 })_f.$$
	
	We have a rational functional 
	$\Tr \circ \iota \colon H^1(X, \mathcal E_{\underline 1, 1})_f \otimes E \to E$ such that:
	$$(\Tr \circ \iota)(\lambda^+ \omega^+) \sim_{E^\times} \lambda^+ \int\limits_{C_\C^\an}  \iota^\ast( \omega_f^{\sigma_1}) ,\quad (\Tr \circ \iota)(\lambda^- \omega^-) = 0 $$
	by Lemma~\ref{lemma:restrict_classes}. Conjecture~\ref{conj:cons_of_conj_quadratic} then shows that we may take
	$$\lambda^+ = \frac{1}{\log|u_{f_0}|}.$$
	
	Finally, by Corollary~\ref{cor:d=2det}, we know that the determinant of the basis
	$$\frac{\omega^+}{\log|u_{f_0}|}, \frac{\omega^-}{\log|u_{f_0}^F|}$$
	is $E[\sqrt{2}]$-rational, and hence
	$$\lambda^+ \cdot \lambda^- \sim_{E[\sqrt{2}]^\times} \frac{1}{\log|u_{f_0}| \cdot \log|u_{f_0}^{F}|},$$
	showing that we may take $\lambda^- = \frac{1}{\log|u_{f_0}^F|}$.
\end{proof}

\begin{remark}
	The idea to use the map $\varphi$ was communicated to us by the referee for a previous version for this manuscript. We thank them for this suggestion.
\end{remark}

\begin{remark}\label{rmk:assume_chi_0_quadratic}
	We expect that the condition (2) in Proposition~\ref{prop:d=2_equivalent} (i.e.\  $\iota^\ast(\omega_f^{\sigma_1}) \neq 0$) is equivalent to the character $\chi_0$ of $f_0$ being quadratic. One implication is clear: $\iota^\ast(\omega_f^{\sigma_1}) $ transforms by the character $\chi_0^2$ under the action of $\Gamma_0(N)$, and hence $\mathrm{Tr}(\iota^\ast(\omega_f^{\sigma_1})) = 0$ unless $\chi_0^2 = 1$. Conversely, if $\chi_0^2 = 1$, then the global analogue of Jacquet's conjecture~\cite{Kable:Jacquet, Prasad:Jacquet} implies that the automorphic representation $\pi$ generated by $f$ contains a non-zero $\GL_2(\A_\Q)$-invariant functional. We predict that $f \mapsto \iota^\ast(\omega_f^{\sigma_1})$ is this functional, i.e.\ $\iota^\ast(\omega_f^{\sigma_1}) \neq 0$.
	
	Finally, we expect that Proposition~\ref{prop:d=2_equivalent} has a refinement when $\chi_0^2 \neq 1$. If $\omega_0$ is the character of $\A_\Q^\times$ corresponding to $\chi_0$ by class field theory  and $\widetilde{\omega_0}$ is its extension to $\A_F^\times$ (which always exists), then Jacquet's conjecture predicts that the representation $\pi \otimes \widetilde{\omega_0}^{-1}$ has a non-zero $\GL_2(\A_\Q)$-invariant functional. One could hope to translate this to a classical statement analogous to Conjecture~\ref{conj:cons_of_conj_quadratic}.
\end{remark}

\subsection{Computing the integrals numerically}

The next goal is to provide numerical evidence of Conjecture~\ref{conj:cons_of_conj_quadratic}, i.e.\ check that
\begin{equation}\label{eqn:integral_over_C}
	\int\limits_{C_\C^\an} f(\epsilon_1   \overline z, \epsilon_2 z) y \frac{d z \wedge d\overline z}{y^2} \sim_{E^\times} \log  |u_{f_0}|.
\end{equation}
We will assume that $\chi_0^2 = 1$ (c.f.\ Remark~\ref{rmk:assume_chi_0_quadratic}), and hence the integral may be taken over $\Gamma_0(N) \backslash \H$ instead of $\Gamma_1(N) \backslash \H$. Indeed, equation~\eqref{eqn:integral_over_C} is equivalent to:
\begin{equation}
	\int\limits_{\Gamma_0(N) \backslash \H} f(\epsilon_1   \overline z, \epsilon_2 z) y \frac{d z \wedge d\overline z}{y^2} \sim_{E^\times} \log  |u_{f_0}|,
\end{equation}
because the two integrals differ by a factor of $\varphi(N)$.

We first derive a formula (Theorem~\ref{thm:explicit_formula_for_integral}) for the integral on the left hand side using Nelson's technique \cite{Nelson} for evaluating integrals on modular curves. 

Let $\Gamma \subseteq \SL_2(\Z)$ be a finite index subgroup and let $F \colon \Gamma \backslash \mathcal H \to \C$ be a $\Gamma$-invariant function on the upper half plane $\H$. Suppose we have its $q$-expansions, i.e.\ for all $\tau \in \SL_2(\Z)$, we have
\begin{equation}\label{eqn:q_exp_general}
	F(\tau z) = \sum_{n \in \Q} a_F(n, y; \tau) e(nx)
\end{equation}
where $e(nx) = e^{2 \pi i n x}$.

\begin{theorem}[{\cite[Theorem 5.6]{Nelson}}]\label{thm:Nelson}
	Suppose $F$ is bounded, measurable, and satisfies $F(\tau z) \ll y^{- \alpha}$ for some fixed $\alpha > 0$, almost all $z = x + i y$ with $y \geq 1$, and all $\tau \in \SL_2(\Z)$. Then for $0 < \delta < \alpha$ we have that:
	$$\int\limits_{\Gamma \backslash \H} F(z) \frac{dx dy}{y^2} = \int\limits_{(1+\delta)} (2s-1) 2 \xi(2s) \sum_{\tau \in \Gamma \backslash \SL_2(\Z)} a_F(0, \cdot; \tau)^\wedge(1-s)\, \frac{ds}{2\pi i}$$
	where
	\begin{align*}
		\xi(2s) & = \frac{\Gamma(s)}{\pi^s} \zeta(2s), \\
		a_F(0, \cdot; \tau)^\wedge(1-s) & = \int\limits_{0}^\infty a_F(0, y; \tau) y^{s-1} \, \frac{dy}{y}.
	\end{align*}
\end{theorem}

Applying this to $F(z) = f_0(z) \cdot \overline{f_0(z)} \cdot y^k$ gives an explicit expression for the Petersson inner product $\langle f_0, f_0 \rangle$.

\begin{corollary}[{Nelson, \cite[Theorem 4.2]{Collins}}]\label{cor:compute_Pet}
	Suppose $f_0$ is a cusp form in $S_k(N, \chi)$. For a cusp $s$, let $\sum\limits_n a_{n, s} q^n$ be the $q$-expansion at $\infty$ of $f_0|[\tau_{s, h}]_k$, where $\tau_{s, h} = \tau_s \begin{pmatrix}
		h_s & 0 \\
		0 & 1
	\end{pmatrix}$ and $\tau_s \infty = s$. Then we have that:
	$$\langle f_0, f_0 \rangle = \frac{4}{\mathrm{vol}(\Gamma \backslash \H)} \sum_{s \in \Gamma \backslash \mathbb P^1(\Q)} \frac{h_{s, 0}}{h_{s}} \sum_{m=1}^\infty \frac{|a_{m, s}|^2}{m^{k-1}} \sum_{n=1}^\infty \left(\frac{x}{8 \pi} \right)^{k-1} (x K_{k-2}(x) - K_{k-1}(x)), \quad x = 4 \pi n \sqrt{\frac{m}{h_{s}}},$$
	where $K_v$ is a $K$-Bessel function, $h_{s, 0}$ is the classical width of the cusp $s$, and $h_s$ is the width described in \cite[Lemma 2.1]{Collins}.
\end{corollary}

\begin{remark}
	An algorithm to compute these Petersson inner products was developed and implemented by Collins~\cite[Algorithm 4.3]{Collins}.
\end{remark}

The goal for this section is to prove the following theorem, which is an explicit form of Theorem~\ref{thm:Nelson} in our case. 

Recall that for $\alpha \in \SL_2(\O_F)$, we write $\alpha_i = \sigma_i(\alpha)$ and
$$f|[\alpha]_{\underline k}(z_1, z_2) = f(\alpha_1 z_1, \alpha_2 z_2) j(\alpha_1, z_1)^{-k_1} j(\alpha_2, z_2)^{-k_2}$$
where
$$j(g, z) = \det(g)^{-1/2} (cz+d).$$

By definition, if $f$ is a Hilbert modular form of weight $(k_1, k_2)$ and level $\Gamma$ and character $\chi$, then $f|[\alpha]_{\underline k} = \chi(d) \cdot f$ for $\alpha = \begin{pmatrix}
	a & b \\
	c & d
\end{pmatrix} \in \Gamma$.

\begin{theorem}\label{thm:explicit_formula_for_integral}
	Let $f$ be a normalized parallel weight $k$ Hilbert modular newform of level~$\mathfrak N$ and character $\chi$. For each cusp $s \in \mathbb P^1(\Q)/\Gamma_0(N)$, let $\tau \in \SL_2(\Z)$ satisfy $\tau \infty = s$. Let $h_s$ be the width of the cusp as described in \cite[Lemma 2.1]{Collins}, and
	\begin{align*}
		\tau^\epsilon & = \begin{pmatrix}
			\epsilon & 0 \\
			0 & 1
		\end{pmatrix} \tau \begin{pmatrix}
			\epsilon & 0 \\
			0 & 1
		\end{pmatrix}^{-1}, \\
		\tau^{\epsilon}_h & = \tau^{\epsilon} \begin{pmatrix}
			h_s & 0 \\
			0 & 1
		\end{pmatrix}.
	\end{align*}
	If $\sum\limits_{m \gg 0} a_{(m), s} q^m$ is the $q$-expansion of $f|[\tau^{\epsilon}_h]_k$ at $\infty$, then
	$$\int\limits_{\Gamma_0(N) \backslash \H} f(\epsilon_1   \overline z, \epsilon_2 z) y^k\, \frac{d z \wedge d\overline z}{y^2} = 4 \sum_{s} \frac{h_{s,0}}{h_s} \sum_{m=1}^\infty\frac{a_{(m), s}}{(m/ \sqrt{d})^{k-1}} \sum_{n=1}^\infty \left(\frac{x}{2^{3-i}\pi}\right)^{k-1} \left(xK_{k-2}(x) - K_{k-1}(x) \right)$$
	where $x = 2^{2-i/2}\pi n \sqrt{\frac{m}{h_s \sqrt{d}}}$ and $h_{s, 0}$ is the classical width of the cusp~$s$, and $i = 0$ if $d \equiv 1 \ (4)$ or $i = 1$ if $d \equiv 3\ (4)$.
\end{theorem}

\begin{remark}
	This formula is very similar to the formula for $\langle f_0, f_0 \rangle$ in Corollary~\ref{cor:compute_Pet}. We can hence adapt the algorithm~\cite[Algorithm 4.3]{Collins} to compute the integral. The computation of $q$-expansions of $f$ at other cusps given the $q$-expansion at $\infty$ is discussed in the next section (\ref{sec:q-exp_other_cusps}).
\end{remark}

We devote the rest of this section to the proof of this theorem. We want to apply Theorem~\ref{thm:Nelson} to the function
\begin{equation}
	F(z) = F_f^{\sigma_1}(z) = f(\epsilon_1  \overline z, \epsilon_2  z) \cdot y^k
\end{equation}
where $f$ is a Hilbert modular form of parallel weight $k$.

We will need $q$-expansions of $F(z)$ at other cusps, i.e.\ $q$-expansions of $F(\tau z)$ for $\tau \in \SL_2(\Z)$, as in equation~(\ref{eqn:q_exp_general}). The idea is to express them in terms of $q$-expansions at $\infty$ of another Hilbert modular form.

\begin{lemma}\label{lemma:other_cusps}
	Suppose $f$ is a Hilbert modular form of weight $(k,k)$. For a cusp $s$, let $\tau \in \SL_2(\O_F)$ be such that $\tau \infty = s$ and set
	$$\tau^\epsilon = \begin{pmatrix}
		\epsilon & 0 \\
		0 & 1
	\end{pmatrix} \tau \begin{pmatrix}
		\epsilon & 0 \\
		0 & 1
	\end{pmatrix}^{-1}.$$
	Then we have that:
	$$F_f^{\sigma_1}(\tau z) = F_{f|[\tau^{\epsilon}]_k}^{\sigma_1}(z).$$
\end{lemma}
\begin{proof}
	For $\tau \in \SL_2(\Z)$, we have that:
	\begin{align*}
		f(\epsilon_1 (\tau z_1), \epsilon_2 (\tau z_2)) & = f \bigg|\left[\begin{pmatrix}
			\epsilon & 0 \\
			0 & 1
		\end{pmatrix} \tau \right]_k (z_1, z_2) \cdot (N_{F/\Q}(\epsilon))^{-k/2} \cdot j(\tau, z_1)^{k} j(\tau, z_2)^{k} \\
		& = f \bigg|\left[ \tau^{\epsilon} \begin{pmatrix}
			\epsilon & 0 \\
			0 & 1
		\end{pmatrix} \right]_k (z_1, z_2) \cdot (N_{F/\Q}(\epsilon))^{-k/2}  \cdot j(\tau, z_1)^{k} j(\tau, z_2)^{k} \\ 
		& = f|[\tau^{\epsilon}]_k (\epsilon_1  z_1, \epsilon_2  z_2)  \cdot j(\tau, z_1)^k j(\tau, z_2)^k.
	\end{align*}
	
	Therefore,
	\begin{align*}
		F_f^{\sigma_1}(\tau z) & = f(\epsilon_1 (\tau \overline {z}), \epsilon_2  (\tau z)) \cdot (\mathrm{Im} (\tau z))^k \\
		& = f|[\tau^{\epsilon }]_k (\epsilon_1  \overline z, \epsilon_2  z)  \cdot |j(\tau, z)|^{2k} \cdot (\mathrm{Im}(\tau z))^k \\
		& = f|[\tau^{\epsilon }]_k (\epsilon_1  \overline z, \epsilon_2  z)  \cdot \mathrm{Im}(z)^k \\
		& = F_{f|[\tau^{\epsilon }]_k } (z),
	\end{align*}
	since $\mathrm{Im}(\tau z) =  |j(\tau, z)|^{-2} y$.
\end{proof}

\begin{lemma}\label{lemma:smooth_q-exp_general}
	For a cusp $s$, consider $\tau \in \SL_2(\Z)$ such that $\tau \infty= s$. Let $h_s$ be the width of cusp $s$ (as in \cite[Lemma 2.1]{Collins}) and
	$$\tau_h^{\epsilon} = \tau^{\epsilon} \begin{pmatrix}
		h_s & 0 \\
		0 & 1
	\end{pmatrix}.$$
	
	The $q$-expansion coefficients of $F(\tau z)$ (as in equation~(\ref{eqn:q_exp_general})) are given by
	$$a_F(n/h_s, y; \tau) = (y/h_s)^k \cdot \sum_{\substack{m \gg 0 \\ \mathrm{Tr}(\epsilon m) = n}} a_{(m), s} \cdot e^{- 2\pi (\epsilon_2 m_2/\delta_2 - \epsilon_1 m_1/\delta_1) y/h_s},$$
	where $a_{(m), s}$ are Fourier coefficients of $f|[\tau_h^{\epsilon}]_k$. In particular,
	$$a_F(0, y; \tau) = (y/h_s)^k \cdot \sum_{m = 1}^\infty a_{(m), s} \cdot e^{-2 \pi \frac{2^{1-i} m}{\sqrt{d}} (y / h_s)}$$
	where $i = 0$ if $d \equiv 1 \ (4)$ and $i = 1$ if $d \equiv 3\ (4)$.
\end{lemma}
\begin{proof}
	We write $h = h_s$ for simplicity. Suppose the $q$-expansion of $f|[\tau_h^{\epsilon}]_k$ is:
	$$f|[\tau_h^{\epsilon}]_k(z_1, z_2) = \sum_{m \gg 0} a_{(m), s} q^{m/\delta}.$$
	Then:
	$$f|[\tau^{\epsilon}]_k(z_1, z_2) = h^{-k} \sum_{m \gg 0} a_{(m), s} q^{m/ (\delta h)}.$$
	By Lemma~\ref{lemma:other_cusps},
	\begin{align*}
		F(\tau z) & = f|[\tau^{\epsilon}]_k (\epsilon_1 \overline z, \epsilon_2 z) \cdot y^k \\
		& = (y/h)^{k} \sum_{m \in \O_F^+} a_{(m), s} \cdot e^{2\pi i \left(\epsilon_1 m_1/\delta_1  (\overline{z}/h) + \epsilon_2 m_2/\delta_2 (z/h)  \right)} \\
		& = (y/h)^{k} \sum_{m \in \O_F^+} a_{(m), s} \cdot  e^{- 2\pi (\epsilon_2 m_2/\delta_2 - \epsilon_1 m_1/\delta_1) (y/h)} e^{2\pi i (\mathrm{Tr} \epsilon m/\delta) (x/h)} \\
		& =  (y/h)^{k} \sum_{n \in \Z} \left( \sum_{\substack{m \in \O_F^+ \\ \mathrm{Tr}(\epsilon m/\delta) = n }} a_{(m), s} \cdot e^{- 2\pi (\epsilon_2 m_2/\delta_2 - \epsilon_1 m_1/\delta_1) (y/h)} \right) e((n/h) x).
	\end{align*}
	
	Hence
	$$a_F(n/h, y; \tau) = (y/h)^k \cdot \sum_{\substack{m \gg 0 \\ \mathrm{Tr}(\epsilon m/\delta) = n}} a_{(m), s} \cdot e^{- 2\pi (\epsilon_2 m_2/\delta_2 - \epsilon_1 m_1/\delta_1) (y/h)},$$
	and in particular,
	$$a_F(0, y; \tau) = (y/h)^k \cdot  \sum_{\substack{m \gg 0 \\ \mathrm{Tr}(\epsilon m/\delta) = 0}} a_{(m), s} e^{- 2\pi (\epsilon_2 m_2/\delta_2 - \epsilon_1 m_1/\delta_1) (y/h)}.$$
	
	To make this last formula more explicit, we write $m = \alpha + \beta \sqrt{d}$. We may choose $\delta = 2^i \sqrt{d} \cdot \epsilon$ to be the totally positive generator of the different ideal. Then
	$$\epsilon m/\delta = \frac{\beta}{2^i} + \frac{\alpha}{2^i d} \sqrt{d}.$$
	If $\mathrm{Tr}(\epsilon m/ \delta) = 0$, then $\beta = 0$, so $m = \alpha \in \Z_{> 0}$. Moreover:
	$$\epsilon_2 m_2/\delta_2 - \epsilon_1 m_1/\delta_1 = \frac{2^{1-i} m}{\sqrt{d}}.$$
	
	We may hence rewrite the above sum as
	$$a_F(0, y; \tau) = (y/h)^k \cdot \sum_{m = 1}^\infty a_{(m), s} \cdot e^{-2 \pi \frac{2^{1-i} m}{\sqrt{d}} (y / h)},$$
	proving the lemma.
\end{proof}

We finally complete the proof of Theorem~\ref{thm:explicit_formula_for_integral}.

\begin{proof}[Proof of Theorem~\ref{thm:explicit_formula_for_integral}]
	We will apply Theorem~\ref{thm:Nelson} to the invariant function $F(z) = F_f^{\sigma_1}(z)$. By Lemma~\ref{lemma:smooth_q-exp_general},
	$$a_F(0, y; \tau) = (y/h_s)^k \cdot \sum_{m = 1}^\infty a_{(m), s} \cdot e^{-2 \pi \frac{2^{1-i} m}{\sqrt{d}} (y / h_s)}.$$
	
	Hence:
	\begin{align*}
		a_F(0, \cdot; \tau)^\wedge(1-t) & = \int\limits_{0}^\infty a_F(0, y; \tau) \, y^{t-1} \, \frac{dy}{y}. \\
		& = \sum_{m = 1}^\infty a_{(m), s}  \int\limits_0^\infty  e^{-2 \pi \frac{2^{1-i} m}{\sqrt{d}} (y / h_s)} \, y^{t-1} \, (y/h_s)^{k}\, \frac{dy}{y} \\
		& = \sum_{m = 1}^\infty a_{(m), s} h_s^{-k} \int\limits_0^\infty e^{-2 \pi \frac{2^{1-i} m}{h_s \sqrt{d}} y} \, y^{t + k - 1} \,\frac{dy}{y} \\
		& = \sum_{m = 1}^\infty a_{(m), s} h_s^{-k} \frac{\Gamma(t + k -1)}{( 2 \pi \frac{2^{1-i} m}{h_s \sqrt{d}})^{t + k - 1}}  \\
		& = \sum_{m = 1}^\infty \frac{a_{(m), s}}{(2^{2-i} \pi m / \sqrt{d})^{k-1} h_s} \frac{\Gamma(t+k-1)}{(2^{2-i} \pi \frac{m}{h_s \sqrt{d}})^{t}}.
	\end{align*}
	
	According to \cite[Lemma A.4]{Nelson}:
	$$\int\limits_{(1+\delta)} (t-1/2) \frac{\Gamma(t)\Gamma(t+\nu)}{(x/2)^{2t + \nu}} \frac{dt}{2\pi i} = x K_{\nu - 1}(x) - K_\nu(x)$$
	for $\nu \in \C$ with $\mathrm{Re}(\nu) \geq 0$.
	
	By Theorem~\ref{thm:Nelson},
	\begin{align*}
		\int\limits_{\Gamma \backslash \H} F(z) \frac{dx dy}{y^2} & = \int\limits_{(1+\delta)} (2t-1) 2 \xi(2t) \sum_{\tau} a_F(0, \cdot; \tau)^\wedge(1-t)\, \frac{dt}{2\pi i} \\
		& = 4  \sum_{s} h_{s,0} \sum_{m=1}^\infty 
		\frac{a_{(m), s}}{(2^{2-i} \pi m / \sqrt{d})^{k-1} h_s}
		\sum_{n=1}^\infty  
		\int\limits_{(1+\delta)}  
		(t-1/2) \frac{\Gamma(t) \Gamma(t + k -1)}{(2^{2-i} \pi^2 \frac{m}{h_s \sqrt{d}})^t} \frac{1}{n^{2t}} 
		\frac{dt}{2\pi i} \\
		& = 4 \sum_{s}  \frac{h_{s,0}}{h_s} \sum_{m=1}^\infty 
		\frac{a_{(m), s}}{(2^{2-i} \pi m / \sqrt{d})^{k-1}}
		\sum_{n=1}^\infty  
		\int\limits_{(1+\delta)} 
		(t-1/2) \frac{ \Gamma(t) \Gamma(t+k-1)}{(2^{2-i} \pi^2 \frac{m n^2}{h_s \sqrt{d}})^t} 
		\frac{ds}{2\pi i} \\
		&  = 4  \sum_{s}  \frac{h_{s,0}}{h_s} \sum_{m=1}^\infty\frac{a_{(m), s}}{(m/ \sqrt{d})^{k-1} h_s} \sum_{n=1}^\infty \left(\frac{x}{2^{3-i}\pi}\right)^{k-1} (xK_{k-2}(x) - K_{k-1}(x))
	\end{align*}
	where we set $x = 2^{2-i/2}\pi n \sqrt{m/h_s \sqrt{d}}$ in the last line.
\end{proof}

In order to use Theorem~\ref{thm:explicit_formula_for_integral}, we need to compute the $q$-expansions of the Hilbert modular form~$f$ at other cusps, i.e. $q$-expansions of $f|[\alpha]_k$ at $\infty$ for a matrix $\alpha$. We discuss this problem in the next section.

\subsection{$q$-expansions at other cusps}\label{sec:q-exp_other_cusps}

In this section, we address the following question: given the $q$-expansion of a Hilbert modular form $f(z)$ at the cusp $\infty$, what is the $q$-expansion of $f(z)$ at any cusp of $\Gamma_0(N) \backslash \H^2$?

We take two methods available for modular forms and discuss their generalization to Hilbert modular forms:
\begin{itemize}
	\item Asai's explicit formula~\cite{Asai:Fourier} (Theorem~\ref{thm:explicit_qexp_other_cusps}),
	\item Collins computational method based on a least-squares algorithm~\cite{Collins} (Algorithm~\ref{alg:q-exp}).
\end{itemize}

The first one is much faster in practice but only works for square-free level. The second one works for any level, but our implementation is too slow in practice to compute the above integrals. We include it here since it might be of independent interest.

Collins also introduced an improved computational method for modular forms using twists of eigenforms~\cite[Algorithm 2.6]{Collins}. This is also discussed in Chen's thesis~\cite[Chapter~4]{Chen}. 

An alternative approach is to use the adelic language. The Fourier coefficients of a modular form are given by value of the Whittaker newform of $f$ at certain matrices. Loeffler--Weinstein~\cite{Loeffler_Weinstein} give an algorithm to compute the local representations, so one just needs an algorithm to compute the local newforms. For more details, see~\cite[Section 3]{Corbett_Saha} .

\subsubsection{Explicit formula, following~\cite{Asai:Fourier}}

Let $F$ be a totally real field of narrow class number~1 (of arbitrary degree $d$). Suppose $f$ is a Hilbert modular eigenform of level $\mathfrak N$ with character $\chi \colon (\O_F/\mathfrak N)^\times \to \C^\times$ and parallel weight $k$. Suppose the level $\mathfrak N$ is square-free. We write $\Gamma = \Gamma_0(\mathfrak N)$ throughout this section.

The goal is to prove an explicit formula (Theorem~\ref{thm:explicit_qexp_other_cusps}) for the $q$-expansion of a Hilbert modular form~$f$ at a cusp $C = a/b \in F$ in terms of the $q$-expansion at $\infty$, generalizing the main result of \cite{Asai:Fourier} to the Hilbert modular case.

Since $\mathfrak N$ is square-free, the cusps $C = a/b$ of $\Gamma \backslash \H^2$ are in bijection with decompositions $\mathfrak N = \mathfrak A \cdot \mathfrak B$, where $\mathfrak B = ((b), \mathfrak N)$. For each divisor $\mathfrak A$, we consider the matrix
$$W_{\mathfrak A} = \begin{pmatrix}
	A\alpha & \beta \\
	N\gamma & A \delta
\end{pmatrix} = \begin{pmatrix}
	\alpha & \beta \\
	B \gamma & A \delta
\end{pmatrix}  \begin{pmatrix}
	A & 0 \\
	0 & 1
\end{pmatrix}$$
such that:
\begin{itemize}
	\item $A$, $N$ are totally positive generators of $\mathfrak A$, $\mathfrak N$, respectively; then $B = N/A$ is a totally positive generator of $\mathfrak B$,
	\item $\det W_{\mathfrak A} = A$,
	\item $\alpha, \beta, \gamma, \delta \in \O_F$.
\end{itemize}

Such a matrix always exists: since $\mathfrak A = (A)$ and $\mathfrak B = (B)$ are coprime, we have that $1 = \lambda A + \mu B,$ for some $\lambda, \mu \in \O_F$, so $A = \lambda A^2 + \mu N,$
and we may take $\alpha = \beta = 1$ and $\gamma = -\mu$, $\delta = \lambda$  to obtain such a matrix:
$$W_{\mathfrak A} = \begin{pmatrix}
	A & 1 \\
	-N \mu & A \lambda
\end{pmatrix}.$$
Conversely, for a matrix $W_{\mathfrak A}$,
$$W_{\mathfrak A}^{-1} \infty = \frac{\delta}{- B \gamma}$$
is a cusp with $((B \gamma), \mathfrak N) = \mathfrak B$, because
$$1 = A \alpha \delta - B \beta \gamma \equiv - B \beta \gamma \mod \mathfrak A,$$
so $(\gamma)$ is coprime to $\mathfrak A$.

Such a matrix $W_{\mathfrak A}$ associated to $\mathfrak A$ is well-defined up multiplication by elements of $\Gamma$. Moreover, $W_{\mathfrak A}$ normalizes $\Gamma$ and $A^{-1} W_{\mathfrak A}^2 \in \Gamma$.

The $q$-expansion of $f$ at the cusp corresponding to $\mathfrak N = \mathfrak A \mathfrak B$ is the $q$-expansion of the Hilbert modular form $f_{\mathfrak A} = f|W_{\mathfrak A}$ at $\infty$.

For a prime ideal $\p = (\varpi)$ of $\O_F$, coprime to $\mathfrak N$, with totally positive generator $\varpi$, the action of the Hecke operator $T(\mathfrak p)$ on the space of cusp forms $S_k(\mathfrak N, \chi)$ is given by
\begin{equation}\label{eqn:Hecke}
	f|T(\p) = N_{F/\Q}(\p)^{k/2-1} \left(\chi(\varpi) f \bigg|_k \begin{pmatrix}
		\varpi & 0\\
		0 & 1
	\end{pmatrix}  + \sum_{\nu \in \O_F/\p} f \bigg|_k \begin{pmatrix}
		1 & \nu \\
		0 & \varpi
	\end{pmatrix} \right).
\end{equation}

For example, when $d = 2$, this simplifies to the more familiar expression:
$$f|T(\mathfrak p) = N_{F/\Q}(\p)^{k-1} \left(\chi(\varpi) f(\varpi_1 z_1, \varpi_2 z_2) + N_{F/\Q} \p^{-k} \sum_{\nu \in \O_F/\p} f\left(\frac{z_1 + \nu_1}{\varpi_1}, \frac{z_2 + \nu_2}{\varpi_2}  \right)  \right).$$

We will write $T(\p, \chi)$ for the action of the Hecke operator $T(\mathfrak p)$ on $S_k(\mathfrak N, \chi)$. 

\begin{remark}
	This normalization of Hecke operators is consistent with $T'(\mathfrak p)$ in \cite{Shimura_HMF}.
\end{remark}

For simplicity, whenever we write down a generator of an ideal, it is assumed to be totally positive. The main result of this section is the following.

\begin{theorem}\label{thm:explicit_qexp_other_cusps}
	Let $f$ be a newform in $S_k(\mathfrak N, \chi)$ and $f|T(\mathfrak p, \chi) = a_{\mathfrak p} f$. For each decomposition $\mathfrak N = \mathfrak A \mathfrak B$, let $f_{\mathfrak A} = f| {W_{\mathfrak A}}$. Then $f_{\mathfrak A}$ is a newform in $S_k(\mathfrak N, \!\,^{\mathfrak A}\chi)$ and
	$$f_{\mathfrak A} | T(\mathfrak p, \!\,^{\mathfrak A}\chi) = a_{\mathfrak p}^{(\mathfrak A)} f_{\mathfrak A}$$ 
	for every prime $\mathfrak p = (\varpi)$, where
	$$a_{\mathfrak p}^{(\mathfrak A)} = \begin{cases}
		\overline{\chi_{\mathfrak A}}(\varpi) a_{\mathfrak p} & \text{if }\mathfrak p \not|\  \mathfrak A, \\
		\chi_{\mathfrak B}(\varpi) \overline{a_{\mathfrak p}} & \text{if } \mathfrak p \not|\ \mathfrak B,
	\end{cases}$$
	and
	\begin{align*}
		\chi_{\mathfrak A} \colon (\O_F/\mathfrak A \O_F)^\times & \to \C^\times, \\
		m & \mapsto \chi((-B\beta \gamma) m + (A\alpha \delta)), \\
		\chi_{\mathfrak B} \colon (\O_F/\mathfrak B \O_F)^\times & \to \C^\times, \\
		m & \mapsto  \chi((A\alpha \delta)m +  (-B\beta \gamma) ), \\
		{ \!\,^{\mathfrak A}\chi } \colon (\O_F/\mathfrak N\O_F)^\times & \to \C^\times, \\
		m & \mapsto \chi((A\alpha \delta) m +  (- B \beta \gamma) m^{-1}).
	\end{align*}
\end{theorem}
\begin{proof}
	The proof is a straightforward generalization of \cite[Theorem~1]{Asai:Fourier}, so we just give a sketch.
	
	We first check that $f_{\mathfrak A}$ has character ${ \!\,^{\mathfrak A}\chi }$ described above. Write
	\begin{align*}
		d \colon \Gamma = \Gamma_0(\mathfrak N) & \to (\O_F/\mathfrak N)^\times, \\
		\begin{pmatrix}
			a & b \\
			c & d
		\end{pmatrix} & \mapsto d \mod \mathfrak N.
	\end{align*}
	Then we just need to check that
	$$d \left(W_{\mathfrak A} g W_{\mathfrak A}^{-1}  \right) = {^{\mathfrak A}\chi} (d(g)),$$
	where
	$${^{\mathfrak A}\chi}(m) \equiv (A\alpha \delta) m +  (-B\beta \gamma) m^{-1} \mod \mathfrak N.$$
	
	For $g = \begin{pmatrix}
		a & b \\
		c & d
	\end{pmatrix}$, we have that
	\begin{align*}
		W_{\mathfrak A} \gamma W_{\mathfrak A}^{-1} & =  \begin{pmatrix}
			A\alpha & \beta \\
			N\gamma & A \delta
		\end{pmatrix} \begin{pmatrix}
			a & b \\
			c & d
		\end{pmatrix}  \begin{pmatrix}
			\delta & -\beta/A \\
			-B\gamma & \alpha
		\end{pmatrix} \\
		& = \begin{pmatrix}
			A \alpha & \beta \\
			N \gamma & A \delta
		\end{pmatrix} \begin{pmatrix}
			a\delta - b B \gamma & -a \beta/A + b \alpha \\
			c \delta - d B \gamma & -c\beta/A+ d \alpha
		\end{pmatrix}
	\end{align*}
	so
	\begin{align*}
		d(W_{\mathfrak A} \gamma W_{\mathfrak A}^{-1}) & = -a \beta \gamma B + bN\alpha \gamma - c \beta \delta + dA\alpha \delta \\
		& \equiv  (- \beta \gamma B)a + (A\alpha \delta) d \mod \mathfrak N & \text{since }c \equiv 0 \mod \mathfrak N
	\end{align*}
	which proves the above result, since $a d \equiv  1 \mod \mathfrak N$.
	
	One then computes a formula for how the Hecke operator $T(\mathfrak p, \chi)$ commutes with $W_{\mathfrak A}$ using the above expression for Hecke operators (c.f.~\cite[Lemma 2]{Asai:Fourier}). To check that $f_{\mathfrak A}$ is a newform, one shows that $W_{\mathfrak A}$ preserves oldforms (c.f.~\cite[Lemma 1]{Asai:Fourier}).
\end{proof}

The Hecke eigenvalues $a_{\mathfrak n}$ of $T({\mathfrak n})$ may be computed from the eigenvalues $a_{\mathfrak p}$ of $T(\mathfrak p)$ in the standard way~\cite[(2.26)]{Shimura_HMF}. For $\mathfrak n$ coprime to $\mathfrak m$, we have that
$$a_{\mathfrak n \mathfrak m} = a_{\mathfrak n} \cdot a_{\mathfrak m}$$
and for $\mathfrak n = \mathfrak p^r$, we have that
\begin{equation}\label{eqn:other_evals}
	\sum_{r = 0}^\infty a_{\mathfrak p^r} N(\mathfrak p)^{-rs} = [1 - a_{\mathfrak p} N(\mathfrak p)^{-s} + \chi(\mathfrak p) N(\mathfrak p)^{k_0 - 1 - 2s}]^{-1}
\end{equation}
where $k_0 = \max\{k_1, \ldots, k_n\}$.

We can then recover the $q$-expansion of $f_{\mathfrak A}$, up to a constant $\lambda$, from the Hecke eigenvalues~$a_{\mathfrak p}^{(\mathfrak A)}$ given by Theorem~\ref{thm:explicit_qexp_other_cusps}. There is an explicit expression for $\lambda$, described in the next theorem.

\begin{theorem}\label{thm:lambda}
	Let $f$ be a normalized Hilbert newform with character $\chi$ and level $\mathfrak N$. Then there is a constant $\lambda$ such that
	$$f_{\mathfrak A} = \lambda \cdot \underbrace{\sum\limits_{\nu \gg 0} a_{(\nu)}^{(\mathfrak A)} q^\nu}_{f^{(\mathfrak A)}} $$
	where we define:
	\begin{align*}
		a_{(1)}^{(\mathfrak A)} & = 1 \\
		a_{(\nu)}^{(\mathfrak A)} & = \overline{\chi_{\mathfrak A}(\nu)} a_{(\nu)} & \text{if } ((\nu), \mathfrak A) = \O_F, \\
		a_{(\nu)}^{(\mathfrak A)} & = \chi_{\mathfrak B}(\nu) \overline{a_{(\nu)}} & \text{if } ((\nu), \mathfrak B) = \O_F, \\
		a_{(\nu \mu)}^{(\mathfrak A)} & = a_{(\nu)}^{(\mathfrak A)} a_{(\mu)}^{(\mathfrak A)}  & \text{if } (\nu, \mu) = \O_F.
	\end{align*}
	
	Moreover, there is an explicit formula for $\lambda$, analogous to~\cite[Theorem 2]{Asai:Fourier}. First, for a decomposition $\mathfrak N = \mathfrak p \mathfrak B$ for a prime ideal $\mathfrak p = (\varpi)$, let
	$$W_{\mathfrak p} = \begin{pmatrix}
		\varpi & 1 \\
		N \gamma & \varpi \delta
	\end{pmatrix}$$
	be a matrix of determinant $\varpi$ with $\gamma, \delta \in \O_F$. Then
	$$f|W_{\mathfrak p} = \lambda_{\mathfrak p} f^{(\mathfrak p)}$$
	with
	$$\lambda_{\mathfrak p} = \begin{cases}
		C(\chi_{\mathfrak p}) \cdot N\mathfrak p^{-k/2} \cdot \overline{a_{\mathfrak p}} & \text{if }\mathfrak p \text{ divides }\mathrm{cond}(\chi), \\
		-N\mathfrak p^{1-k/2} \cdot \overline{a_{\mathfrak p}} & \text{otherwise},
	\end{cases}$$
	where
	$$C(\chi_{\mathfrak p}) = \sum_{h \!\!\!\mod \mathfrak p} \chi_{\mathfrak p}(h) \cdot e^{2\pi i \mathrm{Tr}(h/\varpi)}$$
	is a Gauss sum associated to $\chi_{\mathfrak p}$.
	
	In general, for any $\mathfrak N = \mathfrak A \mathfrak B$ with an associated matrix $W_{\mathfrak A} = \begin{pmatrix}
		A \alpha & \beta \\
		N \gamma & A \delta
	\end{pmatrix}$, we have that
	$$\lambda = \chi(A \delta - B \gamma) \prod_{(\varpi) = \mathfrak p | \mathfrak A} \chi_{\mathfrak p} (A/\varpi) \lambda_{\mathfrak p}.$$
\end{theorem}
\begin{proof}
	Once again, the proof generalizes the proof of~\cite[Theorem 2]{Asai}. Since for $\mathfrak A$ coprime to $\mathfrak A'$, we may take $W_{\mathfrak A\mathfrak A'} = W_{\mathfrak A} W_{\mathfrak A'}$, it is enough to check the assertion for a prime ideal $\mathfrak A = \mathfrak p$.
	
	By definition of $a_{(\nu)}^{(\p)}$ and $\lambda_{\mathfrak p}$, we have that:
	\begin{equation}\label{eqn:lambda}
		f | T(\mathfrak p) \circ W_{\mathfrak p} = a_{\mathfrak p} f|W_{\p} = a_\p \lambda_\p \sum_{\nu \gg 0} a_{(\nu)}^{(\p)} q^{\nu/\delta}.
	\end{equation}
	We compute the left hand side in another way to get the result.
	
	Since $\det W_{\mathfrak p} = \varpi$, we have that
	$$B\gamma \equiv B \gamma - \varpi \delta = -1\mod \p.$$
	Hence for $j \not \equiv 1 \!\!\mod \p$, 
	$$1 + B \gamma j \equiv 1 - j \not\equiv 0 \mod \p,$$
	so there exists $\ell \not\equiv 0\!\! \mod \p$ such that
	$$(1 + B\gamma j) \ell \equiv 1 \mod \p.$$
	Moreover, this defines a bijection
	$$\{j \in \O_F/\p \ | \ j \not\equiv 1 \mod \p \} \leftrightarrow \{\ell \in \O_F/\p \ | \ j \not\equiv 0 \mod \p \}.$$
	One can then check that for $j \not \equiv 1 \!\! \mod \p$
	$$\begin{pmatrix}
		1 & j \\
		& \varpi
	\end{pmatrix} W_{\mathfrak p} = \sigma_1 \begin{pmatrix}
		1 & \ell \\
		& \varpi
	\end{pmatrix}  \begin{pmatrix}
		\varpi \\
		& 1
	\end{pmatrix}$$
	for some $\sigma_1 \in \Gamma_0(\mathfrak N)$ such that $\chi(d(\sigma_1)) = \chi_{\mathfrak p}(\ell)$.
	
	For $j = 1$, we have that:
	$$\begin{pmatrix}
		1 & 1 \\
		& \varpi
	\end{pmatrix} W_{\mathfrak p} = \sigma_2 W_{\mathfrak p} \begin{pmatrix}
		\varpi &  \\
		& 1
	\end{pmatrix}$$
	for some $\sigma_2 \in \Gamma_0(\mathfrak N)$ such that $\chi(d(\sigma_2)) = \chi_{\mathfrak B}(\varpi)$.
	
	Using the expression~(\ref{eqn:Hecke}) for $T(\mathfrak p)$:
	\begin{align*}
		f | T(\mathfrak p) \circ W_{\mathfrak p} & = (N_{F/\Q} \p)^{k/2 - 1} \left(\sum_{j \in \O_F/\p} f\bigg|_k \begin{pmatrix}
			1 & j \\
			& \varpi
		\end{pmatrix} W_{\mathfrak p} \right) \\
		& = (N_{F/\Q} \p)^{k/2 - 1} \left(\sum_{\ell \not\equiv 0}  \chi_{\mathfrak p}(\ell) f\bigg|_k \begin{pmatrix}
			1 & \ell \\
			& \varpi
		\end{pmatrix} \begin{pmatrix}
			\varpi & \\
			& 1
		\end{pmatrix} \right) + \chi_{\mathfrak B}(\varpi) f\bigg|_k W_{\mathfrak p} \begin{pmatrix}
			\varpi &  \\
			& 1
		\end{pmatrix}.
	\end{align*}
	Using the $q$-expansions:
	$$f = \sum_{\nu \gg 0} a_{(\nu)} q^{\nu/\delta}, \quad f|_kW_{\mathfrak p} = \lambda_{\mathfrak p} \sum_{\nu \gg 0} a_{(\nu)}^{(\mathfrak p)} q^{\nu/\delta},$$
	we have that
	\begin{align*}
		f\bigg|_k \begin{pmatrix}
			1 & \ell \\
			& \varpi
		\end{pmatrix} \begin{pmatrix}
			\varpi & \\
			& 1
		\end{pmatrix} & = \sum_{\nu \gg 0} a_{(\nu)} e^{2\pi i \Tr(\nu \ell/ \delta \varpi)} q^{\nu/\delta}, \\
		f\bigg|_k W_{\mathfrak p} \begin{pmatrix}
			\varpi &  \\
			& 1
		\end{pmatrix} & = (N_{F/\Q} \mathfrak p)^{k/2} \lambda_{\mathfrak p} \sum_{\nu \gg 0} a_{(\nu)}^{(\p)} q^{\nu \varpi/\delta}.
	\end{align*}
	Hence
	\begin{align*}
		f | T(\mathfrak p) \circ W_{\mathfrak p}  & = (N_{F/\Q} \p)^{k/2 - 1} \sum_{\nu \gg 0} a_{(\nu)} \left( \sum_{\ell \not\equiv 0}  \chi_{\mathfrak p}(\ell) e^{2\pi i \Tr(\nu \ell/ \delta\varpi)} \right) q^{\nu/\delta}  \\
		& \quad + (N_{F/\Q} \p)^{k-1} \chi_{\mathfrak B}(\varpi) \lambda_{\mathfrak p} \sum_{\nu \gg 0} a_{(\nu)}^{(\p)} q^{\nu \varpi/\delta}.
	\end{align*}
	
	If $\chi_{\mathfrak p}$ is primitive, then
	$$\sum_{\ell \not\equiv 0} \chi_{\mathfrak p}(\ell) e^{2 \pi i \Tr(\nu \ell/\delta \varpi)} = \overline {\chi_{\mathfrak p}(\nu)} \chi_{\mathfrak p}(\delta) C(\chi_{\mathfrak p})$$
	since $\delta$ is coprime to $\varpi$, and hence 
	$$f | T(\mathfrak p) \circ W_{\mathfrak p}  =  (N_{F/\Q} \p)^{k/2 - 1} \chi_\p(\delta) C(\chi_\p) \sum_{\nu \gg 0} \overline{\chi_\p(\nu)} a_{(\nu)} q^{\nu/\delta} + (N_{F/\Q} \p)^{k-1} \chi_{\mathfrak B}(\varpi) \lambda_{\mathfrak p} \sum_{\nu \gg 0} a_{(\nu)}^{(\p)} q^{\nu \varpi/\delta}.$$
	
	If $\chi_{\mathfrak p}$ is not primitive, then $\chi_{\mathfrak p} = \mathbbm{1}_{\p}$ is the trivial character modulo $\mathfrak p$. Then, since $\varpi$ is coprime to~$\delta$,
	$$\sum_{\ell \not\equiv 0} \chi_{\mathfrak p}(\ell) e^{2 \pi i \Tr(\nu \ell/\delta \varpi)} = \sum_{\ell \not\equiv 0} e^{2 \pi i \Tr(\nu \ell/\delta \varpi)} = \begin{cases}
		N(\mathfrak q) - 1 & \p | (\nu), \\
		-1 & \text{otherwise}.
	\end{cases}$$
	Hence:
	$$f | T(\mathfrak p) \circ W_{\mathfrak p}  =  -(N_{F/\Q} \p)^{k/2-1} \sum_{\nu \gg 0} a_{(\nu)} q^{\nu/\delta} + \sum_{\nu \gg 0} \left( (N\mathfrak p)^{k/2} a_{(\nu\varpi)} + (N_{F/\Q}\p)^{k-1} \chi_{\mathfrak B}(\varpi) \lambda_\p a_{(\nu)}^{(\p)}  \right) q^{\nu \varpi/\delta}.$$
	
	Comparing the expression for $f|T(\mathfrak p) \circ W_{\mathfrak p}$ in each case with equation~(\ref{eqn:lambda}) gives the result.
\end{proof}

\subsubsection{Numerical method, following~\cite{Collins}}

The explicit formulas above only apply to Hilbert modular forms of square-free level. We discuss how one could generalize a method of Collins to compute $q$-expansions at other cusps for general levels.

As in \cite[Section 2]{Collins}, we consider a matrix $\alpha$ which takes infinity to the cusp and
$$\alpha_h = \begin{pmatrix}
	a & b \\
	c & d
\end{pmatrix} \begin{pmatrix}
	h & 0 \\
	0 & 1
\end{pmatrix}.$$
For $f \in S_k(\Gamma_0(\mathfrak N))$,
$$f|[\alpha_h]_k \in S_k(\Gamma_0(\mathfrak N h))$$
and we want to compute its $q$-expansion:
\begin{equation}
	f|[\alpha_h]_k = \sum_{\nu \gg 0} a_{(\nu), \alpha} q^m = \sum_{\mathfrak n} a_{\mathfrak n, \alpha} \left( \sum_{m \in \Z} q^{u^m \nu} \right)
\end{equation}
where $q^m = e^{2\pi i \mathrm{Tr}(m/\delta)}$ and $u \in (\O_F)^{\times}_+$ is a fundamental unit.

The idea of Collins~\cite[Section 2.3]{Collins} is to sample points $z_1, \ldots, z_M \in \H^2$ and use the $q$-expansion at $\infty$ of $f$ to compute $f[\alpha_h]_k(z)$ for these values. Then to use a least squares algorithm to approximate the constants $a_{\mathfrak n, \alpha}$ which satisfy
$$f[\alpha_h]_k \approx \sum_{\mathfrak n} a_{\mathfrak n, \alpha} \left(\sum_{m \in \Z} q^{u^m \nu} \right).$$

\begin{algorithm}[$q$-expansion at other cusps, adapted from {\cite[Algorithm 2.3]{Collins}}]\label{alg:q-exp}
	Given:
	\begin{itemize}
		\item a Hilbert modular form $f$ of level $\mathfrak N$, weight $(k,k)$, with an algorithm to compute its Fourier coefficients $a_{\mathfrak n}$ for arbitrarily large $\mathfrak n$,
		\item a cusp $a/c \in \Q$ of width $h$,
		\item a maximal norm $K$ of Fourier coefficients needed,
		\item a desired accuracy $10^{-E}$,
		\item an exponential decay factor $e^{-C_0}$,
	\end{itemize}
	we can compute the Fourier coefficients $a_{\mathfrak n, \alpha}$ for $\mathrm{Norm}(\mathfrak n) < N$, accurate up to $10^{-E}$ as follows:
	
	\begin{enumerate}
		\item Either increase $K = K_0$ or decrease $C = C_0$ so that $KC \approx \log(10) E$ and work with interpolating
		$$\sum\limits_{\substack{\mathfrak n\\ N\mathfrak n \leq K}} a_{\mathfrak n, \alpha} \left(\sum_{m \in \Z} q^{u^m \nu} \right).$$
		
		
		%
		%
		\item Choose $M$ (for example, $2K_0$) and pick points $z_1, \ldots, z_M \in \H^2$ with both imaginary parts equal to $C/2\pi$ and $\mathrm{Re}(z_j)$ randomly in $(-d/ch - 1/2, -d/ch + 1/2)^2$.
		
		\item Numerically compute the values $f|[\alpha_h](z_j) = h^{k/2}(ch(z_{j,1})+d)^{-k}(ch(z_{j,2})+d)^{-k} f(\alpha_h z_j)$ using the $q$-expansion of $f$, truncating until we have reached an accuracy a little greater than $10^{-E}$, and fill these into a vector $b$.
		
		\item Numerically compute the values $\sum\limits_{m \in \Z} q^{u^m \nu}$ for each $z = z_1, \ldots, z_M$ with an accuracy a little greater than $10^{-E}$, and store them in a matrix $A$.
		
		\item Numerically find the least squares solution to $Ax = b$ as the exact solution to $(A^\ast A) x = A^\ast b$. The solution vector is our approximation to the coefficients $a_{\mathfrak n, \alpha}$ for each $\mathfrak n$ of norm at most~$K$.
	\end{enumerate}
\end{algorithm}

We implemented this algorithm, but step (3) is very slow in practice. Since we need a lot of Fourier coefficients in our case, it is not realistic to apply this algorithm for our purposes.

\subsection{Numerical evidence}

We can use Theorems~\ref{thm:explicit_formula_for_integral}, \ref{thm:explicit_qexp_other_cusps}, and~\ref{thm:lambda} to compute the integral and verify that:
\begin{equation}\label{eqn:const_c}
	\int\limits_{\Gamma_0(N) \backslash \H} f(\epsilon_1   \overline z, \epsilon_2 z) y \frac{d z \wedge d\overline z}{y^2} = c \cdot \log |u_{f_0}|
\end{equation}
for some $c \in E^\times$. This numerically verifies Conjecture~\ref{conj:cons_of_conj} which we showed is equivalent to Conjecture~\ref{conj:motivic_action} in base change cases.

%

\subsubsection{Modular forms associated to cubic extensions}

In Example~\ref{ex:base_change_units}, the unit group $U_{f_F}$ is described explicitly, so this is the first case we consider. This is the base change of Example~\ref{ex:cubic_fields} to a real quadratic extension~$F = \Q(\sqrt{d})$ of $\Q$.

We briefly recall Example~\ref{ex:cubic_fields} to set up the notation. Let $K = \Q(\alpha)$ be a cubic field of signature $[1,1]$, obtained by adjoining a root $\alpha$ of a cubic polynomial $P(x)$. The splitting field $L$ of $P(x)$ is the Galois closure of $K$ and $G_{L/\Q} \iso S_3$. We consider the irreducible odd Artin representation 
$$G_{L/\Q} \iso S_3 \to \GL_2(\Z).$$
It has an associated modular form $f_0$ and we consider its base change $f$ to~$F = \Q(\sqrt{d})$. The associated unit group is $U_{f_0} \iso U_K^{(1)}$, the norm 1 units of $K$, and we consider a generator $u = u_{f_0}$ of this group.

Table~\ref{table} shows constants $c \in \Q$ such that the equality~\eqref{eqn:const_c} holds up to at least 15 digits. The computations were performed on the High Performance Computing cluster Great Lakes at the University of Michigan.

\begin{table}[H]
	\centering
	\begin{tabular}{c|c|c|c|c|c|c}
		$d$ & polynomial $P(x)$ & \url{lmfdb.org} label & level $N$ & unit $u$ & constant $c$ & time taken \\ \hline 
		& & & & &  \\ [-0.3cm]
		5 & $x^3 - x^2 + 1$ & \verb|23.1.b.a| &23 & $\alpha^2 - \alpha$ & $2$ & 00:09:34  \\
		5 & $x^3 + x - 1$ & \verb|31.1.b.a| & 31 & $\alpha$ & $-4$ & 00:13:36 \\
		5 & $x^3 + 2x - 1$ & \verb|59.1.b.a| & 59 & $\alpha^2$ & $-8$ & 01:56:22 \\
		5 & $x^3 - x^2 + 2x + 1$ & \verb|87.1.d.b| & 87 & $\alpha$ & $-2$ & 04:15:09 \\[0.1cm] \hdashline
		& & & & &  \\ [-0.3cm]
		13 & $x^3 - x^2 + 1$ &\verb|23.1.b.a|  & 23 & $\alpha^2 - \alpha$ & $8$ & 00:10:19 \\
		13 & $x^3 + x - 1$ & \verb|31.1.b.a|   &31 & $\alpha$ & $-2$ & 00:49:47 \\
		13 &	$x^3 + 2x - 1$ & \verb|59.1.b.a| &59 & $\alpha^2$ & $-22$ & 29:47:44  \\
		13 & $x^3 - x^2 + 2x + 1$ & \verb|87.1.d.b| & 87 & $\alpha$ & $-4$ & 04:23:13 \\[0.1cm] \hdashline
		& & & & &  \\ [-0.3cm]
		17 & $x^3 - x^2 + 1$ & \verb|23.1.b.a|  & 23 & $\alpha^2 - \alpha$ & $14$ & 00:16:52 \\
		17 & $x^3 + x - 1$ & \verb|31.1.b.a|  &31 &  $\alpha$ & $-18$ & 01:01:15 \\
		17 & $x^3 - x^2 + 2x + 1$ & \verb|87.1.d.b| & 87 & $\alpha$ & $-14$ & 19:40:11 \\[0.1cm] \hdashline
		& & & & &  \\ [-0.3cm]
		29 & $x^3 - x^2 + 1$ & \verb|23.1.b.a| &23 & $\alpha^2 - \alpha$ & $4$ & 00:32:08 \\
		29 & $x^3 + x - 1$ & \verb|31.1.b.a|  & 31 &  $\alpha$ & $-14$ & 02:38:12 \\[0.1cm] \hdashline
		& & & & &  \\ [-0.3cm]
		37 & $x^3 - x^2 + 1$ &\verb|23.1.b.a|  & 23 & $\alpha^2 - \alpha$ & $10$ & 00:25:45 \\
		37 & $x^3 + x - 1$ & \verb|31.1.b.a|  &31 & $\alpha$ & $-6$ & 01:41:38  \\
	\end{tabular}
	\caption{
		This table presents constants $c$ such that equation~(\ref{eqn:const_c}) holds for the unit $u$ and the base change to $\Q(\sqrt{d})$ of the modular form of level $N$ associated to the polynomial $P(x)$. We give the \url{lmfdb.org} label of the modular form. The time taken to perform the computation with at least 15 digits of accuracy is displayed in the format hh:mm:ss.
	}
	\label{table}
\end{table}


It is quite remarkable that all the constants $c$ are even integers and not just rational numbers.  Rubin's integral refinement of Stark's conjecture~\cite{Rubin:Stark} could provide an explanation. Understanding this phenomenon may also be related to studying congruence numbers for~$f$~\cite{Doi_Hida_Ishii} and a potential integral refinement of Conjecture~\ref{conj:Harris_periods} would have to take them into account.

\subsubsection*{Weight one form of level $47$}

We give an example where the coefficients of $f_0$ are not rational and hence Stark's conjecture~\ref{conj:Stark} is not known for the base change form $f$. Let $f_0$ be the modular form of weight one, level $47$, label \verb|47.1.b.a| in \url{lmfdb.org}, and $q$-expansion:
$$f_0 = q + (-1+ \beta)q^2 - \beta q^3 + (1-\beta)q^4 + \cdots$$
where $\beta = \frac{1}{2}(1+ \sqrt 5)$.

The associated Galois representation is:
\begin{align*}
	\varrho \colon \Gal(L/\Q) \iso D_5 & = \langle s, r \ | \ s^2 =1,  r^5 = 1, srs = r^4\rangle \to \GL_2(\Z[\zeta_5]) \\
	s & \mapsto \begin{pmatrix}
		0 & 1 \\
		1 & 0
	\end{pmatrix}, \\
	r & \mapsto \begin{pmatrix}
		\zeta_5 & 0 \\
		0 & \zeta_5^{4}
	\end{pmatrix},
\end{align*}
where we choose $s \in D_5$ corresponding to the complex conjugation $c_0 \in \Gal(L/\Q)$ associated to $L \hookrightarrow \C$. For the basis $\begin{pmatrix}
	0 & 1 \\
	1 & 0
\end{pmatrix}$, $\begin{pmatrix}
	0 & 1 \\
	-1 & 0
\end{pmatrix}$, $\begin{pmatrix}
	1 & 0 \\
	0 & -1
\end{pmatrix}$ of $\Ad^0 \varrho$, the adjoint representation is:
\begin{align*}
	\varrho \colon \Gal(L/\Q) \iso D_5 & \to \GL_3(\Z[\zeta_5]) \\
	s & \mapsto \begin{pmatrix}
		1 \\
		& -1\\
		& & -1
	\end{pmatrix}, \\
	r & \mapsto \begin{pmatrix}
		(\zeta_5^{2} + \zeta_5^{-2})/2 & (\zeta_5^{2} - \zeta_5^{-2})/2 & 0  \\
		(\zeta_5^{2} - \zeta_5^{-2})/2 & (\zeta_5^{2} + \zeta_5^{-2})/2 &0 \\
		0& 0 & 1
	\end{pmatrix}.
\end{align*}
Finally, this shows that:
\begin{equation}\label{eqn:uN=47}
	u = u_{f_0} = \prod_{i=0}^4 (\epsilon^{r^{-i}})^{\zeta^{2i} + \zeta^{-2i}},
\end{equation}
where $\epsilon$ is the Minkowski unit (Definition~\ref{def:Minkowski}) for the embedding $\tau \colon L \hookrightarrow \C$ such that $s$ is the complex conjugation associated to $\tau$.

Note that $\beta = \zeta^2 + \zeta^{-2}$, so the coefficients $\zeta^{2i} + \zeta^{-2i}$ lie in the coefficient field $\Q(\sqrt{5})$ of $f$.

\begin{table}[H]
	\centering
	\begin{tabular}{c|c|c|c|c|c}
		$d$ & \url{lmfdb.org} label & level $N$ & unit $u_{f_0}$ & constant $c \in \Q(\sqrt{5})$ & time taken \\ \hline 
		& & & & &  \\ [-0.4cm]
		5 &  \verb|47.1.b.a| & 47 & \eqref{eqn:uN=47} & $ 1 - \frac{ \sqrt{5}}{5}$ & 04:44:15 \\ 
		& & & & &  \\ [-0.4cm]
		13 & \verb|47.1.b.a| & 47 &  \eqref{eqn:uN=47} & $5 - \sqrt{5} $ &  09:20:12 \\  
		& & & & &  \\ [-0.4cm]
		17 & \verb|47.1.b.a| & 47 & \eqref{eqn:uN=47} & $8 - 8\frac{\sqrt{5}}{5}$ & 02:04:28  \\ 
		& & & & &  \\ [-0.4cm]
		29 & \verb|47.1.b.a| & 47 &  \eqref{eqn:uN=47} & $3 - 3 \frac{\sqrt{5}}{5}$ &  15:47:31
	\end{tabular}
	\caption{
		This table presents constants $c$ such that equation~(\ref{eqn:const_c}) holds for the unit $u_{f_0}$ and the base change to $\Q(\sqrt{d})$ of the modular form $f_0$ of level $47$. The time taken to perform the computation with at least 15 digits of accuracy is displayed in the format hh:mm:ss.
	}
	\label{table2}
\end{table}

Interestingly, in this case, the right hand side seems to always be an integer multiple of $1 - \frac{\sqrt{5}}{5}$. Once again, this may be related to congruence numbers for $f$.

\appendix

\section{Comparison to Prasanna--Venkatesh~\cite{Prasanna_Venkatesh}}\label{app:Prasanna-Venkatesh}

Prasanna--Venkatesh gave a conjectural definition~\cite[Definition 4.2.1]{Prasanna_Venkatesh} of the adjoint motive. Beilinson's regulator defines a map
\begin{equation}\label{eqn:Beilinson_regulator}
	H^1_{\mathcal M}(M_{\rm coad}, \Q(1)) \to H_B(M_{\rm coad, \C}, \R)^{W_{\R}} \iso \hat{\g}^{W_\R}.
\end{equation}
For a cohomological, tempered automorphic representation, they define an action of $\bigwedge^\ast (\hat {\mathfrak g}^{W_\R})$ on Betti cohomology of the associated symmetric space, and conjecture that the action is rational for the rational structure given by motivic cohomology.

In this appendix, we explain that Conjecture~\ref{conj:motivic_action} is the natural analogue of this for coherent cohomology. In our case,
$$\hat \g \iso \bigoplus_{j=1}^d \mathfrak{sl}_{2, \C}.$$
The archimedean Langlands parameter associated to a Hilbert modular form $f$ of weight $(\underline k, r)$ is given by
\begin{align*}
	\varphi \colon W_\R = \C^\times \cup \C^\times j & \to \bigoplus_{j=1}^d \GL_2(\C) \\
	\C^\times \ni s e^{i \theta} & \mapsto  \begin{pmatrix}
		s^{2r} e^{i (k_j - 1) \theta} \\
		& s^{2r} e^{- i (k_j - 1) \theta} 
	\end{pmatrix} \\
	j  & \mapsto \begin{pmatrix}
		0 & (-1)^{k_j - 1}	 \\
		1 & 0
	\end{pmatrix}
\end{align*}
(see \cite{Knapp:LLC}). 

A simple computation of the adjoint action gives the following lemma.

\begin{lemma}\label{lemma:Deligne_coh}
	For a Hilbert modular form of weight $(\underline k, r)$, we have that:
	\begin{equation}\label{eqn:Deligne_coh_explicit}
		\hat \g^{W_{\mathbb R}} \iso \bigoplus_{j \text{ s.t.\ } k_j = 1} \R \begin{pmatrix}
			1 & 0 \\
			0 & -1
		\end{pmatrix}_j
	\end{equation}
\end{lemma}

This allows us to define the action of this Deligne cohomology group on coherent cohomology.

\begin{definition}\label{def:action_of_Deligne}
	Let $f$ be a Hilbert modular form of weight $(\underline k, r)$.  We define an action $\star$ of $\bigwedge^\ast \hat \g^{W_\R}$ on $H^\ast(X_\C, \mathcal E_{\underline k, r})_f$ by letting $ \begin{pmatrix}
		1 & 0 \\
		0 & -1
	\end{pmatrix}_j$for $j$ such that $k_j = 1$ act by:
	\begin{align*}
		H^j(X_\C, \mathcal E_{\underline k , r})_f & \to H^{j+1}(X_\C, \mathcal E_{\underline k , r})_f \\
		\omega_f^J & \mapsto \begin{cases} \omega_f^{J \cup \{\sigma_j\}} & \sigma_j \not\in J \\
			0 & \sigma_j \in J \end{cases} 
	\end{align*}
	Here, we use the bases of cohomology groups given in Corollary~\ref{cor:dim_of_coh_group}.
\end{definition}

This is precisely the action we defined in Definition~\ref{def:action}.

\begin{remark}
	Recall from Remark~\ref{rmk:cc_automorphically} that the cohomology class $\omega_f^J$ is associated to the action of right translation by the matrix $g_J \in G(\R)$ where
	$$(g_J)_j = \begin{cases}
		\begin{pmatrix}
			1 & 0 \\
			0 & -1
		\end{pmatrix} & \sigma_j \in J, \\
		\begin{pmatrix}
			1 & 0 \\
			0 & 1
		\end{pmatrix} & \sigma_j \not\in J.
	\end{cases}$$
	Although the elements in equation~\eqref{eqn:Deligne_coh_explicit} belong to the Lie algebra $\hat \g$ and not $G(\R)$, this seems like a natural way to define this action.
\end{remark}

In the case $(\underline k, r) = (\underline 1, 1)$, we expect from Proposition~\ref{prop:motivic_coh} that $U_f^\vee \iso H^1_{\mathcal M}(M_{\mathrm{coad}}, \Q(1))$. Proposition~\ref{prop:dual_units_over_C} gives an explicit expression for the (inverse of the) Beilinson regulator~\eqref{eqn:Beilinson_regulator}. Therefore, Conjecture~\ref{conj:motivic_action} amounts to the fact that the action of $H^1_{\mathcal M}(M_{\mathrm{coad}}, \Q(1))$ preserves the rational structure on coherent cohomology.

Finally, we briefly discuss the motivic action conjecture for partial weight one Hilbert modular forms. 
Suppose $f$ is a Hilbert modular form of weight $(\underline k, r)$ and let  $M = M_{\rm coad}$ be the conjectural coadjoint motive of weight zero associated to $f$. The Beilinson short exact sequence for $M$ is:
\begin{equation}\label{eqn:Beilinson_ses}
	0 \to F^1(H_{\rm dR}(M)) \otimes_\Q \R \to H_B(M_\R, \R) \to H^1_{\mathcal D}(M_\R, \R(1)) \to 0.
\end{equation}
A simple calculation using the Hodge decomposition of $H_B(M)$ gives:
$$\dim F^1(H_{\rm dR}(M)) = \#\{j \ | \ k_j > 1 \},$$
and hence
$$\dim H^1_{\mathcal D}(M_\R, \R(1)) = \#\{j \ | \ k_j = 1\}.$$
The last assertion is consistent with Lemma~\ref{lemma:Deligne_coh}.

%

Consider the rational structure on $H^1_{\mathcal D}(M_\R, \R(1))$ given by the motivic cohomology group $H^1_{\mathcal M}(M, \Q(1))$ via Beilinson's regulator~\eqref{eqn:Beilinson_regulator}. This gives an action $\star$ of $H^1_{\mathcal M}(M, \Q(1))$ on coherent cohomology $H^\ast(X_\C, \mathcal E_{\underline k, r})_f$ via Definition~\ref{def:action_of_Deligne}. 

\begin{conjecture}\label{conj:partial_weight_one}
	The action $\star$ of $\bigwedge^\ast H^1_{\mathcal M}(M, \Q(1))$ on $H^\ast(X_\C, \mathcal E_{\underline k, r})_f$ preserves the rational structure $H^\ast(X, \mathcal E_{\underline k, r})_f$.
\end{conjecture}

The action of top-degree elements, i.e.\ the group $\bigwedge^{\ell}  H^1_{\mathcal M}(M, \Q(1))$ where $\ell = \# \{j \ | \ k_j = 1\}$, has a particularly nice description in terms of  Beilinson's conjecture for the adjoint $L$-function. For $m \in \bigwedge^\ell H^1_{\mathcal M}(M, \Q(1))$, we have that
$$m\star f = \frac{\omega_f^{J_1}}{r_{\mathcal D}(m)} \in H^{\ell}(X_\C, \mathcal E_{\underline k, r})_f,$$
where $J_1 = \{j \ | \  k_j = 1\}$. This final space is one-dimensional according to Theorem~\ref{thm:Harris_Su}~(2) and hence we may check the rationality of $m \star f$ using Serre duality. We consider the rational element
$$\frac{\omega_f^{\Sigma_\infty \setminus J_1}}{\nu^{\Sigma_\infty \setminus J_1}(f)} \in H^{d - \ell}(X, \mathcal E_{\underline 2 - \underline k, r})_f$$
(see Definition~\ref{def:nu_Harris}). Then:
$$\left\langle \frac{\omega_f^{J_1}}{r_{\mathcal D}(m)}, \frac{\omega_f^{\Sigma_\infty \setminus J_1}}{\nu^{\Sigma_\infty \setminus J_1}(f)}  \right\rangle \sim_{E_f(\underline k)^\times} \frac{\langle f, f \rangle}{r_{\mathcal D}(m) \cdot\nu^{\Sigma_\infty \setminus J_1}(f)}$$
by Proposition~\ref{prop:fact_of_Pet}. Using Theorem~\ref{thm:adj_and_pet}, this amounts to the statement
$$L(1, f, \Ad) \sim_{E_f(\underline k)^\times} r_{\mathcal D}(m) \nu^{\Sigma_\infty \setminus J_1}(f),$$
up to appropriate powers of $\pi$.

Finally, Beilinson's conjecture implies that:
$$r_{\mathcal D}(\det H^{1}_{\mathcal M}(M, \Q(1))) = L(1, f, \Ad) \det H_B(M_\R, \Q)$$
as rational structures on $H^1_{\mathcal D}(M_\R, \R(1))$. Assuming this, Conjecture~\ref{conj:partial_weight_one} is equivalent to the statement:
$$\nu^{\Sigma_\infty \setminus J_1}(f) =  \det H_B(M_\R, \Q),$$
which we would expect to be true. It would be interesting to verify this final equality, but we decided to pursue this problem elsewhere.

\bibliographystyle{amsalphaurl}
\bibliography{bibliography}

\end{document}